\documentclass[reqno]{amsart}

\numberwithin{equation}{section}
\newtheorem{theorem}{\sc Theorem}[section]
\newtheorem{lemma}{\sc Lemma}[section]

\newtheorem{remark}{\sc Remark}

\def\os{\Omega_0^s}
\def\of{\Omega_0^f}

\def\R{\mathbb R}

\def\n{\nonumber}
\def\d{\displaystyle}
\def\w{\tilde w}

\def\a{\tilde a}

\def\q{\tilde q}
\def\e{\tilde\eta}
\def\z{\zeta}
\def\P{\Psi}
\def\w{\tilde w}

\def\ao{\alpha_1}
\def\at{\alpha_2}
\def\att{\alpha_3}

\title[Motion of quasilinear elastic solids in a viscous fluid]
      {On the interaction between quasilinear elastodynamics and the 
       Navier-Stokes equations}

\author[D. Coutand and S. Shkoller]{}

\email{coutand@math.ucdavis.edu}
\email{shkoller@math.ucdavis.edu}

\begin{document}

\maketitle

\centerline{\scshape   Daniel Coutand and Steve Shkoller}
%\medskip

  {\footnotesize \centerline{ Department of Mathematics }
  \centerline{ University of California at Davis } \centerline{Davis, CA
   95616 } }

\begin{abstract}
The interaction between a viscous fluid and an elastic solid is modeled by
a system of parabolic and hyperbolic equations, coupled to one another  
along the moving
material interface through the continuity of the velocity and traction vectors.
We prove the existence and uniqueness (locally in time) of strong
solutions in Sobolev spaces for  quasilinear elastodynamics coupled to the
incompressible Navier-Stokes equations.
Unlike our approach in \cite{CoSh2004} for the case of linear elastodynamics, 
we cannot employ a fixed-point argument on the nonlinear system itself, and 
are instead forced to regularize it by a {\it particular} parabolic artificial  
viscosity term. We proceed to show that with this specific regularization,  
we obtain a time interval of existence which is independent of the
artificial viscosity; together with a priori estimates, we identify the global 
solution (in both phases), as well as the interface motion, as a weak limit in 
srong norms of our sequence of regularized problems.
\end{abstract}

\section{Introduction}
We establish the existence and uniqueness in Sobolev spaces of strong solutions
to the unsteady fluid-structure interaction problem consisting of a nonlinear 
large-displacement elastic solid coupled to a viscous incompressible Newtonian 
fluid.  The fluid motion is governed by the incompressible Navier-Stokes equations,
while the solid, which can be either compressible or incompressible, is modeled by
the celebrated St. Venant-Kirchhoff constitutive law (although our method can be
applied to more general quasilinear hyperelastic models).

The first fluid-solid interaction problems solved were for the case of a
rigid body inside of a viscous flow in a bounded domain (see \cite{DeEs1999}, 
\cite{GrMa}, \cite{Conca}, \cite{Gunz}), and the case of a rigid body inside of
a viscous flow in an infinite domain (\cite{Wein}, \cite{Serre}, \cite{HiSe}). 
Later, the elastic body was modeled with the restriction of either a finite number
of modes (\cite{DeEsGrLe}) or with a hyperviscous type law for the solid 
(\cite{ChDeEsGr}, \cite{FlOr}), essentially by the same type of Eulerian global 
variational methods  developed in \cite{DeEs1999}. For the steady-state problem,
which is elliptic in both phases, \cite{Gr} solved the case of solid modeled as a
St. Venant-Kirchhoff material. In \cite{LW}, an Eulerian approach was used for
the case in which the solid is a visco-hyperelastic material, which is a 
regularization of a hyperbolic model of solid deformation.

With the exception of our recent well-posedness result for the case of a 
linear elastic solid in \cite{CoSh2004}, there are no known existence results 
for fluid-structure interaction  when the solid is modeled by a standard 
second-order hyperbolic equation.  This may be attributed to
the difficulties associated with coupling a parabolic PDE for the fluid with
a hyperbolic PDE for the solid through continuity of the velocity and traction
vectors across the moving material interface.  As we explained in 
\cite{CoSh2004}, an iteration scheme between fluid and solid phases fails to 
converge  due to a regularity loss induced by the hyperbolic phase (this 
divergent behavior has
been computationally noted as well in \cite{Heil1998}), and so we developed a 
method comprised of the following new ideas:
first, a functional framework which scales in a {\it hyperbolic} fashion for
{\it both}  the fluid and solid phases.  This scaling leads to
additional compatibility conditions in the fluid phase (when compared to the
use of the classical parabolic framework), and is absolutely crucial for
obtaining consistent energy estimates.  Second, we developed a regularity theory 
founded upon central {\it trace} estimates for the velocity vector restricted to
the interface, rather than traditional interior regularity arguments which do
not work for our problem.  Third, we were forced to bypass the use of the
frozen (or constant) coefficient basic linear problem, which requires
estimates on one more time derivative of the pressure function than the
initial data allows, and created a new method wherein the solution was
found as a limit of a sequence of penalized problems set in the Lagrangian 
framework.  The penalization scheme approximates the divergence-free constraint,
whereas the Lagrangian framework alleviates the difficulties associated with the
lack of a priori estimates in the solid phase for the frozen coefficient 
problem; this method indeed differs significantly from the classical methods used
in fluid-fluid interface problems (see for instance \cite{Sol1992}, \cite{Beale}). 

The fundamental difficulty in extending our result to the case of nonlinear elasticity is the absence of any method of analysis for
quasilinear elastodynamics which is compatible with the general scheme of 
\cite{CoSh2004}, involving a global Lagrangian variational formulation and the 
use of difference quotients to track the regularity of interface data. 
We remind the reader that unlike the analysis of elastostatic motion,
direct inverse function theorem arguments cannot be applied directly to the case 
of quasilinear elastodynamics due to the fact that the perturbation term arising
from the nonlinear operator is not an element of the appropriate function space 
for optimal regularity.  Alternatively, one might attempt a fixed-point approach,
wherein a portion of the nonlinear elasticity operator
is viewed as a forcing function coming from a given velocity $v$, and then try 
to solve a linear problem for an unknown $w$. The difficulty in this approach stems
from the fact that one has to find exact time derivatives of elastic energies 
for the forcing term associated to the elasticity operator, which is complicated by
the inner-product of a term involving $\int_0^t v$ and a term involving $w$.
This difficulty is overcome in \cite{DH}, by a clever and  essential use of the 
Dirichlet boundary condition in order to reformulate the problem in a non-standard
way.  As it turns out, the variety of known methods that have
been  used in the well-studied  area of quasilinear elasticity, such as those 
in \cite{DH} and \cite{HR} for the Dirichlet boundary condition, or in \cite{ST}
and \cite{ES} for the Neumann boundary condition, require a priori knowledge of
the boundary data regularity, and are hence instrinsically incompatible with 
fluid-structure interaction analysis (in fact, the methods devised for Dirichlet
conditions do not work for Neumann conditions and vice versa).
Indeed, of these various methodologies, only \cite{DH} and \cite{HR} use a variational
approach, the others employing either semi-group techniques as in
the early work of \cite{HKM} in the full space, or the technical 
paradifferential calculus in \cite{ST} for the two dimensional Neumann case. 

In this paper, we develop a new method for quasilinear elastodynamics,
variational in nature, that is compatible with our method in \cite{CoSh2004};
we proceed in two steps. First, we add a {\it specific} artificial viscosity to the solid
phase which regularizes the system, thus converting our hyperbolic PDE into
a parabolic one,  transforming the fluid-structure interaction into a
fluid-fluid interface-type problem for 
which existence and uniqueness of solutions is already known on a time interval
that a priori shrinks to zero as the artificial viscosity $\kappa$ tends to 
zero. Second, and this is where the primary difficulty rests, we prove that our 
{\it specific} choice of parabolic smoothing renders the time interval (on 
which a unique solution exists) independent of $\kappa$; furthermore, our a 
priori estimates allow us to construct a solution by weak convergence in strong
norms.  We note that the use of higher-order operators in the artificial
viscosity term, while providing the necessary a priori control of the
regularity of the moving interface,  would not yield $\kappa$-independent 
estimates which are essential here.  Also,
as our parabolic reguralization method is not specialized to any particular
boundary condition, it thus provides a {\it unified} approach to the classical problem
of quasilinear elastodynamics when the solid is not coupled to a fluid.

We now proceed to the formulation of our problem.
Let $\Omega \subset{\mathbb R}^3$ denote an open, 
bounded, connected and smooth domain with smooth boundary $\partial\Omega$ 
which 
represents the fluid container in which both the solid and fluid move.
Let $ \overline{\Omega ^s(t)}\subset  \Omega $ denote the closure of an open
and bounded
subset representing the solid body at each instant of time $t\in [0,T]$ with
$\Omega^f(t):= \Omega / \overline{\Omega^s(t)}$ denoting the fluid domain at 
each $t\in [0,T]$. 
Note that in our analysis $\Omega^s (t)$ is not necessarily connected,
which allows us to handle the case of several elastic bodies moving in the 
fluid.

\begin{remark}
If a function $u$ is defined on all of $\Omega$, we will denote 
$u^f=u\ 1_{\overline \Omega_0^f}$ and
$u^s=u\ 1_{\overline\Omega_0^s}$. This allows us to indicate from which phase
the traces on 
$$\Gamma(0):= \overline{\Omega^f(0)}\cap \overline{\Omega^s(0)} $$ 
of various discontinuous terms arise, and also to
specify functions that are associated with the fluid and solid phases.
\end{remark}

For each $t\in (0,T]$, we wish to find the location of these domains
inside $\Omega$, the
divergence-free velocity field $u^f(t, \cdot)$ of the fluid, the fluid pressure
function $p(t, \cdot)$ on $\Omega^f(t)$,  the fluid Lagrangian 
volume-preserving configuration
$ \eta ^f(t, \cdot ):\Omega ^f(0)=\Omega_0^f \rightarrow \Omega ^f(t)$, 
and the elastic  Lagrangian configuration field 
$\eta^s(t, \cdot):\Omega^s(0)=\Omega_0^s \rightarrow \Omega ^s(t)$ such that
$
\Omega=  \eta^s(t,\overline{\Omega^s_0}) \cup \eta^f(t,\Omega^f_0) $,
where
$\eta^f_t(t,x) = u^f(t, \eta^f(t,x))$ and $u^f$ solves the Navier-Stokes equations in $\Omega^f(t)$:
\begin{align*}
u^f_t + (u^f\cdot \nabla)u^f  &=  \operatorname{div} T^f + f_f \,,
        \\
   \operatorname{div} u^f &= 0     \,,
\end{align*}
with 
\begin{equation} 
\label{fluidlaw}
T^f=\nu\ \operatorname{Def} u^f -p\ \text{I} \,,
\end{equation}
and $\eta^s$ solves the elasticity equations on $\Omega^s(0)$
$$\ddot{\eta}^s  = \operatorname{div} T^s + f_s, $$
         with $T^s= \frac{\lambda}{2}\ \text{Tr}((\nabla\eta^s)^T\nabla\eta^s-I)\ I\ +\mu\ ((\nabla\eta^s)^T\nabla\eta^s-\ I) \,,$
and where the equations are coupled together by the continuity of
the normal component of stress along the material interface 
$\Gamma(t):= \overline{\Omega^s}(t) \cap \overline{\Omega^f}(t)$ expressed 
in the Lagrangian representation on
$\Gamma_0:=\Gamma(0)$
as
$$T^s \ N
=  
[T^f \circ \eta^f] \ [( \nabla \eta^f)^{-1} \ N ]
\,,$$
and the continuity of particle displacement fields along
$\Gamma_0$
$$       
\eta^f =  \eta^s
\,,$$
together with the initial conditions $u(0,x) = u_0(x)$, and  $\eta(0,x) = x$ 
and the Dirichlet (no-slip) condition on the boundary 
$\partial\Omega$ of the container $u^f=0$,
where $\nu>0$ is the kinematic viscosity of the fluid, $\lambda>0$ and $\mu>0$ 
denote the Lam\'e constants of the elastic material, $N$ is the outward unit 
normal to $\Gamma_0$ and $\operatorname{Def}u$ is twice the rate of deformation
tensor of $u$, given in coordinates by $u^i,_j + u^j,_i$.  All Latin indices 
run through $1,2,3$, the Einstein  summation convention is employed, and indices
after commas denote partial derivatives.

We now briefly outline the proof. As the solid and fluid phases are naturally 
expressed in the Lagrangian and Eulerian framework, respectively, we begin by transforming the fluid phase into Lagrangian 
coordinates, leading us to the system of equations (\ref{nsl}) and, as in \cite{CoSh2004}, we work in an hyperbolic framework in order to accomodate the dual nature of the problem (parabolic in 
the fluid and hyperbolic in the solid).  

In order to solve (\ref{nsl}), in Section \ref{7}, we first add a {\it particular} form of 
artificial viscosity to the quasilinear hyperbolic equation in the solid, transforming the hyperbolic phase into a parabolic one;  
specifically, we add the term $-\kappa L(\eta_t)$, where $L$ denotes 
the linearized (about the identity) elasticity operator and $\eta_t$ is
the material velocity. We hence obtain an interface problem that is
{\it parabolic} in nature in both phases, and can be thought as a fluid-fluid 
parabolic interface problem for which well-posedness is classical, save for the fact that both phases are required to scale in an hyperbolic fashion. 
The time interval of existence $[0,T_\kappa]$ for this parabolic system
a priori {\it shrinks} to zero as $\kappa\rightarrow 0$.

In Section \ref{8}, we establish $\kappa$-independent estimates on the
solutions $v_\kappa$ to the regularized parabolic problem on the time
interval $[0,T_\kappa]$ by identifying exact time derivatives of elastic
energies, and establish regularity of the 
{\it interface}. A direct fixed-point approach for (\ref{nsl}) does not appear to
yield these  exact time derivatives for the elastic energy, 
whereas the the regularized problem (\ref{smoothie}) does indeed lead to them. 
An essential {\it key} for obtaining
estimates independent of $\kappa$ inside the solid is Lemma \ref{key}. Whereas the trace estimates could be carried with other choices of artificial viscosity, we absolutely need the special choice made in our analysis in order to
recover the regularity inside the solid independently of $\kappa$. 
In particular, a different choice of a regularizing operator
either of the same order such as $-\triangle\eta_t$ or of higher order
such as $L^2\eta$ or $L^2\eta_t$ would not provide $\kappa$-independent estimates.

In Section \ref{9}, we then explain how our estimates allow the construction 
of solutions $v_{\kappa}$ on a time interval independent of $\kappa$, still 
with energy estimates independent of $\kappa$. The existence of a solution to
(\ref{nsl}) then follows by weak convergence as $\kappa\rightarrow 0$.

Uniqueness is established in Section \ref{11} in the same functional framework 
used for existence.

As our method seemingly requires more regularity on the initial data in the solid 
than it should, due to the artificial 
viscosity in the compatibility conditions, we explain in Section \ref{12} how 
this extra regularity can be removed, thus 
leading to the result with optimal regularity. 

Section \ref{13} is dedicated to the case where the incompressibility constraint 
is added to the solid. The additional difficulty with respect to the compressible 
case comes from the fact that we control the velocity uniformly in $\kappa$ in 
function spaces which possess less regularity than in the fluid, whereas the 
pressure is controlled uniformly in the same regularity spaces in both phases. 
Also, we cannot use Lemma \ref{key} in the  most optimal form for the 
regularity of the pressure in the solid phase.

\section{Notational simplification}
\label{2}
Although a fluid with a Neumann (free-slip) boundary condition indeed obeys the 
constitutive law (\ref{fluidlaw}), we will replace for notational convenience (\ref{fluidlaw}) with
\begin{equation}
\label{newfluid}
T^f=\nu \nabla u^f-p \text{I};
\end{equation}
this amounts to replacing the energy 
$\displaystyle\int_{\Omega_0^f}\operatorname{Def}u^f : \operatorname{Def} v$ 
by $\displaystyle\int_{\Omega_0^f}\nabla u^f:\nabla v$, 
which is not a problem mathematically due to the
well-known Korn inequality. Henceforth, we shall take (\ref{newfluid}) as the
fluid constitutive law.

\section{Lagrangian formulation of the problem}
\label{3}
In regards to the forcing functions, we shall use the convention of denoting 
both the fluid forcing $f_f$ and the solid forcing $f_s$ by the same letter $f$. 
Since $f_f$ has to be defined in 
$\Omega$ (because of the composition with $\eta$), and $f_s$ must be defined 
in $\Omega_0^s$, we will assume that the forcing $f$ is defined over the entire
domain $\Omega$.
 
Let 
\begin{equation}\label{a}
a(x) = [\text{Cof} \nabla\eta^f(x)]^T, 
\end{equation}
where $(\nabla \eta^f(x))^i_j = \partial (\eta^f)^i/\partial x^j(x)$
denotes the matrix of partial derivatives of $\eta^f$.  
Clearly, the matrix $a$ depends on $\eta$ and we shall sometimes use the 
notation $a^i_j(\eta)$ to denote the formula (\ref{a}).

Let $v=u\circ \eta$ denote
the Lagrangian or material velocity field, $q=p \circ \eta$ is the Lagrangian 
pressure function (in the fluid), 
and $F= f^f \circ \eta^f$ is the fluid forcing function in
the material frame.  
Then, as long as no collisions occur between the solids (if there are initially
more than one) or between a solid and $\partial\Omega$, the problem can be reformulated as 
\begin{subequations}
  \label{nsl}
\begin{alignat}{2}
\eta_t &=v&\ &\text{in} \ \ (0,T)\times \Omega \,, 
         \label{nsl.a}\\
v^i_t - \nu (a^j_l a^k_l v^i,_k),_j + (a^k_i q),_k &= F^i 
&&\text{in} \ \ (0,T)\times \Omega_0^f \,, 
         \label{nsl.b}\\
   a^k_i v^i,_k &= 0     &&\text{in} \ \ (0,T)\times \Omega_0^f \,, 
         \label{nsl.c}\\
v_t - c^{mjkl} [(\eta,_m\cdot\eta,_j-\delta_{mj}) \eta,_k],_l &= f
&&\text{in} \ \ (0,T)\times \Omega_0^s \,, 
         \label{nsl.d}\\
c^{mjkl}\ (\eta,_m\cdot\eta,_j-\delta_{mj}) \eta^i,_k\ N_l&=
\nu\ v^i,_k a^k_l a^j_l N_j - q a^j_i N_j  
&&\text{on} \ \ (0,T)\times \Gamma_0 \,, 
         \label{nsl.e}\\
v(t,\cdot)& \in H^1_0 (\Omega;{\mathbb R}^3) &&\text{a.e. in } \ \ (0,T)\,, \label{nsl.f}\\
   v &= u_0  
 &&\text{on} \ \ \Omega\times \{ t=0\} \,, 
         \label{nsl.g}\\
   \eta &= \text{Id}     
 &&\text{on} \ \ \Omega\times \{ t=0\} \,, 
         \label{nsl.h}
\end{alignat}
\end{subequations}
where $N$ denotes the outward-pointing unit normal to $\Gamma_0$ (pointing
into the solid phase), and
$$c^{ijkl}=\lambda \delta^{ij}\delta^{kl}+\mu (\delta^{ik}\delta^{jl}
+\delta^{il}\delta^{jk})\ . $$
Throughout the paper, all Greek indices run through $1,2$ and 
all Latin indices run through $1,2,3$.  Note that the continuity of the
velocity along the interface is satisfied in the sense of
traces on $\Gamma_0$ by condition (\ref{nsl.f}), whereas the continuity of 
the normal stress along the interface is represented by (\ref{nsl.e}).

%\begin{remark}
%The case in which  the viscosity or Lam\'e coefficients are variable functions
%depending on $x \in \Omega$ and satisfying the usual 
%assumptions, can be handled by our methodology without any supplementary 
%mathematical difficulties.
%\end{remark}

%\begin{remark}
%Remember that we use (\ref{newfluid}) as the constitutive law for the fluid.
%\end{remark} 

\section{Notation and conventions}
\label{4}
We begin by specifying our notation for certain vector and matrix operations.
\begin{itemize}
\item[] We write the Euclidean inner-product between two vectors $x$ and $y$ 
as $x\cdot y$, so that $x\cdot y=x^i\ y^i$.
\item[] The transpose of a matrix $A$ will be denoted by $A^T$, {\it i.e.}, 
$(A^T)^i_j=A^j_i$.
\item[] We write the product of a matrix $A$ and a vector $b$ as $A\ b$, 
{\it i.e}, $(A\ b)^i=A^i_j b^j$.
\item[] The product of two matrices $A$ and $S$ will be denoted by 
$A\cdot S$, {\it i.e.}, $(A\cdot S)^i_j=A^i_k\ S^k_j$.
\end{itemize}

For $T>0$ and $k\in {\mathbb N}$, we set
\begin{align*}
V^k_f(T)&=\{ w \in L^2(0,T;  H^k(\Omega_0^f;{\mathbb R}^3)) \ | \ 
\partial_t^n w \in L^2(0,T;  H^{k-n}(\Omega_0^f;{\mathbb R}^3)),\ n=1,...,k\},  
\end{align*}
with
$V^k_s(T)$ defined with $\os$ replacing $\of$.

In order to specify the initial data for the weak formulation, we 
introduce the space
$$L^2_{div,f}
=\{ \psi \in L^2(\Omega; {\mathbb R}^3) \ | \ \operatorname{div}\psi=0 \ 
\text{in}\ \Omega_0^f, \ \  \psi\cdot N=0\ \text{on}\ \partial \Omega \}\ ,$$
which is endowed with the $L^2(\Omega; {\mathbb R}^3)$ scalar product.

The space of velocities, $X_T$, where the solution to (\ref{nsl}) lives, is defined as the following separable Hilbert space:
\begin{align}
X_T= \{ v\in L^2(0,T; H^1_0(\Omega;{\mathbb R}^3))\ 
|\  \ (v^f, \int_0^\cdot v^s)\in V^4_f (T)
\times V^4_s (T)\}\ ,
\label{XT}
\end{align}
 endowed with its natural Hilbert norm
\begin{align*}
\|v\|^2_{X_T}&=\|  v \|^2_ {L^2(0,T;H^1_0(\Omega;{\mathbb R}^3))}\\
&\ \ +\sum_{n=0}^3 [\|\partial_t^n  v \|^2_ {L^2(0,T;H^{4-n}(\Omega_0^f;{\mathbb R}^3))}+ \|\partial_t^n  \int_0^\cdot v \|^2_ {L^2(0,T;H^{4-n}(\Omega_0^s;{\mathbb R}^3))}]
\,.
\end{align*}
We also need  the space
\begin{align*}
Y_T= \{ (v,q)\in X_T\times L^2(0,T;H^3(\Omega_0^f;{\mathbb R}))|\  \partial_t^n q\in L^2(0,T;H^{3-n}(\Omega_0^f;{\mathbb R})) (n=1,2)\} , 
\end{align*}
 endowed with its natural Hilbert norm
\begin{align*}
\|(v,q)\|^2_{Y_T}=&\ \|  v \|^2_ {X_T}
+ \sum_{n=0}^2\|  \partial_t^n q \|^2_ {L^2(0,T;H^{3-n}(\Omega_0^f;{\mathbb R}))}
\,.
\end{align*}

%\begin{remark}
%Note well that our method does not require any a priori knowledge of the 
%regularity of the third time derivative of the pressure function $q_{ttt}$; 
%this is due to the Dirichlet boundary condition on 
%$\partial\Omega$ as well as the Lagrangian representation of the problem
%that we employ.
%\end{remark}

We shall also need $L^{\infty}$-in-time control of certain norms of the
velocity, which necessitates the use of the  following closed subspace of $X_T$:
\begin{align*}
W_T=\{\ v\in X_T|\ &\ v_{ttt}\in L^{\infty}(0,T;L^2(\Omega;{\mathbb R}^3)),\\
& \partial_t^n \int_0^\cdot v\in L^{\infty}(0,T;H^{4-n}(\Omega_0^s;{\mathbb R}^3))( n=0,1,2,3)\}\ ,
\end{align*}
endowed with the following norm
\begin{align*}
\|v\|^2_{W_T}&=\|  v \|^2_ {X_T}
+\|  v_{ttt} \|^2_ {L^{\infty}(0,T;L^2(\Omega;{\mathbb R}^3))}+ 
\sum_{n=0}^3 \|\partial_t^n \int_0^{\cdot} v\|^2_ {L^{\infty}(0,T;H^{4-n}(\Omega_0^s;{\mathbb R}^3))}
\,.
\end{align*}
Finally, we will also make use of the space
\begin{align*}
Z_T= \{ (v,q)\in W_T\times L^2(0,T;H^3(\Omega_0^f;{\mathbb R}))|&\ \partial_t^n q\in L^2(0,T;H^{3-n} (\Omega_0^f;{\mathbb R})), (n=1,2)\\
&|\ q_{tt}\in L^\infty(0,T;L^2(\of;\R))\}\ , 
\end{align*}
endowed with its natural norm
\begin{align*}
\|(v,q)\|^2_{Z_T}=&\ \|  v \|^2_ {W_T}
+ \sum_{n=0}^2 \| \partial_t^n q \|^2_ {L^2(0,T;H^{3-n}(\Omega_0^f;{\mathbb R}))}+  \| \partial_t^2 q \|^2_ {L^\infty(0,T;L^2(\Omega_0^f;{\mathbb R}))}
\,.
\end{align*}

\begin{remark}
Note that our functional framework does not make use of the third time
derivative of the pressure $q_{ttt}$, even though we do use the third time
derivative of velocity $w_{ttt}$; this functional framework is necessitated by
the fact that the Dirichlet boundary condition together with the limited 
regularity of $w_{ttt}$ does not allow us to obtain $q_{ttt}$ with the appropriate
regularity. 
Note also that we have added the $L^\infty$-in-time control of $q_{tt}$ in the 
definition of $Z_T$ mostly for a more convenient way to prove our theorems, 
rather than absolute necessity.
\end{remark}

Throughout the paper, we shall use
$C$ to denote a generic constant, which may possibly depend on the coefficients
$\nu$, $\lambda$, $\mu$,  or on the initial geometry given by $\Omega$ and
$\Omega_0^f$ (such as a Sobolev constant or an elliptic constant).
For the sake of notational convenience, we will also write
$u(t)$ for $u(t,\cdot)$. 

\section{A first theorem} 
\label{5}

We now state our first theorem. We impose greater regularity requirements on
the initial data than is optimal so as to avoid technical difficulties associated
with a particular type of initial data regularization that would otherwise be
necessitated.  We consider the case of optimal regularity on the initial data in 
Theorem \ref{optimal}.

\begin{theorem}\label{main}
Let $\Omega\subset {\mathbb R}^3$ be a bounded domain of class $H^4$, and let 
$\Omega_0^s$ be an open set (with a finite number $\ge 1$ of 
connected components) of class $H^4$ such that 
$\overline{\Omega_0^s}\subset \Omega$ and such that the distance between two 
distinct connected components of $\os$ (if there are multiple solid components) 
is greater than zero.
Let us denote $\Omega_0^f=\Omega\cap (\overline{\Omega_0^s})^c$.
Let $\nu>0$, $\lambda> 0$, $\mu>0$ be given. 
Let 
\begin{subequations}
\label{f_regularity}
\begin{align} 
&(f, f_t, f_{tt}, f_{ttt}) \in L^2(0,\bar T; H^3(\Omega;{\mathbb R}^3)\times H^2(\Omega;{\mathbb R}^3)\times H^1(\Omega;{\mathbb R}^3)\times L^2(\Omega;{\mathbb R}^3)),\\
& f(0)\in H^4(\Omega;{\mathbb R}^3)\ ,\ f_t(0)\in H^4(\Omega;{\mathbb R}^3).
\end{align}
\end{subequations}
Assume that the initial data  satisfies
$$u_0 \in H^6(\Omega_0^f;{\mathbb R}^3)\cap H^6(\Omega_0^s;{\mathbb R}^3) 
\cap H^1_0(\Omega;{\mathbb R}^3)\cap L^2_{{div},f},$$
as well as  the compatibility conditions 
\begin{subequations}
\label{compatibility}
\begin{align}
& [\nabla u_0^f\  N]_{\operatorname{tan}}=0\ \text{ on }\ \Gamma_0,\ \  w_1=0=w_2\ \text{on}\ \partial\Omega,\ \ \nu\triangle u_0^f-\nabla q_0=0\ \text{on}\ \Gamma_0,\label{c1}\\
&\  [ (\nu [(a_l^k a_l^j) {w^f},_k^i]_t(0) N_j)_{i=1}^3]_{\operatorname{tan}}  -[(q_0\ {a_i^j}_t(0) N_j)_{i=1}^3]_{\operatorname{tan}}\n\\
&\qquad\qquad\qquad\qquad=[c^{mjkl} \bigl[(\eta^s,_m\cdot\eta^s,_j-\delta_{mj})\eta^s,_k\bigr]_t(0) N_l]_{\operatorname{tan}}\ \text{on}\ \Gamma_0,\label{c2}\\
&\  [ (\nu [(a_l^k a_l^j) {w^f},_k^i]_{tt}(0) N_j)_{i=1}^3]_{\operatorname{tan}}
 -[(2 q_1\ {a_i^j}_t(0) N_j +q_0\ {a_i^j}_{tt}(0) N_j)_{i=1}^3]_{\operatorname{tan}}\n\\
&\qquad\qquad\qquad\qquad=[c^{mjkl} \bigl[(\eta^s,_m\cdot\eta^s,_j-\delta_{mj})\eta^s,_k\bigr]_{tt}(0) N_l]_{\operatorname{tan}}\ \text{on}\ \Gamma_0,\label{c3}\\
&\ \nu \triangle  {w_1^f} +\nu (( a_l^j a_l^k)_t (0) {u_0^f},_k),_j + F_t(0) - [((a_i^j)_t(0) q_0),_j + q_1,_i]_{i=1}^3 \n\\
 &\qquad\qquad\qquad\qquad  =f_t(0) + c^{mjkl} \bigl[[(\eta^s,_m\cdot\eta^s,_j-\delta_{mj})\eta^s,_k],_l\bigr]_t(0)\ \text{on}\ \Gamma_0 \ ,\label{c4}
\end{align}
\end{subequations}
where the time derivatives appearing in these equations and in the following ones 
are computed from any $w$ satisfying $w(0)=u_0$, $\partial_t^n w (0)=w_n$ 
($n=1,2$), and from any $q$ satisfying $\partial_t^n q (0)=q_n$ ($n=0,1,2$), the 
quantities $w_n$ and $q_n$ being defined as follows. First, 
 $q_0\in H^3(\Omega_0^f;{\mathbb R})$ is defined by 
\begin{subequations}
\label{defq0}
\begin{align}
\triangle q_0&=\operatorname{div} f(0)+(a_i^j)_t (0) u_0^i,_j\ \text{in}\ \Omega_0^f,\\
q_0&=\nu [\nabla u_0^f\ N]\cdot N\ \text{on}\ \Gamma_0,\label{defq0.b}\\
\frac{\partial q_0}{\partial N} &=f(0)\cdot N+\nu \triangle u_0^f\cdot N\ \text{on}\ \partial\Omega ,
\end{align}
\end{subequations}
and $w_1\in H^1_0(\Omega;{\mathbb R}^3)\cap H^4(\Omega_0^s;{\mathbb R}^3)\cap H^4(\Omega_0^f;{\mathbb R}^3)$ by 
\begin{subequations}
\label{def1}
\begin{align}
 w_1&=\nu\triangle u_0-\nabla q_0 +f(0)\ \text{in}\ \Omega_0^f\\
w_1&=f(0)\ \text{in}\ \Omega_0^s\ .
\end{align}
\end{subequations}
Note that $w_1\in H^4(\Omega_0^f;{\mathbb R}^3)$ since $\triangle w_1\in H^2(\of;\R^3)$ and $w_1=0$ on $\partial\Omega$, $w_1=f(0)$ on $\Gamma_0$. We also have 
$q_1\in H^3(\Omega_0^f;{\mathbb R})$ defined by
\begin{subequations}
\label{defq1}
\begin{align}
\triangle q_1= &\ \operatorname{div}[\nu \triangle w_1 + F_t(0)+  [\nu(( a_l^j a_l^k)_t (0) u_0^i,_k),_j - ((a_i^j)_t(0) q_0),_j]_{i=1}^3]\nonumber\\
& + 2 (a_i^j)_t (0) w_1^i,_j+(a_i^j)_{tt} (0) u_0^i,_j
\ \text{in}\ \Omega_0^f\ ,\\
q_1 =&\  \nu [\nabla w_1^f\ N\cdot N + (a_l^k a_l^j)_t(0) {u_0^f}^i,_k N_j N_i]
-q_0\ {a_i^j}_t(0) N_j N_i\n\\
& -c^{mjkl} \bigl[(\eta^s,_m\cdot\eta^s,_j-\delta_{mj})\eta^s,_k\bigr]_t(0) N_l\cdot N\ \text{on}\ \Gamma_0,\\
\frac{\partial q_1}{\partial N}=& F_t(0)\cdot N - [(a_i^j)_t(0) q_0],_j N_i +\nu \triangle w_1\cdot N+\nu ((a_l^j a_l^k)_t(0) u_0^i,_k),_j N_i\ \text{on}\ \partial\Omega\ , 
\end{align}
\end{subequations}
and $w_2\in H^1_0(\Omega;{\mathbb R}^3)\cap H^4(\Omega_0^s;{\mathbb R}^3)\cap H^2(\Omega_0^f;{\mathbb R}^3)$ by
\begin{subequations}
\label{w2def}
\begin{align}
w_2^i= &\ \nu \triangle  w_1^i +\nu (( a_l^j a_l^k)_t (0) u_0^i,_k),_j + F_t^i(0) - ((a_i^j)_t(0) q_0),_j - q_1,_i
\ \text{in}\ \Omega_0^f\ ,\label{w2def.a}\\
{w}_2 =&\  f_t(0) + c^{mjkl} \bigl[[(\eta,_m\cdot\eta,_j-\delta_{mj})\eta,_k],_l\bigr]_t(0)\ \text{in}\ \Omega_0^s\ ,\label{w2def.b}
\end{align}
\end{subequations}
 Finally, $q_2\in H^1(\Omega_0^f;{\mathbb R})$ is defined by
\begin{subequations}
\label{defq2}
\begin{align}
\triangle q_2= &\ \operatorname{div}[(f\circ\eta)_{tt}(0)+  \nu\ [( a_l^j a_l^k w,_k,_j]_{tt}(0) - [((a_i^j)_{tt}(0) q_0+ 2 (a_i^j)_{t}(0) q_1),_j]_{i=1}^3]\nonumber\\
& + 3 (a_i^j)_{t} (0) w_2^i,_j+ 3 (a_i^j)_{tt} (0) w_1^i,_j+(a_i^j)_{ttt} (0) u_0^i,_j
\ \text{in}\ \Omega_0^f\ ,\\
q_2 =&\  \nu [(a_l^k a_l^j) {w^f}^i,_k]_{tt}(0) N_j N_i -c^{mjkl} \bigl[(\eta^s,_m\cdot\eta^s,_j-\delta_{mj})\eta^s,_k\bigr]_{tt}(0) N_l \cdot N\n\\
& -q_0\ {a_i^j}_{tt}(0) N_j N_i-2 q_1\ {a_i^j}_{t}(0) N_j N_i\ \text{on}\ \Gamma_0,\\
\frac{\partial q_2}{\partial N}=& (f\circ\eta)_{tt}(0)\cdot N - 2[(a_i^j)_t(0) q_1],_j N_i - [(a_i^j)_{tt}(0) q_0],_j N_i+\nu \triangle w_2\cdot N\n\\
&\ + 2 \nu ((a_l^j a_l^k)_t(0) w_1^i,_k),_j N_i
 + \nu ((a_l^j a_l^k)_{tt}(0) u_0^i,_k),_j N_i\ \ \text{on}\ \partial\Omega\ . 
\end{align}
\end{subequations}

Then there 
exists $T\in (0,\bar T)$ depending on $u_0$, $f$, and $\Omega_0^f$, such that 
there exists a unique solution $(v,q) \in Z_T$  
of the problem (\ref{nsl}).  Furthermore, 
$\eta \in  C^0([0, T]; H^4(\Omega_0^f; {\mathbb R}^3)
\cap  H^4(\Omega_0^s; {\mathbb R}^3)\cap H^1(\Omega; {\mathbb R}^3))$.
\end{theorem}

\begin{remark}
The remarks appearing in \cite{CoSh2004} at the end of Section 5 concerning the 
compatibility conditions and forcing functions for the linear elasticity case 
still hold in this setting with the necessary adjustments. In particular, we do 
not need the forcing functions to have the same regularity in both phases.
\end{remark}

\section{ Preliminary result}
\label{6}
In the remainder of the paper, we set
$$L(u)^i=[c^{ijkl} (u^k,_l+u^l,_k)],_j.$$

 In our limit process as the artificial viscosity tends to zero, we will make use in a crucial way of the basic following result:
\begin{lemma}
\label{key}
Let $g\in C^0([0,T];L^2(\os;\R^3))$ and $u$ be such that $u_t\in L^2 (0,T;H^2(\os;\R^3))$ and
\begin{equation}
\label{equality}
\epsilon L(u_t)+L(u)=g\ \text{on}\ [0,T]\times\os.
\end{equation}
Then, independently of $\epsilon>0$, $$\|L(u)\|_{L^{\infty}(0,T;L^2(\os;\R^3))}\le   \|g\|_{L^{\infty}(0,T;L^2(\os;\R^3))}+\|L(u_0)\|_{L^2(\os;\R^3)}.$$
\end{lemma}

\begin{proof}
Since $L(u)\in C^0(0,T;L^2(\os;\R^3))$, let $T'\in [0,T]$ be such that
$$\|L(u(T'))\|_{L^2(\os;\R^3)}=\sup_{[0,T]} \|L(u)\|_{L^2(\os;\R^3)}.$$
If $T'=0$, then the statement of the Lemma is satisfied. Now, let
us assume that $T'\in (0,T]$. Let $\delta\in (0,T')$ be arbitrary. From
(\ref{equality}), we infer that
\begin{align*}
\epsilon^2\int_{T'-\delta}^{T'} &\|L(u_t)\|^2_{L^2(\os;\R^3)}+ \int_{T'-\delta}^{T'} \|L(u)\|^2_{L^2(\os;\R^3)}\\
&+\epsilon\ [ \|L(u(T'))\|^2_{L^2(\os;\R^3)}-
\|L(u(T'-\delta))\|^2_{L^2(\os;\R^3)}]=\int_{T'-\delta}^{T'} \|g\|^2_{L^2(\os;\R^3)}.
\end{align*}
From the definition of $T'$ we then infer that for any $\delta\in (0,T')$,
$$ \int_{T'-\delta}^{T'} \|L(u)\|^2_{L^2(\os;\R^3)}\le \int_{T'-\delta}^{T'} \|g\|^2_{L^2(\os;\R^3)},$$ which after division by $\delta$ gives at the limit
$\delta\rightarrow 0$:
$$ \|L(u(T'))\|^2_{L^2(\os;\R^3)}\le \|g(T')\|^2_{L^2(\os;\R^3)},$$
which concludes the proof of the Lemma.
\end{proof}

\begin{remark}
It should be clear that
Lemma \ref{key} applies to a more general class of linear operators than $L$.
\end{remark} 

\section{The smoothed problem and its basic linear problem}
\label{7}
As we described in the introduction, we cannot find an appropriate linear problem 
whose fixed-point provides a solution to (\ref{nsl}).  We are thus lead to 
introduce introduce the following (parabolic) regularization of (\ref{nsl}), 
with the artificial viscosity $\kappa>0$:
\begin{subequations}
  \label{smoothie}
\begin{alignat}{2}
 v^i_t - \nu (a^j_l a^k_l  v^i,_k),_j + (a^k_i q),_k &=  f^i\circ{\eta} 
&&\text{in} \ \ (0,T)\times \Omega_0^f \,, 
         \\
   a^k_i  v^i,_k &= 0    \ &&\text{in} \ \ (0,T)\times \Omega_0^f \,, 
         \\
 v^i_t - \kappa [c^{ijkl} v^k,_l],_j +N(\eta)^i &=  f^i + \kappa h^i
\ &&\text{in} \ \ (0,T)\times \Omega_0^s \,, 
         \\
 \kappa\ c^{ijkl} v^k,_l\ N_j
+ G(\eta)^i&=\nu\ v^i,_k a^k_l a^j_l N_j- q a^j_i N_j+ \kappa g^i
&&\text{on} \ \ (0,T)\times \Gamma_0 \,, 
         \\
v(t,\cdot) &\in H^1_0 (\Omega;{\mathbb R}^3) \ &&\text{a.e. in } \ \ (0,T)\,, \\
   v &= u_0  
\ &&\text{on} \ \ \Omega\times \{ t=0\} \,, 
\end{alignat}
\end{subequations}
where 
\begin{subequations}
\begin{align}
\label{forcetrbis}
N(\eta)&= - c^{mjkl} [(\eta,_m\cdot\eta,_j-\delta_{mj}) \eta^i,_k],_l\ \text{in}\ \os,\\
G(\eta)&=   c^{mjkl} [(\eta,_m\cdot\eta,_j-\delta_{mj}) \eta^i,_k] N_l \ \text{on}\ \Gamma_0\ ,
\end{align}
\end{subequations}
and
\begin{subequations}
\begin{align}
\label{forcetr}
h^i(t,\cdot)&= -  [c^{ijkl} (u_0+t w_1+\frac{t^2}{2} w_2)^k,_l],_j\ \ \text{in}\ \os,\\
g^i(t,\cdot)&=   [c^{ijkl} (u_0+t w_1+\frac{t^2}{2} w_2)^k,_l]\ N_j\ \ \text{on}\ \Gamma_0\ .
\end{align}
\end{subequations}
Solutions to (\ref{nsl}) will be obtained as the limit as $\kappa \rightarrow 0$ of
solutions to (\ref{smoothie}).

Suppose that $ v\in W_T$ is given. 
Let $\displaystyle\eta=Id+\int_0^{\cdot}v$ and let $a_i^j$ be the quantity
associated with $\eta$ through (\ref{a}).

We are concerned with the following time-dependent linear problem, whose 
fixed-point $w=v$ provides a solution to (\ref{smoothie}):
\begin{subequations}
  \label{linear}
\begin{alignat}{2}
 w^i_t - \nu (a^j_l a^k_l  w^i,_k),_j + (a^k_i q),_k &=  f^i\circ{\eta} 
&&\text{in} \ \ (0,T)\times \Omega_0^f \,, 
         \label{linear.b}\\
   a^k_i  w^i,_k &= 0    &&\text{in} \ \ (0,T)\times \Omega_0^f \,, 
         \label{linear.c}\\
 w^i_t - \kappa [c^{ijkl} w^k,_l],_j +N(\eta)^i&=  f^i + \kappa h^i
&&\text{in} \ \ (0,T)\times \Omega_0^s \,, 
         \label{linear.d}\\
 \kappa\ c^{ijkl} w^k,_l\ N_j+ G(\eta)^i&=\nu\ w^i,_k a^k_l a^j_l N_j\n\\
&\ \ \ - q a^j_i N_j+ \kappa g^i
&&\text{on} \ \ (0,T)\times \Gamma_0 \,, 
         \label{linear.e}\\
w(t,\cdot)& \in H^1_0 (\Omega;{\mathbb R}^3) &&\text{a.e. in } \ \ (0,T)\,, \label{linear.f}\\
   w &= u_0  
 &&\text{on} \ \ \Omega\times \{ t=0\} . 
         \label{linear.g}
\end{alignat}
\end{subequations}

\begin{remark}
The two forcing functions (\ref{forcetr}) are introduced for compatibility 
conditions at time $t=0$, allowing the solution of (\ref{linear}) to satisfy 
$(w_t(0),w_{tt}(0))\in H^1_0(\Omega;\R^3)^2$ and even to satisfy the same initial 
conditions as solutions to (\ref{nsl}) would. 
\end{remark}

In the following, for the sake of notational convenience, we will denote by $N(u_0,(w_i)_{i=1}^3)$ a generic smooth
function depending only on $\d\sum_{i=0}^3 [\ \|w_{3-i}\|_{H^i(\os;\R^3)}+
\|w_{3-i}\|_{H^i(\of;\R^3)} \ ]$ (with the convention that $u_0=w_0$), by $N((q_i)_{i=0}^2)$ a generic smooth
function depending only on $\d\sum_{i=0}^2 \ \|q_{2-i}\|_{H^i(\of;\R^3)}$ and by $ M(f,\kappa g,\kappa h)$ a generic smooth function
depending only on $\|f\|_{V^3_f(T)}+\|f\|_{V^3_s(T)}+ \kappa\ [\ \|u_0\|_{H^4(\os;\R^3)}+\|w_1\|_{H^4(\os;\R^3)}+\|w_2\|_{H^4(\os;\R^3)}\ ]$. 
Then, let $w_3\in L^2(\Omega;{\mathbb R}^3)$ be defined by
\begin{subequations}
\label{w3def}
\begin{align}
w_3^i= &\ \nu [( a_l^j a_l^k  w^i,_k),_j]_{tt}(0) + (f\circ\eta)_{tt}(0)  - [a_i^j q,_j]_{tt}(0) \ \text{in}\ \Omega_0^f\ ,\label{w3def.a}\\
{w}_3^i =&\  f_{tt}^i(0) + c^{mjkl} \bigl[[(\eta,_m\cdot\eta,_j-\delta_{mj})\eta,_k],_l\bigr]_{tt}(0)\ \text{in}\ \Omega_0^s\ ,\label{w3def.b}
\end{align}
\end{subequations}
where the time derivatives are computed with any $\eta(0,x)=x$, $w=\eta_t (0)=u_0$,
$\partial_t^n w(0)=w_n, (n=1,2)$,  $\partial_t^n q(0)=q_n, (n=0,1,2)$. 

%Once again, we remind that $(a_i^j)_t(0)$, $(a_i^j)_{tt}(0)$ and $(a_i^j)_{ttt}(0)$ are only functions of $u_0^f$, $w_1^f$ and $w_2^f$, and are also equal to respectively $(\tilde a_i^j)_t(0)$, $(\tilde a_i^j)_{tt}(0)$ and $(\tilde a_i^j)_{ttt}(0)$.

Let us now define 
\begin{subequations}
\begin{align}
%\label{bkappa}
b_{\kappa}(\phi)&= \kappa(c^{ijkl} {{w_{ 2}}}^k,_l, \phi^i,_j)_{L^2(\Omega_0^s;{\mathbb R})},\label{bkappa}\\
%\end{equation}
%\begin{equation}
%\label{ckappa}
c_{\kappa}(\phi)&=  \kappa (c^{ijkl}   w_1^k,_l, \phi^i,_j)_{L^2(\Omega_0^s;{\mathbb R})},\label{ckappa}\\
%\end{equation}
% \begin{equation}
%\label{dkappa}
d_{\kappa} (\phi)&=  \kappa (c^{ijkl}{u_0}^k,_l, \phi^i,_j)_{L^2(\Omega_0^s;{\mathbb R})}\ .\label{dkappa}
\end{align}
\end{subequations}

By proceding as in \cite{CoSh2004}, we can establish the existence of a 
fixed-point to the system (\ref{smoothie}). This follows the lines of 
\cite{CoSh2004} by first approximating by a penalty scheme the divergence-free 
constraint in the fluid in our Lagrangian setting, and by performing a regularity 
analysis of the solution of (\ref{linear}) allowing the use of the Tychonoff 
fixed-point theorem. 
Given the estimates obtained in
\cite{CoSh2004}, no new difficulty arise, since the parabolic
artificial viscosity in the solid controls the forcing coming from the quasilinear
part on a short time which is a priori shrinking to zero, and for this reason the 
proof is omitted here.

This leads us to the following
\begin{lemma}
There exists $T_{\kappa}>0$ depending a priori on $\kappa$ and on a given 
expression of the type $N_0(u_0,(w_i)_{i=1}^3)+N_0((q_i)_{i=0}^2)
+M_0(f,\kappa g,\kappa h)$, so that there exists a unique solution
$(w_{\kappa}, q_{\kappa})\in Z_{T_{\kappa}}$ of the regularized problem 
(\ref{smoothie}). Moreover, $w_\kappa \in V^4_s(T_\kappa)$.
\end{lemma}

In the next section we will study the limit of these solutions of the smoothed
problems as $\kappa\rightarrow 0$ (this being problematic since the solutions to these regularized problems are a priori defined on a time
interval shrinking to zero as $\kappa\rightarrow 0$).

Moreover, the following variational equations (for $n=0,1,2)$ are satisfied for any test function $\phi\in L^2(0,T_\kappa;H^1_0(\Omega;\R^3))$:	
 \begin{align}
& \int_0^{T_\kappa} (\partial_t^n (w_{\kappa})_t, \phi)_{L^2(\Omega;{\mathbb R}^3)}\ dt  
+\nu \int_0^{T_\kappa} ( \partial_t^n (a_k^r a_k^s  w_{\kappa},_r),\  \phi,_s)_{L^2(\Omega_0^f;{\mathbb R}^3)}\ dt \nonumber\\
&+ \kappa\int_0^{T_\kappa} (c^{ijkl} \partial_t^n {w_{\kappa}}^k,_l, {\phi}^i,_j)_{L^2(\Omega_0^s;{\mathbb R})}\ dt  - \int_0^{T_\kappa} ( \partial_t^n ( a_k^l q_{\kappa}),\  \phi^k,_l)_{L^2(\Omega_0^f;{\mathbb R})} dt\nonumber\\
&+ \int_0^{T_\kappa} (c^{ijkl} \partial_t^n [(\eta,_i\cdot\eta,_j-\delta_{ij}) \eta,_l], \phi,_k)_{L^2(\Omega_0^s;{\mathbb R}^3)}\ dt\n\\
&= \int_0^{T_\kappa} (\partial_t^n F, \phi)_{L^2(\Omega_0^f;{\mathbb R}^3)}+ (\partial_t^nf, \phi)_{L^2(\Omega_0^s;{\mathbb R}^3)} + \partial_t^n [\frac{t^2}{2}] b_{\kappa }(\phi)+ \partial_t^n[t] c_{\kappa }(\phi) +\partial_t^n[1] d_{\kappa}(\phi)\ dt\ ,
\label{weakW}
\end{align}
together with the initial conditions $w_{\kappa}(0)=u_0$, $(w_{\kappa})_t(0)=w_1$, $(w_{\kappa})_{tt}(0)=w_{2}$ and $ q_{\kappa}(0)=q_0$, $( q_{\kappa})_t(0)=q_{1}$, $( q_{\kappa})_{tt}(0)=q_{2}$. Moreover for the third time differentiated problem in time, we also have that a.e. in $(0,{T_\kappa})$,
\begin{align}
&  \frac{1}{2}\  \|(w_{\kappa})_{ttt} (t)\|^2_{L^2(\Omega;{\mathbb R}^3)} 
+{\nu} \int_0^{t} ((a_k^r  a_k^s w_{\kappa},_r)_{ttt},\  (w_{\kappa})_{ttt},_s)_{L^2(\Omega_0^f;{\mathbb R}^3)} 
\nonumber\\
&+\kappa \int_0^t (c^{ijkl}   {( w_{\kappa})_{ttt}}^k,_l (t), {( w_{\kappa})_{ttt}}^i,_j (t))_{L^2(\Omega_0^s;{\mathbb R})}\n\\
&- \int_0^{t}\int_{\Omega_0^f} ( q_{\kappa})_{tt}\ [3 ( a_i^j)_{tt} 
 (w_{\kappa})_t^i,_j +3 (a_i^j)_t (w_{\kappa})_{tt}^i,_j+ ( a_i^j)_{ttt} 
 { w_{\kappa}}^i,_j]_t \nonumber\\
&+ \int_{\Omega_0^f} (q_{\kappa})_{tt} (t)\  [3 ( a_i^j)_{tt} 
 (w_{\kappa})_t^i,_j +3 (a_i^j)_t (w_{\kappa})_{tt}^i,_j+ (a_i^j)_{ttt} 
 { w_{\kappa}}^i,_j](t)\n\\
& - \int_0^t \int_{\Omega_0^f} [\ 3 ( a_i^j)_{tt} (q_{\kappa})_t 
 ( w_{\kappa})_{ttt}^i,_j +  3( a_i^j)_{t} ( q_{\kappa})_{tt} (w_{\kappa})_{ttt}^i,_j + ( a_i^j)_{ttt} { q_{\kappa}}( w_{\kappa})_{ttt}^i,_j\ ] \nonumber\\
&+\int_0^t\int_{\os} c^{ijkl} [(\eta,_i\cdot\eta,_j-\delta_{ij}) \eta,_l]_{ttt}\cdot (w_{\kappa})_{ttt},_k\nonumber\\
&\le  N(u_0,(w_i)_{i=1}^3)+ N((q_i)_{i=0}^2) + \int_0^{t}[( F_{ttt}, ( w_{\kappa})_{ttt})_{L^2(\Omega_0^f;{\mathbb R}^3)} 
+  (f_{ttt},( w_{\kappa})_{ttt})_{L^2 (\Omega_0^s;{\mathbb R}^3)}] ,
\label{energywttt}
\end{align}
where, recall that $C$ does not depend on the artificial viscosity $\kappa$. 
The following result will be fundamental to our proof that the time interval
of existence of solutions to (\ref{smoothie}) is in fact $\kappa$-independent.
\begin{lemma}
\label{continuity}
The mapping $\gamma: t\rightarrow \|(w_\kappa,q_\kappa)\|_{Z_t}$ is continuous on $[0,T_\kappa]$.
\end{lemma}
\begin{proof}
The continuity with respect to $t$ of the terms of the type $L^2(0,t;H^s)$ is 
obvious, and since $w_{\kappa}\in V^4_s(T_\kappa)$ (due to our artificial 
viscosity), so is the continuity of
$\sum_{n=0}^3 \|\partial_t^n \eta_\kappa\|^2_{L^\infty(0,t;H^{4-n}(\os;\R^3))}$. 
The only terms that remain are 
$\|\partial_t^3 w_\kappa\|_{L^\infty(0,t;L^2(\Omega;\R^3))}$ and 
$\|\partial_t^2 q_\kappa\|_{L^\infty(0,t;L^2(\Omega_0^f;\R))}$. 

In order to treat them, we will invoke the fact that due to our artificial
viscosity in the solid, we in fact have $\partial_t^4 w_\kappa\in 
L^2(0,T_{\kappa};L^2(\Omega;\R^3))$, which provides $\partial_t^3 w_\kappa\in 
{\mathcal C}^0([0,T_\kappa];L^2(\Omega;\R^3))$. For the second time derivative of 
the pressure, we notice that from the variational form, true almost everywhere on
$[0,T_\kappa]$ for any $\phi\in H^1_0(\Omega;\R^3)$,
 \begin{align*}
&  (\partial_t^3 w_{\kappa}, \phi)_{L^2(\Omega;{\mathbb R}^3)}  
+\nu ( \partial_t^2 (a_k^r a_k^s  w_{\kappa},_r),\  \phi,_s)_{L^2(\Omega_0^f;{\mathbb R}^3)} + \kappa (c^{ijkl} \partial_t^2 {w_{\kappa}}^k,_l, {\phi}^i,_j)_{L^2(\Omega_0^s;{\mathbb R})}\\
&  - ( \partial_t^2 ( a_k^l q_{\kappa})-a_k^l {q_\kappa}_{tt},\  \phi^k,_l)_{L^2(\Omega_0^f;{\mathbb R})}+ (c^{ijkl} \partial_t^2 [(\eta,_i\cdot\eta,_j-\delta_{ij}) \eta,_l], \phi,_k)_{L^2(\Omega_0^s;{\mathbb R}^3)}\\
&- (\partial_t^2 F, \phi)_{L^2(\Omega_0^f;{\mathbb R}^3)}- (\partial_t^2 f, \phi)_{L^2(\Omega_0^s;{\mathbb R}^3)} - b_{\kappa }(\phi)= ( a_k^l {q_\kappa}_{tt},\  \phi^k,_l)_{L^2(\Omega_0^f;{\mathbb R})}\ ,
\end{align*}
and the Lagrange multiplier Lemma 13 of \cite{CoSh2004} associated to the continuity results previously established, we have the continuity of  $t\rightarrow \|{q_\kappa}_{tt}\|_{L^\infty(0,t;L^2(\of;\R))}$ on $[0,T_\kappa]$.

We now explain briefly why such a control on the fourth time derivative of 
$\w_\kappa$ holds, and is possible only with the addition of the artificial 
viscosity in the solid. In particular, this norm cannot be controlled as 
$\kappa\rightarrow 0$, which is not crucial for our purposes in any case. 
In order to understand the idea, we return to the level of the setting of the 
fixed-point argument, where we assume that $v$ in an appropriate convex set of 
$V^4_f(T)\times V^4_s(T)$ is given, and search for a solution $w$ of (\ref{linear})
by a Galerkin approximation on a penalized problem (for the pressure), in a way 
similar to \cite{CoSh2004}.  The penalization parameter  $\epsilon>0$ is given, and
we denote $\d q^n_\epsilon=\sum_{n=0}^2 \frac{t^n}{n!}q_n-\frac{1}{\epsilon} 
a_i^j (w^n_{\epsilon})^i,_j$, where $w^n_{\epsilon}$ is solution of the Galerkin 
approximation at rank $n$, and where $a_i^j$ is computed from $\eta$ associated 
to the given $v$. Our interest will be with the first problem that appears in our 
methodology in \cite{CoSh2004}; namely, the highest order time-differentiated
problem is multiplied by $\partial_t^4 w^n_{\epsilon}$ 
(which is permitted since it belongs to the appropriate finite dimensional space), 
and then integrate from $0$ to $t$.  We obtain
 \begin{align*}
& \int_0^{t} \|\partial_t^4 w^n_{\epsilon} \|^2_{L^2(\Omega;{\mathbb R}^3)}  
+\nu \int_0^{t} ( \partial_t^3 (a_k^r a_k^s  {w^n_{\epsilon}},_r),\  \partial_t^4 (w^n_{\epsilon}),_s)_{L^2(\Omega_0^f;{\mathbb R}^3)} \nonumber\\
&+ \bigl[\frac{\kappa}{2}(c^{ijkl} \partial_t^3 (w^n_{\epsilon})^k,_l, \partial_t^3 (w^n_{\epsilon})^i,_j)_{L^2(\Omega_0^s;{\mathbb R})}\bigr]_0^t - \int_0^{t} ( \partial_t^3 ( a_k^l q^n_\epsilon),\  \partial_t^4 (w^n_{\epsilon})^k,_l)_{L^2(\Omega_0^f;{\mathbb R})}\nonumber\\
&- \int_0^{t} (c^{ijkl} \partial_t^4 [(\eta,_i\cdot\eta,_j-\delta_{ij}) \eta,_l], \partial_t^3 (w^n_{\epsilon}),_k)_{L^2(\Omega_0^s;{\mathbb R}^3)}
\n\\
&+\bigl[\int_{\os} c^{ijkl} \partial_t^3 [(\eta,_i\cdot\eta,_j-\delta_{ij}) \eta,_l]\cdot \partial_t^3 (w^n_{\epsilon}),_k\bigr]_0^t\\
&\qquad\qquad\qquad\qquad\qquad\qquad\qquad= \int_0^{t} [\int_{\of} \partial_t^3 F\cdot \partial_t^4 w^n_{\epsilon}+ \int_{\os} \partial_t^3 f\cdot \partial_t^4 w^n_{\epsilon}],
\end{align*}
  leading us for a time small enough depending on the artificial viscosity $\kappa$ (but not on $n$ and $\epsilon$) to an inequality of the type,
   \begin{align*}
 \int_0^{t} \|\partial_t^4 w^n_{\epsilon} \|^2_{L^2(\Omega;{\mathbb R}^3)}&  
+\sup_{[0,t]} [\| \partial_t^3 {w^n_{\epsilon}}\|^2_{H^1(\Omega_0^f;{\mathbb R}^3)}+ \kappa \| \partial_t^3 {w^n_{\epsilon}}\|^2_{H^1(\Omega_0^s;{\mathbb R}^3)}+\epsilon \| \partial_t^3 {q^n_{\epsilon}}\|^2_{L^2(\Omega_0^f;{\mathbb R})}] \nonumber\\
&\le C_\epsilon [N(u_0,(w_i)_{i=1}^3)+N((q_i)_{i=0}^2)+N(f)],
\end{align*}
where $C_{\epsilon}$ depends a priori on $\epsilon$. By proceding in a way inspired
by our methodology in Section 9 of \cite{CoSh2004}, we can then prove that we have
control, independently of $\epsilon$, on the first three norms. Taking
the limit first as $\n\rightarrow\infty$ and then as $\epsilon\rightarrow 0$, 
indeed provides us with $\partial_t^4 w_\kappa\in L^2(0,T_\kappa;L^2(\of;\R^3))$ 
as announced. 
\end{proof}
We note that this latter regularity property in the solid is only
possible with the artificial viscosity $\kappa>0$. 

\section{Estimate for the solutions of (\ref{linear}) independently of $\kappa$}
\label{8}

In this section, we will denote $(w_{\kappa},q_{\kappa})=(\w,\q)$ and denote the corresponding quantities $a_i^j$ by $\a_i^j$. In what follows, $\delta>0$ is a 
given positive number to be made precise later when it will be chosen to be 
sufficiently small.

\subsection{Energy estimate for $\w_{ttt}$ independently of $\kappa$.} We are now going to use the regularity result $(\w,\q)\in Z_{T_\kappa}$ in the energy inequality (\ref{energywttt}) (which was established independently of the artificial viscosity), this time by interpolating and using the energy properties of the nonlinear elasticity operator, in order to get an estimate
independent of the artificial viscosity.

\noindent Step 1. Let $\d I_1=\int_0^t\int_{\os} c^{ijkl} (\e,_i\cdot\e,_j-\delta_{ij}) {\w_{tt}},_l\cdot \w_{ttt},_k$.
An integration by parts in time shows that
\begin{align*}
\d I_1=\ -\frac{1}{2}\int_0^t\int_{\os} c^{ijkl} (\e,_i\cdot\e,_j-\delta_{ij})_t &{\w_{tt}},_l\cdot \w_{tt},_k\\
& + \frac{1}{2}\int_{\os}  c^{ijkl} (\e,_i\cdot\e,_j-\delta_{ij}) {\w_{tt}},_l\cdot \w_{tt},_k (t)\ ,
\end{align*}
and thus with the properties of the Bochner integral in $H^2(\os;\R)$, 
\begin{align*}
\d\e,_i\cdot\e,_j (t)-\delta_{ij}
=\int_0^t [\e,_i\cdot\w,_j+\w,_i\cdot\e,_j],
\end{align*}
we deduce
\begin{align}
|I_1|\le &  C t\ \|\w_{tt}\|^2_{L^{\infty}(0,t;H^1(\os;\R^3))}\|\w\|_{L^{\infty}(0,t;H^3(\os;\R^3))}\|\e\|_{L^{\infty}(0,t;H^3(\os;\R^3))}\n\\
\le & Ct\ \|\w\|^4_{W_t}  .
\label{I1}
\end{align}

\noindent Step 2. Let $\d I_2=3\int_0^t\int_{\os} c^{ijkl} (\e,_i\cdot\e,_j-\delta_{ij})_t {\w_{t}},_l\cdot \w_{ttt},_k$.
Similarly,
\begin{align*}
\d I_2=&\ -3\int_0^t\int_{\os} c^{ijkl} [(\e,_i\cdot\e,_j-\delta_{ij})_t {\w_{tt}},_l\cdot \w_{tt},_k+ (\e,_i\cdot\e,_j-\delta_{ij})_{tt} {\w_{t}},_l\cdot \w_{tt},_k]\\
& + 3\bigl[\int_{\os}  c^{ijkl} (\e,_i\cdot\e,_j-\delta_{ij})_t {\w_{t}},_l\cdot \w_{tt},_k (t)\bigr]_0^t.
\end{align*}
By the same type of argument as for the previous step, we then get
\begin{align*}
\d |I_2|\le &\ C t\  \|\w_{tt}\|^2_{L^{\infty}(0,t;H^1(\os;\R^3))}\|\w\|_{L^{\infty}(0,t;H^3(\os;\R^3))}\|\e\|_{L^{\infty}(0,t;H^3(\os;\R^3))}\n\\
&\ + C t\  \|\w_{tt}\|_{L^{\infty}(0,t;H^1(\os;\R^3))}\|\w\|^2_{L^{\infty}(0,t;H^3(\os;\R^3))}\|\w_t\|_{L^{\infty}(0,t;H^1(\os;\R^3))}\n\\
&\ + C t\  \|\w_{tt}\|_{L^{\infty}(0,t;H^1(\os;\R^3))}\|\w_t\|^2_{L^{\infty}(0,t;H^2(\os;\R^3))}\|\e\|_{L^{\infty}(0,t;H^3(\os;\R^3))}\n\\
&\ +C\ \|c^{ijkl} [( u_0^i,_j+u_0^j,_i) w_1,_l + \int_0^\cdot
((\e,_i\cdot\e,_j-\delta_{ij})_t {\w_{t}},_l)_t]\|_{L^\infty(0,t;L^2(\os;\R^3))} \n\\
&\qquad\qquad\qquad\qquad\qquad\qquad\qquad\qquad\qquad\qquad\times \|{\w,_k}_{tt}\|_{L^\infty(0,t;L^2(\os;\R^3))}\\
&\ +N(u_0,(w_i)_{i=1}^3), 
\end{align*}
and thus,
\begin{align}
\d |I_2|\le \delta \|\w\|^2_{W_t} +C t \|\w\|^4_{W_t} +C_{\delta} N(u_0,(w_i)_{i=1}^3)\ .
\label{I2}
\end{align}

\noindent Step 3. Let $\d I_3=3\int_0^t\int_{\os} c^{ijkl} (\e,_i\cdot\e,_j-\delta_{ij})_{tt} {\w},_l\cdot \w_{ttt},_k$.
By an integration by parts in time,
\begin{align*}
\d I_3=&\ -3 \int_0^t\int_{\os} c^{ijkl} [ (\e,_i\cdot\e,_j-\delta_{ij})_{tt} {\w_{t}},_l\cdot \w_{tt},_k+ (\e,_i\cdot\e,_j-\delta_{ij})_{ttt} {\w},_l\cdot \w_{tt},_k ] \\
& + \bigl[\int_{\os}  c^{ijkl} (\e,_i\cdot\e,_j-\delta_{ij})_{tt} {\w},_l\cdot \w_{tt},_k (t)\bigr]_0^t.
\end{align*}
Similarly as before, we get
\begin{align*}
\d |I_3|& \le   C t\  \|\w_{tt}\|_{L^{\infty}(0,T;H^1(\os;\R^3))}\|\w\|^2_{L^{\infty}(0,t;H^3(\os;\R^3))}\|\w_t\|_{L^{\infty}(0,t;H^1(\os;\R^3))}\n\\
&\ + C t\  \|\w_{tt}\|_{L^{\infty}(0,T;H^1(\os;\R^3))}\|\w_t\|^2_{L^{\infty}(0,t;H^2(\os;\R^3))}\|\e\|_{L^{\infty}(0,t;H^3(\os;\R^3))}\n\\
&\ + C t\  \|\w_{tt}\|^2_{L^{\infty}(0,T;H^1(\os;\R^3))}\|\w\|_{L^{\infty}(0,t;H^3(\os;\R^3))}\|\e\|_{L^{\infty}(0,t;H^3(\os;\R^3))}\n\\
&\ + C t\  \|\w_{tt}\|_{L^{\infty}(0,T;H^1(\os;\R^3))}\|\w_t\|_{L^{\infty}(0,t;H^1(\os;\R^3))}\|\w\|^2_{L^{\infty}(0,t;H^3(\os;\R^3))}\n\\
&\ +C\ \|c^{ijkl} ((\e,_i\cdot\e,_j)_{tt}(0) w_1,_l + \int_0^\cdot
((\e,_i\cdot\e,_j-\delta_{ij})_{tt} {\w},_l)_t)\|_{L^\infty(0,t;L^2(\os;\R^3))} \n\\
&\qquad\qquad\qquad\qquad\qquad\qquad\qquad\qquad\qquad\qquad\times \|\w_{tt},_k\|_{L^\infty(0,t;L^2(\os;\R^3))}\\
&+N(u_0,(w_i)_{i=1}^3), 
\end{align*}
and therefore
\begin{align}
\label{I3}
\d |I_3|\le &\ \delta \|\w\|^2_{W_t}+C t \|\w\|^4_{W_t} +C_{\delta} N(u_0,(w_i)_{i=1}^3). 
\end{align}

\noindent Step 4. Let $\d I_4=\int_0^t\int_{\os} c^{ijkl} (\e,_i\cdot\e,_j-\delta_{ij})_{ttt} {\e},_l\cdot \w_{ttt},_k$.
By symmetry of $c$, we notice that
\begin{align*}
I_4=&\ \frac{1}{2}\int_0^t\int_{\os} c^{ijkl} (\e,_i\cdot{\e,_j})_{ttt} ({\e},_l\cdot {\e,_k})_{tttt}\\
& -\int_0^t\int_{\os} c^{ijkl}[4 (\e,_i\cdot{\e,_j})_{ttt} ({\w},_l\cdot {\w_{tt}},_k)+ 3 (\e,_i\cdot{\e,_j})_{ttt} ({\w_t},_l\cdot {\w_{t}},_k)],
\end{align*}
and thus,
\begin{align}
\bigl|I_4&-  \frac{1}{4} \int_{\os} c^{ijkl} (\e,_i\cdot{\e,_j})_{ttt} ({\e},_l\cdot {\e,_k})_{ttt} (t)\bigr|\n\\
\le\ & Ct\ \|\nabla\w_{tt}\|^2_{L^{\infty}(0,t;L^2(\os;\R^9))}\|\nabla\w\|_{L^{\infty}(0,t;H^2(\os;\R^9))}\|\nabla\e\|_{L^{\infty}(0,t;H^2(\os;\R^9))}\n\\
& +Ct\ \|\nabla\w_{t}\|_{L^{\infty}(0,t;L^2(\os;\R^9))}\|\nabla\w_{tt}\|_{L^{\infty}(0,t;L^2(\os;\R^9))}\|\nabla\w\|^2_{L^{\infty}(0,t;H^2(\os;\R^9))}\n\\
& +Ct\ \|\nabla\w_t\|^2_{L^{\infty}(0,t;H^1(\os;\R^9))}\|\nabla\w_{tt}\|_{L^{\infty}(0,t;L^2(\os;\R^9))}\|\nabla\e\|_{L^{\infty}(0,t;H^2(\os;\R^9))}\n\\
& +Ct\ \|\nabla\w\|_{L^{\infty}(0,t;H^2(\os;\R^9))}\|\nabla\w_{t}\|^3_{L^{\infty}(0,t;H^1(\os;\R^9))}+N(u_0,(w_i)_{i=1}^3) \n\\
&\le N(u_0,(w_i)_{i=1}^3)+ Ct\ \|\w\|^4_{W_t}.
\label{I4}
\end{align}

Now, Let $\d I=\frac{1}{4} \int_{\os} c^{ijkl} (\e,_i\cdot{\e,_j})_{ttt} ({\e},_l\cdot {\e,_k})_{ttt} (t)$. 
By expanding the integrand with respect to the time derivatives and using the relation in $H^3(\os;\R^3)$: $\d\e(t,\cdot)=\text{Id}+\int_0^t \w(t',\cdot)\ dt'$ and estimates similar as in the previous steps, we find that
\begin{align}
\bigl|I-  \int_{\os} c^{ijkl} \w_{tt}^j,_i\w_{tt}^k,_l (t)\bigr|\le & C_{\delta} N(u_0,(w_i)_{i=1}^3)+\delta \|\w\|^2_{W_t}+C t \|\w\|^4_{W_t}.
\label{I}
\end{align}

\noindent{Step 5.}
By using (\ref{I4}) and (\ref{I}) we find that
\begin{align}
\bigl|I_4-  \int_{\os} c^{ijkl} \w_{tt}^j,_i\w_{tt}^k,_l (t)\bigr|\le & C_{\delta} N(u_0,(w_i)_{i=1}^3)+\delta \|\w\|^2_{W_t}+ C t \|\w\|^4_{W_t}.
\label{I5}
\end{align}

\noindent Step 6. By proceding in a way similar to \cite{CoSh2004}, except that we replace the constants $C(M)$ appearing there by appropriate powers of $\|(\w,\q)\|_{Z_t}$, we find that the integrals set in the fluid domain are bounded by
\begin{align*}
\delta \|(\w,\q)\|^2_{Z_t}+ C_{\delta}  [N((q_i)_{i=0}^2)+N(u_0,(w_i)_{i=1}^3)+M(f,\kappa g,\kappa h)+t^{\frac{1}{4}}\|(\w,\q)\|^6_{Z_t}].
\end{align*}

\noindent{Step 7.} 
Thus, from (\ref{energywttt}), and Steps 1-6, we then get on $[0,T_\kappa]$:
\begin{align}
 & \sup_{[0,t]} \|\tilde w_{ttt} \|^2_{L^2(\Omega;{\mathbb R}^3)} 
+\int_0^{t} [\ \|{\tilde w}_{ttt}\|^2_{H^1(\Omega_0^f;{\mathbb R}^3)}
+\kappa \|{\tilde w}_{ttt}\|^2_{H^1(\Omega_0^s;{\mathbb R}^3)}] 
+ \sup_{[0,t]} \|{{\tilde w}_{tt}}\|^2_{H^1(\Omega_0^s;{\mathbb R}^3)}
\nonumber\\
& \le  C_{\delta} [N(u_0,(w_i)_{i=1}^3)+M(f,\kappa g,\kappa h)+N((q_i)_{i=0}^2)+ t^{\frac{1}{4}}\|(\tilde w,\tilde q)\|^6_{Z_t}] + C {\delta} \|(\tilde w,\tilde q)\|^2_{Z_t} .
\label{I7}
\end{align}

\noindent Step 8. The estimate of $\q_{tt}$ in $L^\infty(L^2)$, independently of $\kappa$, will require some adjustments with respect to the methodology of \cite{CoSh2004}. To this end, we notice that we can apply a Lagrange multiplier Lemma similar to Lemma 13 of \cite{CoSh2004}, but corresponding to the case $a_i^j=\delta_i^j$, to the variational form true on $[0,T_{\kappa}]$: for all $\phi\in H^1_0(\Omega;\R^3)$,	
 \begin{align*}
&  (\w_{ttt}, \phi)_{L^2(\Omega;{\mathbb R}^3)}  
+\nu ( (\a_k^r \a_k^s  \w,_r)_{tt},\  \phi,_s)_{L^2(\Omega_0^f;{\mathbb R}^3)} + \kappa (c^{ijkl}  {\w}_{tt}^k,_l, {\phi}^i,_j)_{L^2(\Omega_0^s;{\mathbb R})}\\
& +  (c^{ijkl} [(\e,_i\cdot\e,_j-\delta_{ij}) \e,_l]_{tt}, \phi,_k)_{L^2(\Omega_0^s;{\mathbb R}^3)}-  (  (\a_k^l-\delta_{kl}) \q_{tt},\  \phi^k,_l)_{L^2(\Omega_0^f;{\mathbb R})}\\
& -  (  (\a_k^l \q)_{tt}-\a_k^l \q_{tt},\  \phi^k,_l)_{L^2(\Omega_0^f;{\mathbb R})}-  ( F_{tt}, \phi)_{L^2(\Omega_0^f;{\mathbb R}^3)}- (f_{tt}, \phi)_{L^2(\Omega_0^s;{\mathbb R}^3)}  -  b_{\kappa }(\phi)\\
&\qquad\qquad\qquad\qquad\qquad\qquad\qquad\qquad\qquad\qquad\qquad\qquad\qquad =(   \q_{tt},\  \operatorname{div} \phi)_{L^2(\Omega_0^f;{\mathbb R})},
\end{align*}
which provides
for any $t\in [0,T_{\kappa}]$:
\begin{align*}
\|\q_{tt}\|_{L^2(\Omega_0^f;\R)}\le&\ C [\|\w_{ttt}\|_{L^2(\Omega;{\mathbb R}^3)}  
+\|(\a_k^r \a_k^s  \w,_r)_{tt}\|_{L^2(\Omega_0^f;{\mathbb R}^3)} + \kappa\| {\w}_{tt}\|_{H^1(\Omega_0^s;{\mathbb R})}\n\\
& +  \|[(\e,_i\cdot\e,_j-\delta_{ij}) \e,_l]_{tt}\|_{L^2(\Omega_0^s;{\mathbb R}^3)}+  \|(\a_k^l \q)_{tt}-\a_k^l \q_{tt}\|_{L^2(\Omega;{\mathbb R})}\n\\
& + \|\a-\text{Id}\|_{H^2(\of;\R^9)}\|\q_{tt}\|_{L^2(\of;\R)} +N(u_0,(w_i)_{i=1}^3)+M(f,\kappa g,\kappa h)].
\end{align*}
By using (\ref{I7}) for the first four terms of the right-hand side of this inequality and remembering that the $L^\infty(0,t;L^2(\of;\R))$ norm of $\q_{tt}$ is part of the norm $Z_t$ for the next two terms of this inequality, we get 
\begin{align}
  \sup_{[0,t]} \|\tilde q_{tt} \|^2_{L^2(\Omega_0^f;{\mathbb R})} 
 \le&\  C_{\delta} [N(u_0,(w_i)_{i=1}^3)+M(f,\kappa g,\kappa h)+N((q_i)_{i=0}^2)]\n\\
& +C_{\delta} t^{\frac{1}{4}} \|(\tilde w,\tilde q)\|^6_{Z_t} + C {\delta} \|(\tilde w,\tilde q)\|^2_{Z_t} .
\label{I8}
\end{align}

\subsection{Estimate on $w_{tt}$ and $w_t$.}
 From the previous estimates, and the arguments that we will see hereafter for the case of $\w$, we have
\begin{align}
\|\w_{tt}\|^2_{L^2(0,t;H^2(\of;\R^3))}&+ \|\q_{tt}\|^2_{L^2(0,t;H^1(\of;\R))}+\|\w_{t}\|^2_{L^{\infty}(0,t;H^2(\os;\R^3))}\n\\
& \le C_{\delta} [N(u_0,(w_i)_{i=1}^3)+M(f,\kappa g,\kappa h)+N((q_i)_{i=0}^2)]\n\\
&\ \ \  +C_{\delta} t^{\frac{1}{4}} \|(\tilde w,\tilde q)\|^6_{Z_t} + C {\delta} \|(\tilde w,\tilde q)\|^2_{Z_t}. 
\label{wtth2}
\end{align}
 Similarly, we infer from (\ref{wtth2}) that 
\begin{align}
\|\w_{t}\|^2_{L^2(0,t;H^3(\of;\R^3))}&+ \|\q_{t}\|^2_{L^2(0,t;H^2(\of;\R))}+\|\w\|^2_{L^{\infty}(0,t;H^3(\os;\R^3))}\n\\
& \le C_{\delta} [N(u_0,(w_i)_{i=1}^3)+M(f,\kappa g,\kappa h)+N((q_i)_{i=0}^2)]\n\\
&\ \ +C_{\delta} t^{\frac{1}{4}} \|(\tilde w,\tilde q)\|^6_{Z_t} + C {\delta} \|(\tilde w,\tilde q)\|^2_{Z_t}. 
\label{wth3f}
\end{align}

\subsection{Estimate on $\w$}\hfill\break
We will denote $\R^3_+=\{x\in\R^3|\ x_3> 0\}$, $\R^3_-=\{x\in\R^3|\ x_3< 0\}$, and $B_-(0,r)=B(0,r)\cap\R^3_-$.
 We denote by $\P$ an $H^4$ diffeomorphism from $B(0,1)$ into a neighborhood $V$ of a point $x_0\in\Gamma_0$ such that
$\P(B(0,1)\cap\R^3_+)=V\cap\of$, $\P(B(0,1)\cap\R^3_-)=V\cap\os$,
$\P(B(0,1)\cap\R^2\times\{0\})=V\cap\Gamma_0$, with $\det\nabla\P=1$. We consider a cut-off function $\z$ compactly supported in $B(0,1)$, and equal to $1$ in $B(0,\frac{1}{2})$. 
%For $h\ne 0$, we denote the standard difference quotient $\d D_h u(x)=\frac{u(x+h)-u(x)}{|h|}$ and denote $D_{hk}u=D_h[D_k u]$. Moreover, $h_i$ ($i=1,2,3)$ being an arbitrary non-zero vector, we will denote $D_{h_i}=D_i$ and $D_{-h_i}=D_{-i}$.

With the use of test functions $\phi_p=-[\rho_p\star (\z^2\ \w\circ\P)],_{\ao\ao\at\at\att\att}\circ\P^{-1}$ (which is in $L^2(0,T_\kappa; H^1_0(\Omega;\R^3))$) in (\ref{weakW}) for $n=0$, and by denoting $W=\w\circ\P$, $Q=\q\circ\P$, $E=\e\circ\P$, we get after integrating by parts appropriately and letting $p\rightarrow\infty$,
\begin{align}
&  \frac{1}{2}\  \|\z   W,_{\ao\at\att} (t)\|^2_{L^2(\R^3;{\mathbb R}^3)} + \int_0^t  ( {W_{t}},_{\ao\at\att},\ [\z^2 W],_{\ao\at\att}-\z^2 W,_{\ao\at\att})_{L^2(\R^3;\R^3)}\n\\
&+{\nu} \int_0^{t} ([\tilde b_k^r \tilde b_k^s W,_r],_{\ao\at\att},\  [\z^2 W],_{ s\ao\at\att})_{L^2(\R^3_+;{\mathbb R}^3)}\n\\
&- \int_0^{t}\int_{\R^3_+} [Q\ \tilde b_i^j],_{\ao\at\att} [\z^2 W]^i,_{j\ao\at\att}\n\\
&+\kappa\int_0^t ([C^{ijkl} (W,_i\cdot \P,_j)\ \P,_l],_{\ao\at\att},\ [\z^2 W],_{ k\ao\at\att})_{L^2(\R^3_-;{\mathbb R}^3)}\nonumber\\ 
&+\int_0^t ([C^{ijkl} (E,_i\cdot E,_j-\P,_i\cdot\P,_j)\ E,_l],_{\ao\at\att},\ [\z^2 W],_{ k\ao\at\att})_{L^2(\R^3_-;{\mathbb R}^3)}\nonumber\\
&\qquad\qquad \le C\  N(u_0,(w_i)_{i=1}^3) + \lim_{p\rightarrow\infty}\int_0^{t} [\ ( F, \phi_p)_{L^2(\Omega_0^f;{\mathbb R}^3)} 
+  (f,\phi_p)_{L^2 (\Omega_0^s;{\mathbb R}^3)}] \n\\
&\qquad\qquad\  + \lim_{p\rightarrow\infty}\int_0^t [\ \frac{t'^2}{2} b_{\kappa}(\phi_p)+t' c_{\kappa}(\phi_p)+d_{\kappa}(\phi_p)\ ],
\label{energyw}
\end{align}
where 
$C^{ijkl}= c^{mnop} g^i_m g^j_n g^k_o g^l_p\in H^3(B(0,1);\R)$, $g=[\nabla\P]^{-1}\in H^3(B(0,1);\R^9)$, $\tilde b_l^j= \a_l^k(\P) g_k^j$.

\begin{remark}
Note that this limit process as $p\rightarrow\infty$ for the nonlinear elastic energy is possible because $\partial_t^n\w\in L^2(0,T_\kappa; H^{4-n}(\os;\R^3))$ ($n=0,1$) due to our artificial viscosity in the solid. Whereas we could also use difference quotients, it appears that the product rules are less cumbersome with the use of horizontal derivatives instead, which is permitted since we already know at this stage the regularity of $\w$ and $\q$. Also, the limits on the
right-hand side of (\ref{energyw}) do not present any difficulty, given the
regularity of the forcing functions and three integrations by parts with respect 
to horizontal variables.
\end{remark}

\begin{remark}
Since $\z$ is compactly supported in $B(0,1)$, the integrals set on $\R^3$, $\R^3_-$, $\R^3_+$ do not depend on the extension that we chose for $W$, $E$ or $Q$, and 
simply represent a more convenient way to write these integrals.
\end{remark}

\noindent Step 1. Let $\d L_1=\kappa\int_0^t(\ [C^{ijkl} W,_i\cdot\P,_j\ \P,_l],_{\ao\at\att}, [\z^2\ W],_{ k\ao\at\att})_{L^2(\R^3_-;{\mathbb R}^3)}$.
By using the $H^3$ regularity of the coefficients $C^{ijkl}$, 
\begin{align}
 L_1&=\kappa \int_0^t( C^{ijkl} W,_{i\ao\at\att}\cdot\P,_j\ \P,_l,\ \z^2\  W,_{k\ao\at\att})_{L^2(\R^3_-;{\mathbb R}^3)}\n\\
&\ \ +\kappa \int_0^t([C^{ijkl} W,_i\cdot\P,_j\ \P,_l],_{\ao\at\att} - C^{ijkl}  W,_{i\ao\at\att}\cdot\P,_j\ \P,_l,\ \z^2\  W,_{ k\ao\at\att})_{L^2(\R^3_-;{\mathbb R}^3)}\n\\
&\ \ +\kappa \int_0^t( [C^{ijkl} W,_i\cdot\P,_j\ \P,_l],_{\ao\at\att},\ [\z^2\ W,_{k}],_{\ao\at\att} - \z^2\ W,_{ k\ao\at\att})_{L^2(\R^3_-;{\mathbb R}^3)}\n\\
& \ge C \kappa \int_0^t ( C^{ijkl} W,_{ i\ao\at\att}\cdot\P,_j\ \P,_l,\ \z^2\ W,_{ k\ao\at\att})_{L^2(\R^3_-;{\mathbb R}^3)} -C \kappa \int_0^t \|\w\|^2_{H^3(\os;\R^{3})}\n\\
& \ge C \kappa \int_0^t \|W,_{\ao\at\att}\|^2_{H^1(B_-(0,\frac{1}{2});\R^3)}-C \kappa t \|\w\|^2_{W_t}.
\label{L1}
\end{align}

\noindent Step 2. Let 
$$
L_2=\d \int_0^t ([C^{ijkl} (E,_i\cdot E,_j-\P,_i\cdot\P,_j)\ E,_l],_{\ao\at\att},\ [\z^2 E_{t}],_{ k\ao\at\att})_{L^2(\R^3_-;{\mathbb R}^3)}
.$$  $\Sigma_3$ denoting the set of permutations of
$\{1,2,3\}$, 
\begin{align*}
L_2=& \int_0^t (C^{ijkl} (E,_i\cdot E,_j-\P,_i\cdot\P,_j)\ E,_{ l\ao\at\att},\ \z^2 E_{t},_{ k\ao\at\att })_{L^2(\R^3_-;{\mathbb R}^3)}\\
&+ 2 \int_0^t (\z^2 C^{ijkl} (E,_{ i\ao\at\att}\cdot E,_j)\ E,_l,\  E_t,_{ k\ao\at\att })_{L^2(\R^3_-;{\mathbb R}^3)}\\
&+ \int_0^t (\z^2 [ [C^{ijkl} (E,_{i}\cdot E,_j)],_{\ao\at\att}-
2 C^{ijkl} (E,_{ i\ao\at\att}\cdot E,_j)] E,_l
,\  E_t,_{ k\ao\at\att })_{L^2(\R^3_-;{\mathbb R}^3)}\\
&+  \int_0^t ([\z^2 (C^{ijkl} \P,_i\cdot\P,_j),_{\ao\at\att}\ E,_l]_t,\   E,_{ k\ao\at\att })_{L^2(\R^3_-;{\mathbb R}^3)}\\
& -\bigl[(\z^2 (C^{ijkl} \P,_i\cdot\P,_j),_{\ao\at\att}\ E,_l,\   E,_{ k\ao\at\att })_{L^2(\R^3_-;{\mathbb R}^3)}\bigr]_0^t\\
&- \sum_{\sigma\in\Sigma_3} \int_0^t ([\z^2[C^{ijkl} (E,_i\cdot E,_j-\P,_i\cdot\P,_j)],_{\alpha_{\sigma(1)}} E,_{\alpha_{\sigma(2)}\alpha_{\sigma(3)} l}],_{\ao},\  E_{t},_{\at\att k})_{L^2(\R^3_-;{\mathbb R}^3)}\\
&- \sum_{\sigma\in\Sigma_3} \int_0^t ([\z^2 [C^{ijkl} (E,_i\cdot E,_j-\P,_i\cdot\P,_j)],_{\alpha_{\sigma(1)}\alpha_{\sigma(2)}} E,_{\alpha_{\sigma(3)} l}],_{\ao},\  E_{t},_{\at\att k})_{L^2(\R^3_-;{\mathbb R}^3)}\\
&+ \int_0^t ([C^{ijkl} (E,_i\cdot E,_j-\P,_i\cdot\P,_j)\ E,_{l}],_{\ao\at\att}, [\z^2 E_{t}],_{\ao\at\att k}-\z^2 E_{t},_{\ao\at\att k })_{L^2(\R^3_-;{\mathbb R}^3)}\ .
\end{align*}
From the regularity of $\w$ and the $H^4$ regularity of $\P$, we then infer
\begin{align*}
L_2=& \int_0^t (C^{ijkl} (E,_i\cdot E,_j-\P,_i\cdot\P,_j)\ E,_{\ao\at\att l},\ \z^2 E_{t},_{\ao\at\att k })_{L^2(\R^3_-;{\mathbb R}^3)}\n\\
& + 2\int_0^t (C^{ijkl} E,_{\ao\at\att i}\cdot E,_j\ \z^2 E,_l,\ E_{t},_{\ao\at\att k})_{L^2(\R^3_-;{\mathbb R}^3)}\n\\
&  + \int_0^t (C,_{\ao\at\att}^{ijkl} E,_{i}\cdot E,_j\ \z^2 E,_l,\ E_{t},_{\ao\at\att k})_{L^2(\R^3_-;{\mathbb R}^3)} + L_2^r,
\end{align*}
with 
\begin{equation}
\label{L22}
|L_2^r|\le \delta\|\w\|^2_{W_t} + C t\ \|\w\|^4_{W_t}+C_{\delta} N(u_0,(w_i)_{i=1}^3).
\end{equation}
By integrating by parts in time, we deduce
\begin{align*}
L_2=& -\frac{1}{2}\int_0^t (C^{ijkl} (E,_i\cdot E,_j-\P,_i\cdot\P,_j)_t\ E,_{\ao\at\att l},\ \z^2 E,_{\ao\at\att k })_{L^2(\R^3_-;{\mathbb R}^3)}\\
&+\bigl[\frac{1}{2} (C^{ijkl} (E,_i\cdot E,_j-\P,_i\cdot\P,_j)\ E,_{\ao\at\att l},\ \z^2 E,_{\ao\at\att k })_{L^2(\R^3_-;{\mathbb R}^3)}\bigr]_0^t\\
& - 2 \int_0^t (C^{ijkl} E,_{\ao\at\att i}\cdot E,_{j}\ \z^2 E_t,_l,\ E,_{\ao\at\att k})_{L^2(\R^3_-;{\mathbb R}^3)}\\
& + \bigl[ (C^{ijkl} E,_{\ao\at\att i}\cdot E,_{ j}\ \z^2 E,_l,\ E,_{\ao\at\att k})_{L^2(\R^3_-;{\mathbb R}^3)}\bigr]_0^t \n\\
&  - \int_0^t (C,_{\ao\at\att}^{ijkl} (E,_{i}\cdot E,_j\ \z^2 E,_l)_t,\ E,_{\ao\at\att k})_{L^2(\R^3_-;{\mathbb R}^3)}\n\\
&  + \bigl[(C,_{\ao\at\att}^{ijkl} E,_{i}\cdot E,_j\ \z^2 E,_l,\ E,_{\ao\at\att k})_{L^2(\R^3_-;{\mathbb R}^3)}\bigr]_0^t + L_2^r,
\end{align*}
which implies in turn
\begin{align}
\bigl|L_2- (C^{ijkl} \z^2 E,_{\ao\at\att i}(t)&\cdot \P,_j\ \P,_l,\ E,_{\ao\at\att k}(t))_{L^2(\R^3_-;{\mathbb R}^3)}\bigr| \n\\
&\le \delta \|\w\|^2_{W_t}+C_{\delta} N(u_0,(w_i)_{i=1}^3)+ C t\ \|\w\|^4_{W_t}.
\label{L2}
\end{align}
With $e_k$ ($k=1,2,3$) denoting the canonical vectors of $\R^3$, let $$\d P(t)=\bigl\|\ 2 c^{ijkl} (\e,_{mnil}\cdot e_j)\ e_k- 
[c^{ijkl} (\e,_i\cdot \e_j-\delta_{ij})\e,_k],_{lmn} \bigr\|_{L^2(\os;\R^3)}(t),$$  where $m$ and $n$ are arbitrarily fixed in $\{1,2,3\}$. We then have $$P(t)\le P_1(t)+P_2(t),$$ with
\begin{align*}
P_1(t)&= \bigl\|\ 2c^{ijkl} (\e,_{ilmn}\cdot e,_j)\ e_k-2 
[c^{ijkl} (\e,_{mni}\cdot \e,_j)\e,_k],_l \bigr\|_{L^2(\os;\R^3)}(t),\\
P_2(t)&= \bigl\| 
-[c^{ijkl} (\e,_i\cdot \e,_j-\delta_{ij})\e,_k],_{lmn}+ 2[c^{ijkl} (\e,_{mni}\cdot \e,_j)\e,_k],_l \bigr\|_{L^2(\os;\R^3)}(t).
\end{align*}
We first notice that 
\begin{align*}
P_1(t) \le &  \bigl\|\ 2c^{ijkl} [ (\e,_{ilmn}\cdot e_j)\ e_k- 
(\e,_{ilmn}\cdot \e,_j)\e,_k] \bigr\|_{L^2(\os;\R^3)}(t)\\
&+ \bigl\|\ 2c^{ijkl} [ \e,_{imn}\cdot (\e,_j]\ \e,_k),_l \bigr\|_{L^2(\os;\R^3)}(t).
\end{align*}
Next, by writing $\d\e(t)=\text{Id}+\int_0^t \w$ and $\d [\e,_{mni}\cdot(\e,_j] \ \e,_{k}),_l (t)=[\e,_{mni}\cdot(\e,_j]\ \e,_{k}),_l (0)+\int_0^t [[\e,_{mni}\cdot(\e,_j]\ \e,_{k}),_l]_t$ respectively in $H^3(\os;\R^3)$ and $L^2(\os;\R^3)$, we obtain 
\begin{align*}
P_1(t) & \le C [\int_0^t\|\w\|_{H^3(\os;\R^3)}]\ \sup_{[0,t]}[\|\e\|^2_{H^4(\os;\R^3)}+\|\e\|_{H^4(\os;\R^3)}]\n\\
&\ \ + N(u_0,(w_i)_{i=1}^3)+ C\int_0^t \|\w\|_{H^3(\os;\R^3)}\sup_{[0,t]}\|\e\|^2_{H^4(\os;\R^3)}\n\\
&\le  N(u_0,(w_i)_{i=1}^3)+ C t\ \|\w\|^3_{W_t}.
\end{align*}
Next, we see that
\begin{align*}
P_2(t)\le&\ \bigl\| 
c^{ijkl} [ (\e,_i\cdot \e,_j-\delta_{ij})\e,_{kmn}+ (\e,_{i}\cdot \e,_j),_m\e,_{kn}+(\e,_{i}\cdot \e,_j),_n\e,_{km}],_l \bigr\|_{L^2(\os;\R^3)}(t)\\
&+ \bigl\| 
c^{ijkl} [ (\e,_{mi}\cdot \e,_{nj}+ \e,_{in}\cdot \e,_{jm})\e,_{k}],_l \bigr\|_{L^2(\os;\R^3)}(t),
\end{align*}
and by the same type of arguments as for $P_1(t)$,
\begin{align*}
P_2(t)\le  N(u_0,(w_i)_{i=1}^3)+ C t\ \|\w\|^3_{W_t},
\end{align*}
implying
\begin{equation}
\label{P}
P(t)\le  N(u_0,(w_i)_{i=1}^3)+ C t\ \|\w\|^3_{W_t}.
\end{equation}

Now, from the definition of a solution to the smoothed problem (\ref{smoothie}), 
\begin{align*}
\d \bigl\|\  \kappa c^{ijkl} (\e_{t},_{ilmn}\cdot e_j)\ e_k+ 
[c^{ijkl} (\e,_i\cdot \e_j-\delta_{ij})&\e,_k],_{mnl} \bigr\|_{L^2(\os;\R^3)}(t)\n\\
& =\|(\w_{t}-f-\kappa h),_{mn}(t)\|_{L^2(\os;\R^3)} ,
\end{align*} which implies with (\ref{P})
\begin{align*}
\|\frac{\kappa}{2} L(\w,_{nm})+  L(\e,_{nm})\|_{L^2(\os;\R^3)}(t)\le &\ \|(\w_{t}-f-\kappa h),_{nm}(t)\|_{L^2(\os;\R^3)}\n\\
&+ N(u_0,(w_i)_{i=1}^3)+ C t\ \|\w\|^3_{W_t}.
\end{align*}
Since this inequality also holds for any $t'\in(0,t)$, Lemma \ref{key} provides
\begin{align*}
\| L(\e,_{nm})\|_{L^{\infty}(0,t;L^2(\os;\R^3))}\le &\ C \|\w_{t},_{mn}\|_{L^{\infty}(0,t;L^2(\os;\R^3))}+C M(f,\kappa g,\kappa h)\n\\
&+ N(u_0,(w_i)_{i=1}^3)+ C t\ \|\w\|^3_{W_t},
\end{align*}
which with the estimate on $w_t$ from the previous subsection leads to  
\begin{align}
\| L(\e,_{nm})\|_{L^{\infty}(0,t;L^2(\os;\R^3))}\le & C_{\delta} [N(u_0,(w_i)_{i=1}^3)+M(f,\kappa g,\kappa h)+N((q_i)_{i=0}^2)]\n\\
&\ \ +C_{\delta} t^{\frac{1}{4}} \|(\tilde w,\tilde q)\|^6_{Z_t} + C {\delta} \|(\tilde w,\tilde q)\|^2_{Z_t}. 
\label{etah4s}
\end{align}

\noindent Step 3. 
From the estimates on $L_1-L_{2}$, and similar estimates that we could get in the fluid as in \cite{CoSh2004}, but this time by replacing $C(M)$ by appropriate powers of $\|(\w,\q)\|_{Z_t}$, we then deduce that for all $t\in [0,\check {T}]$,
\begin{align*}
  \frac{1}{2}\  \|\z  W,_{\ao\at\att} (t)\|^2_{L^2(\R^3_+;{\mathbb R}^3)} &
+{\nu} \int_0^{t} (\z^2 b_k^r  W,_{\ao\at\att r},\  b_k^s W,_{\ao\at\att s})_{L^2(\R^3_+;{\mathbb R}^3)}\n\\
& \le C_{\delta} [N(u_0,(w_i)_{i=1}^3)+M(f,\kappa g,\kappa h)+N((q_i)_{i=0}^2)]\n\\
&\ \ \ +C_{\delta} t^{\frac{1}{4}} \|(\tilde w,\tilde q)\|^6_{Z_t} + C {\delta} \|(\tilde w,\tilde q)\|^2_{Z_t}.
\end{align*}
By the trace theorem, we then get
\begin{align*}
\int_0^t \|W\|^2_{H^{3.5}(S;\R^3)} \le & C_{\delta} [N(u_0,(w_i)_{i=1}^3)+M(f,\kappa g,\kappa h)+N((q_i)_{i=0}^2)]\n\\
&\ \ \ +C_{\delta} t^{\frac{1}{4}} \|(\tilde w,\tilde q)\|^6_{Z_t} + C {\delta} \|(\tilde w,\tilde q)\|^2_{Z_t} ,
\end{align*}
where
$S=\{(x_1,x_2,x_3)\in\R^3|\ |x_1|\le\frac{1}{2},\ |x_2|\le\frac{1}{2},\ x_3=0\}$.
By a finite covering argument, we then get
\begin{align}
\int_0^t \|\w\|^2_{H^{3.5}(\Gamma_0;\R^3)} \le &\  C_{\delta} [N(u_0,(w_i)_{i=1}^3)+M(f,\kappa g,\kappa h)+N((q_i)_{i=0}^2)]\n\\
&\ \ \ +C_{\delta} t^{\frac{1}{4}} \|(\tilde w,\tilde q)\|^6_{Z_t} + C {\delta} \|(\tilde w,\tilde q)\|^2_{Z_t}.
\label{tracew}
\end{align}

 From the estimate (\ref{wth3f}) on $\w_t$ and the trace estimate (\ref{tracew}), we infer in a way similar to \cite{CoSh2004} by elliptic regularity arguments that
\begin{align}
\|\w\|^2_{L^2(0,t;H^4(\of;\R^3))}&+ \|\q\|^2_{L^2(0,t;H^3(\of;\R))}\n\\
& \le C_{\delta} [N(u_0,(w_i)_{i=1}^3)+M(f,\kappa g,\kappa h)+N((q_i)_{i=0}^2)]\n\\
&\ \ \ +C_{\delta} t^{\frac{1}{4}} \|(\tilde w,\tilde q)\|^6_{Z_t} + C {\delta} \|(\tilde w,\tilde q)\|^2_{Z_t}. 
\label{wh4f}
\end{align}

Similarly, from (\ref{etah4s}), and the trace estimate (\ref{tracew}), elliptic regularity provides
\begin{align}
\|\e\|^2_{L^{\infty}(0,t;H^4(\os;\R^3))} \le &\  C_{\delta} [N(u_0,(w_i)_{i=1}^3)+M(f,\kappa g,\kappa h)+N((q_i)_{i=0}^2)]\n\\
&\ +C_{\delta} t^{\frac{1}{4}} \|(\tilde w,\tilde q)\|^6_{Z_t} + C {\delta}\ \|(\tilde w,\tilde q)\|^2_{Z_t}. 
\label{etah4sbis}
\end{align}

\section{Time of existence independent of $\kappa$} 
\label{9}
From (\ref{I7}), (\ref{I8}), (\ref{wtth2}), (\ref{wth3f}), (\ref{wh4f}) and (\ref{etah4sbis}), we then have for any $t\in[0,T_{\kappa}]$, 
\begin{align*}
\|(\w,\q)\|^2_{Z_{t}}  \le &\  C_{\delta} [N_0(u_0,(w_i)_{i=1}^3)+M_0(f,\kappa g,\kappa h))+N_0((q_i)_{i=0}^2)]\n\\
&\ +C_{\delta} t^{\frac{1}{4}} \|(\tilde w,\tilde q)\|^6_{Z_t} + C_0 {\delta}\ \|(\tilde w,\tilde q)\|^2_{Z_t}.
\end{align*}
The subscripts $0$ in $C_0$, $N_0$, $M_0$ mean that we no longer consider generic constants from now on.

Now, let $\delta_0>0$ be such that $C_0\delta_0=\frac{1}{2}$. For $\kappa>0$ small enough, for any $t\in (0, T_{\kappa})$ we have
\begin{align}
\|(\w,\q)\|^2_{Z_{t}}  \le &\  4 C_{\delta_0} [N_0(u_0,(w_i)_{i=1}^3)+M_0(f)+N_0((q_i)_{i=0}^2)] +2 C_{\delta_0} t^{\frac{1}{4}}\ \|(\tilde w,\tilde q)\|^6_{Z_T},
\label{t1} 
\end{align}
where $M_0(f)=M_0(f,0,0)$. 
For conciseness, we will denote $C_1= 2C_{\delta_0}$ and $N_1=4 C_{\delta_0} [N_0(u_0,(w_i)_{i=1}^3)+M_0(f)+N_0((q_i)_{i=0}^2)]$.

 Now for $t\in (0, T_{\kappa})$ fixed, let $\d\alpha_t(x)=x^3-\frac{x}{C_1 t^{\frac{1}{4}}} + \frac{N_1}{C_1 t^{\frac{1}{4}}}$, so that $$\alpha_t(\|(\w,\q)\|^2_{Z_t})\ge 0.$$ Now let $t_1=[\frac{2}{27 C_1 N_1^2}]^4>0$, {\it which does not depend on} $\kappa$, and let $\check{T}=\min(T_{\kappa}, t_1)$. From now on, we assume that $t\in (0,\check T)$. We then have $\alpha_t((3C_1 t^{\frac{1}{4}})^{-\frac{1}{2}})<0$ which implies that $\alpha_t$ has three real roots $z_1$, $z_2$, $z_3$, with $z_1<-(3C_1 t^{\frac{1}{4}})^{-\frac{1}{2}}<z_2<(3C_1 t^{\frac{1}{4}})^{-\frac{1}{2}}<z_3$. From the product $\d z_1 z_2 z_3=-\frac{N_1}{C_1 t^{\frac{1}{4}}}$ and $\alpha_t(3 N_1)<0$, we infer that $0<z_2< 3 N_1<z_3$. From (\ref{t1}) and the {\it continuity} of $t\rightarrow \|(\w,\q)\|_{Z_t}$ (established in Lemma \ref{continuity}) we then infer since $\|(\w,\q)\|^2_{Z_0}\le N_1<z_3$ that we have
\begin{equation}
\label{t2}
\forall t\in (0,\check T],\ \|(\w,\q)\|^2_{Z_t}\le z_2\le 3 N_1\ .
\end{equation}
This implies that $\e(\check T)\in H^4(\of;\R^3)\cap H^4(\os;\R^3)$, $\w(\check T)\in H^1_0(\Omega;\R^3)\cap H^3(\os;\R^3)\cap H^3(\os;\R^3)$, $\w_t(\check T)\in H^1_0(\Omega;\R^3)\cap H^2(\os;\R^3)\cap H^2(\os;\R^3)$, $\w_{tt}(\check T)\in H^1_0(\Omega;\R^3)$, $\w_{ttt}(\check T)\in L^2(\Omega;\R^3)$, $\q(\check T)\in H^2(\Omega_0^f;\R)$, $\q_t(\check T)\in H^1(\Omega_0^f;\R)$, $\q_{tt}(\check T)\in L^2(\Omega_0^f;\R)$, with a bound that depends only on the right-hand side of (\ref{t2}). The compatibility conditions for the smoothed problem (\ref{smoothie}) at $\check T$ are also satisfied by definition of a solution, which means that we do not have any new term of the type of $b_{\kappa}$, $c_{\kappa}$ or $d_{\kappa}$ associated to $\w(\check{T})$ to add to the already existing forcing terms coming from $t=0$. 

We can thus build a solution of the smoothed problem (\ref{smoothie}) defined on $[\check T, \check T+\delta T]$, $\delta T$ depending solely on the right-hand side of (\ref{t2}), that we will still denote $(\w,\q)$. It is then readily seen that $(\w,\q)\in Z_{\check T+\delta T}$ and is a solution of the approximated problem (\ref{smoothie}) on $[0,\check T+\delta T]$. If $\check T=t_1$, we have our solution defined on the $\kappa$ independent time interval $[0,t_1]$, with the $\kappa$ independent estimate (\ref{t2}). Otherwise, if $\check T<t_1$, we can also assume that $\check T+\delta T\le t_1$, which implies, in the same fashion as we got (\ref{t1}),
\begin{equation}
\label{t3}
\forall t\in [0,\tilde T+\delta T],\ \|(\w,\q)\|^2_{Z_{t}}\le  N_1 +C_1 t^{\frac{1}{4}} \ \|(\w,\q)\|^6_{Z_{t}}.
\end{equation}
This implies in turn that $\e(\check T+\delta T)$, $\w(\check T+\delta T), \w_t(\check T+\delta T), \w_{tt}(\check T+\delta T),  \w_{ttt}(\check T+\delta T) $, $\q(\check T+\delta T), \q_t(\check T+\delta T), \q_{tt}(\check T+\delta T)$ are in the same spaces as their respective
counterparts at time $\check T$, with the same bound as well, since we could from (\ref{t3}) repeat the same argument leading to (\ref{t2}), this time on $[0,\check T+\delta T]$. The compatibility conditions at $\check T+\delta T$ being also automatically satisfied, we can thus build a solution of the approximated problem (\ref{smoothie}) defined on $[\check T+\delta T, \check T+2 \delta T]$, the time of existence being the same as starting from $\check T$ from the similarity of the bound that we obtain
on $\e(\check T+\delta T), \partial_t^n \w(\check T+\delta T) (n=0,1,2,3), \partial_t^n \q(\check T+\delta T) (n=0,1,2)$ and their respective counterparts at time $\check T$. We will still denote this solution $(\w,\q)$. It is then readily seen that $(\w,\q)\in Z_{\check T+2 \delta T}$ and is a solution of the approximated problem on $[0,\check T+2\delta T]$. We then have in the same fashion as we got (\ref{t1}), 
\begin{equation*}
\forall t\in [0,\tilde T+2 \delta T],\ \|(\w,\q)\|^2_{Z_{t}}\le  N_1 +C_1 t^{\frac{1}{4}} \ \|(\w,\q)\|^6_{Z_{t}}.
\end{equation*}
By induction, we then see that we get a solution $(\w,\q)$ defined on $[0,t_1]$,satisfying the estimate
\begin{equation}
\label{timeindep}
\forall t\in [0,t_1],\ \|(\w,\q)\|_{Z_{t}}\le 3 N_1= 12 C_{\delta_0} [N_0(u_0,(w_i)_{i=1}^3)+M_0(f)+N_0((q_i)_{i=0}^2)],
\end{equation}
 establishing the independence of the time of existence respectively to $\kappa$, since $t_1$ does not depend on $\kappa$.
In the following we will note $T=t_1$.

\section{Existence for (\ref{nsl})}
\label{10}
\begin{proof}
We can here choose to take $\kappa=\frac{1}{n}$,  and let $n\rightarrow\infty$. By the bound (\ref{timeindep}) independent of $\kappa$ on $[0,T]$, we then have the existence of a weakly convergent subsequence of $(\w,\q)$ in the reflexive Hilbert space $Y_{T}$, to a limit that we
call $(v,q)$, which also belongs to $Z_{T}$ and satisfies the estimate
\begin{equation*}
\ \|(v,q)\|_{Z_{T}}\le 3 N_1= 12 C_{\delta_0} [N_0(u_0,(w_i)_{i=1}^3)+M_0(f)+N((q_i)_{i=0}^2)].
\end{equation*}
The usual compactness theorems ensure at this stage that $(v,q)$ is a solution of
(\ref{nsl}) on $[0,T]$. The smoothness of our solution ensures that the solids do not collide with each other (if there is more than one) or the boundary (for an eventually smaller time), which establishes the existence part of Theorem.
\ref{main}
\end{proof}

\section {Uniqueness for (\ref{nsl})}
\label{11}

\begin{proof}
Since we cannot use a contractive mapping scheme for our problem, we have
to establish uniqueness separately. Let then $(\bar{v},\bar{q})$ denote another solution of (\ref{nsl}) in $Z_T$. Then, taking $v-\bar{v}$ as a test function in the variational formulation of the difference between the systems (\ref{nsl}) associated to each solution provides for $t\in [0,T]$:
\begin{align}
&  \frac{1}{2}\  \|(v-\bar{v}) (t)\|^2_{L^2(\Omega;{\mathbb R}^3)} 
+{\nu} \int_0^{t} (a_k^r a_k^s v,_r-{\bar{a}_k}^r {\bar{a}_k}^s \bar{v},_r,\   {v},_s-\bar{v},_s)_{L^2(\Omega_0^f;{\mathbb R}^3)} 
\nonumber\\
& +\int_0^t (c^{ijkl} [(\eta,_i\cdot\eta,_j-\delta_{ij}) \eta,_k - (\bar{\eta},_i\cdot\bar{\eta},_j-\delta_{ij}) \bar{\eta},_k],\ v,_l-\bar{v},_l)_{L^2(\os;\R^3)}\nonumber\\
& -\int_0^t ( a_i^j  q - {\bar a}_i^j \bar{q}, v^i,_j-\bar{v}^i,_j)_{L^2(\Omega_0^f;{\mathbb R})}
= \int_0^{t}( f\circ\eta-f\circ\bar{\eta}, v-\bar{v})_{L^2(\Omega_0^f;{\mathbb R}^3)} 
.
\label{u0}
\end{align}
For the viscous term in the fluid, we write $$a_k^r a_k^s v,_r-{\bar{a}_k}^r {\bar{a}_k}^s \bar{v},_r=a_k^r a_k^s (v,_r- \bar{v},_r)+ (a_k^r a_k^s-{\bar{a}_k}^r {\bar{a}_k}^s) \bar{v},_r,$$
which with the $L^{\infty}(0,T;H^3(\of;\R^3))$ control of $\bar{v}$ and $v$ provides us with an estimate of the type (where $C$ denotes once again a generic constant)
\begin{align}
\label{u1}
\int_0^{t} (a_k^r a_k^s v,_r-{\bar{a}_k}^r {\bar{a}_k}^s \bar{v},_r,\   {v},_s-\bar{v},_s)_{L^2(\Omega_0^f;{\mathbb R}^3)}\ge &\ C \int_0^t \|v-\bar v\|^2_{H^1(\Omega_0^f;{\mathbb R}^3)}\n\\
& -C\int_0^t\int_0^{t'}\|v-\bar v\|^2_{H^1(\Omega_0^f;{\mathbb R}^3)}.
\end{align}  
Concerning the forcing term in the fluid, we first notice that if we still denote $E(\Omega)(f)$ as $f$, 
$$\d f(t,\bar{\eta}(t,x))-f(t,{\eta}(t,x))=\int_0^1 f,_i (t,(\eta+t'(\bar{\eta}-\eta))(t,x)) dt'\ (\bar{\eta}^i(t,x)-\eta^i(t,x)),$$
leading us to
\begin{align*}
\| f(t,\bar{\eta}(t,\cdot))&-f(t,{\eta}(t,\cdot))\|_{L^{1.5}(\of;\R^3)}\\
& \le C
\|\bar{\eta}(t,\cdot)-\eta(t,\cdot)\|_{L^6(\of;\R^3)} \bigl[ \sum_{i=1}^3\int_0^1 \int_{\of} f,_i^2 (t, \phi(t',t,x)) dx dt'\bigr]^{0.5},
\end{align*}
with $\phi(t',t,x)=\eta(t,x)+t'(\bar{\eta}(t,x)-\eta(t,x))$. We have $\phi(t',t,\cdot)\in \mathcal C^0(\overline{\Omega};\R^3)\cap \mathcal C^1 (\Omega\cap\Gamma_0^c;\R^3)$. Moreover $\phi(t',t,\partial\Omega)=\partial\Omega$. We then have
by invariance by homotopy of the Brouwer degree (for the parameter $t$) 
$$\forall z\in \Omega,\ \text{deg}(\phi(t',t,\cdot),\Omega,z)=
\text{deg}(\phi(t',0,\cdot),\Omega,z)=\text{deg}(\text{Id},\Omega,z)=1,$$
which together with the regularity of $\phi(t',t,\cdot)$ establishes that $\phi(t',t,\cdot)(\Omega)=\Omega$ and that $\text{Card}\{\phi^{-1}(t',t,\cdot)(x)\}=1$ for
almost all $x\in\Omega$. Thus, 
$$\int_{\of} f,_i^2 (t, \phi(t',t,x)) dx =\int_{\phi(t',t,\of)} f,_i^2 (t,y) |\text{det} \nabla\phi(t',t, \phi^{-1}(t',t,y))|^{-1}\ dy,$$
which with the $L^{\infty}(0,T;H^4(\of;\R^3))$ control of $\eta$ and $\bar{\eta}$ provides
$$\int_{\of} f,_i^2 (t, \phi(t',t,x)) dx \le C \int_{\Omega} f,_i^2 (t,y) dy.$$
Consequently,
\begin{align*}
\| f(t,\bar{\eta}(t,\cdot))&-f(t,{\eta}(t,\cdot))\|_{L^{1.5}(\of;\R^3)}\le C
\|\bar{\eta}(t,\cdot)-\eta(t,\cdot)\|_{H^1(\of;\R^3)} \|f\|_{H^1(\Omega;\R^3)},
\end{align*}
implying
\begin{align}
\label{u2}
\bigl|\int_0^{t}( f\circ\eta-f\circ\bar{\eta}, v-\bar{v})_{L^2(\Omega_0^f;{\mathbb R}^3)}\bigr|\le  C \sqrt{t} \|f\|_{L^2(0,t;H^1(\Omega;\R^3))}
   \|v-\bar{v}\|^2_{L^2(0,t;H^1(\of;\R^3))} .
\end{align}
Concerning the elastic term,
\begin{align*}
\int_0^t (c^{ijkl} [(\eta,_i\cdot\eta,_j&-\delta_{ij}) \eta,_k - (\bar{\eta},_i\cdot\bar{\eta},_j-\delta_{ij}) \bar{\eta},_k],\ v,_l-\bar{v},_l)_{L^2(\os;\R^3)}=I_1+I_2+I_3,
\end{align*}
with
\begin{align*}
 I_1 &=\int_0^t (c^{ijkl} (\eta,_i\cdot\eta,_j-\delta_{ij}) (\eta,_k - \bar{\eta},_k) ,\ v,_l-\bar{v},_l)_{L^2(\os;\R^3)}\\
&= \frac{1}{2}(c^{ijkl} (\eta,_i\cdot\eta,_j-\delta_{ij}) (\eta,_k - \bar{\eta},_k) ,\ \eta,_l-\bar{\eta},_l)_{L^2(\os;\R^3)}(t)\\
&\ \ -\frac{1}{2}\int_0^t (c^{ijkl} (\eta,_i\cdot\eta,_j-\delta_{ij})_t (\eta,_k - \bar{\eta},_k) ,\ \eta,_l-\bar{\eta},_l)_{L^2(\os;\R^3)} \\
&\ge - C t \|\eta(t)-\bar{\eta}(t)\|^2_{H^1(\os;\R^3)} -C \int_0^t \|\eta-\bar{\eta}\|^2_{H^1(\os;\R^3)},  
\end{align*}
where we have used the $L^{\infty}(0,T;H^3(\os;\R^3))$ control of $v$ and
$\bar{v}$ for the inequality. 

Next, for the same reasons,
\begin{align*}
 I_2 &=\int_0^t (c^{ijkl} (\eta,_i-\bar{\eta},_i)\cdot\eta,_j \ \bar{\eta},_k ,\ v,_l-\bar{v},_l)_{L^2(\os;\R^3)}\\
&= \int_0^t (c^{ijkl} (\eta,_i-\bar{\eta},_i)\cdot(\eta,_j-\bar{\eta},_j) \ \bar{\eta},_k ,\ v,_l-\bar{v},_l)_{L^2(\os;\R^3)}\\
&\ \ + \int_0^t (c^{ijkl} (\eta,_i-\bar{\eta},_i)\cdot\bar{\eta},_j \ \bar{\eta},_k ,\ v,_l-\bar{v},_l)_{L^2(\os;\R^3)}\\
&= \int_0^t (c^{ijkl} (\eta,_i-\bar{\eta},_i)\cdot(\eta,_j-\bar{\eta},_j) \ \bar{\eta},_k ,\ v,_l-\bar{v},_l)_{L^2(\os;\R^3)}\\
&\ \ + \frac{1}{2}(c^{ijkl} (\eta,_i-\bar{\eta},_i)\cdot\bar{\eta},_j \ \bar{\eta},_k ,\ \eta,_l-\bar{\eta},_l)_{L^2(\os;\R^3)}(t)\\
&\ \ - \int_0^t (c^{ijkl} (\eta,_i-\bar{\eta},_i)\cdot\bar{\eta},_j \ \bar{v},_k ,\ \eta,_l-\bar{\eta},_l)_{L^2(\os;\R^3)}.
\end{align*}
We then write for the second term on the right-hand side of the last equality 
$\d \bar \eta,_i(t,\cdot)=e_i+\int_0^t \bar{v},_i$, to get by Korn's inequality 
 \begin{align*}
I_2 \ge  C [ \|\eta(t)-\bar{\eta}(t)\|^2_{H^1(\os;\R^3)} - \|\eta(t)-\bar{\eta}(t)\|^2_{L^2(\os;\R^3)}] -C t \sup_{[0,t]}\|\eta-\bar{\eta}\|^2_{H^1(\os;\R^3)}.  
\end{align*}
Similarly,
\begin{align*}
 I_3 &=\int_0^t (c^{ijkl} (\eta,_j-\bar{\eta},_j)\cdot\bar{\eta},_i \ \bar{\eta},_k ,\ v,_l-\bar{v},_l)_{L^2(\os;\R^3)}\\
&\ge  C [\ \|\eta(t)-\bar{\eta}(t)\|^2_{H^1(\os;\R^3)}-\|\eta(t)-\bar{\eta}(t)\|^2_{L^2(\os;\R^3)}] -C \int_0^t \|\eta-\bar{\eta}\|^2_{H^1(\os;\R^3)}.  
\end{align*}
Thus, 
\begin{align}
\int_0^t& (c^{ijkl} [(\eta,_i\cdot\eta,_j-\delta_{ij}) \eta,_k - (\bar{\eta},_i\cdot\bar{\eta},_j-\delta_{ij}) \bar{\eta},_k],\ v,_l-\bar{v},_l)_{L^2(\os;\R^3)}\n\\
&\ge  C [\ \|\eta(t)-\bar{\eta}(t)\|^2_{H^1(\os;\R^3)}-\|\eta(t)-\bar{\eta}(t)\|^2_{L^2(\os;\R^3)}] -C \int_0^t \|\eta-\bar{\eta}\|^2_{H^1(\os;\R^3)}.
\label{u3}
\end{align}
Concerning the pressure term, with $a_i^j  q - {\bar a}_i^j \bar{q}=(a_i^j-\bar{a}_i^j) q+ \bar{a}_i^j (q-\bar{q})$ and the $L^{\infty}(0,T;H^2(\of;\R))$ control of the pressure, we get
\begin{align}
&-\int_0^t ( a_i^j  q - {\bar a}_i^j \bar{q}, v^i,_j-\bar{v}^i,_j)_{L^2(\Omega_0^f;{\mathbb R})}\n\\
&\ge -C [ \sqrt{t}\|q-\bar{q}\|_{L^{\infty}(0,t;L^2(\of;\R))} \|v-\bar{v}\|_{L^2(0,t;H^1(\of;\R^3))}+t \|v-\bar{v}\|^2_{L^2(0,t;H^1(\of;\R^3))}].
\label{u4}  
\end{align}

In order to get the estimate of $q-\bar{q}$ in $L^2(\of;\R)$, we have to introduce the time differentiated problem. By taking $v_t-\bar{v}_t$ in the variational formulation associated to the difference between the time differentiated systems, we obtain
\begin{align}
&  \frac{1}{2}\  \|(v_t-\bar{v}_t) (t)\|^2_{L^2(\Omega;{\mathbb R}^3)} 
+{\nu} \int_0^{t} ([ a_k^r a_k^s v,_r-{\bar{a}_k}^r {\bar{a}_k}^s \bar{v},_r]_t,\   [{v},_s-\bar{v},_s]_t)_{L^2(\Omega_0^f;{\mathbb R}^3)} 
\nonumber\\
& +\int_0^t (c^{ijkl} [(\eta,_i\cdot\eta,_j-\delta_{ij}) \eta,_k - (\bar{\eta},_i\cdot\bar{\eta},_j-\delta_{ij}) \bar{\eta},_k]_t,\ [v,_l-\bar{v},_l]_t)_{L^2(\os;\R^3)}\nonumber\\
& -\int_0^t ( [ a_i^j  q - {\bar a}_i^j \bar{q}]_t, [v^i,_j-\bar{v}^i,_j]_t)_{L^2(\Omega_0^f;{\mathbb R})}
= \int_0^{t}( [ f\circ\eta-f\circ\bar{\eta}]_t, v_t-\bar{v_t})_{L^2(\Omega_0^f;{\mathbb R}^3)} 
.
\label{u5}
\end{align}
For the fluid viscous term, we easily find
 with the $L^2(0,T;H^3(\of;\R^3))$ control of the first time derivative of the velocity that 
\begin{equation}
\label{u6}
\int_0^{t} ([a_k^r a_k^s v,_r-{\bar{a}_k}^r {\bar{a}_k}^s \bar{v},_r]_t,\   {v_t},_s-\bar{v_t},_s)_{L^2(\Omega_0^f;{\mathbb R}^3)}\ge \ C (1-t) \int_0^t \|v_t-\bar v_t\|^2_{H^1(\Omega_0^f;{\mathbb R}^3)}.
\end{equation}  
Concerning the forcing term in the fluid, since
$(f\circ\eta)_t=(f_t+v^i f,_i)(\eta)$ (with a similar formula for $\bar{v}$), we
then deduce in a way similar to the steps leading to (\ref{u2}) that
\begin{align}
\label{u7}
\bigl|\int_0^{t}( [f&\circ\eta -f\circ\bar{\eta}]_t, [v-\bar{v}]_t)_{L^2(\Omega_0^f;{\mathbb R}^3)}\bigr|\n\\
&\le  C \sqrt{t}\ [ \|f_t\|_{L^2(0,t;H^1(\Omega;\R^3))}+ \|f\|_{L^2(0,t;H^2(\Omega;\R^3))}]\ \|v_t-\bar{v_t}\|^2_{L^2(0,t;H^1(\of;\R^3))} .
\end{align}
For the elastic term, we can also essentially reproduce the arguments leading to (\ref{u3}), leading us to
\begin{align}
\int_0^t & (c^{ijkl} [(\eta,_i\cdot\eta,_j-\delta_{ij}) \eta,_k - (\bar{\eta},_i\cdot\bar{\eta},_j-\delta_{ij}) \bar{\eta},_k]_t,\ [v,_l-\bar{v},_l]_t)_{L^2(\os;\R^3)}\n\\
&\ge C [\ \| v(t)-\bar{v}(t)\|^2_{H^1(\os;\R^3)}-\| v(t)-\bar{v}(t)\|^2_{L^2(\os;\R^3)}] -C t \sup_{[0,t]} \|v-\bar{v}\|^2_{H^1(\os;\R^3)}.
\label{u8}
\end{align}
The pressure term will require more care since we want to avoid the introduction 
of $q_t-\bar{q}_t$ that the most direct method would lead to. To do so, we 
notice that
\begin{align*}
\int_0^t ( [a_i^j  q - {\bar a}_i^j \bar{q}]_t, [v^i,_j-\bar{v}^i,_j]_t)_{L^2(\Omega_0^f;{\mathbb R})}=I_4+I_5+I_6, 
\end{align*}
with 
\begin{align*}
I_4=& \int_0^t ( [(a_i^j)_t  q - ({\bar a}_i^j)_t \bar{q}], [v^i,_j-\bar{v}^i,_j]_t)_{L^2(\Omega_0^f;{\mathbb R})},\\
I_5=& \int_0^t ( a_i^j ( q_t - \bar{q_t}), [v^i,_j-\bar{v}^i,_j]_t)_{L^2(\Omega_0^f;{\mathbb R})},\\ 
I_6=& \int_0^t ( [a_i^j - {\bar a}_i^j ] \bar{q_t}, [v^i,_j-\bar{v}^i,_j]_t)_{L^2(\Omega_0^f;{\mathbb R})}.
\end{align*}
For $I_4$, we have in a way similar to (\ref{u4}),
\begin{align*}
|I_4| \le C [  \sqrt{t} \|q-\bar{q}\|_{L^{\infty}(0,t;L^2(\of;\R))}& \|v_t-\bar{v_t}\|_{L^2(0,t;H^1(\of;\R^3))}\\
& +t \|v_t-\bar{v_t}\|^2_{L^2(0,t;H^1(\of;\R^3))}]. 
\end{align*}
For $I_6$, the $L^2(0,T;H^2(\of;\R))$ control of $\bar{q}_t$ provides us with
$$|I_6|\le C t\ \|v_t-\bar{v_t}\|^2_{L^2(0,t;H^1(\of;\R^3))}.$$
For $I_5$ we have:
\begin{align*}
I_5=& \int_0^t (  q_t - \bar{q_t}, a_i^j {v^i_t},_j -\bar{a}_i^j {\bar{v}^i_t},_j)_{L^2(\Omega_0^f;{\mathbb R})}- \int_0^t (  q_t - \bar{q_t}, (a_i^j  -\bar{a}_i^j) {\bar{v}^i_t},_j)_{L^2(\Omega_0^f;{\mathbb R})}\\
=& \int_0^t ( \bar{q_t}-q_t , (a_i^j)_t {v^i},_j -(\bar{a}_i^j)_t {\bar{v}^i},_j)_{L^2(\Omega_0^f;{\mathbb R})}- \int_0^t (  q_t - \bar{q}_t, (a_i^j  -\bar{a}_i^j) {\bar{v}^i_t},_j)_{L^2(\Omega_0^f;{\mathbb R})},
\end{align*}
where we have used the relations $a_i^j v^i,_j=0=\bar{a}_i^j \bar{v}^i,_j$ in
$\of$ for the first integral. By integrating by parts in time,
\begin{align*}
I_5=& \int_0^t ( {q}-\bar{q} , [(a_i^j)_t {v^i},_j -(\bar{a}_i^j)_t {\bar{v}^i},_j]_t )_{L^2(\Omega_0^f;{\mathbb R})}+ \int_0^t (  q - \bar{q}, [(a_i^j  -\bar{a}_i^j) {\bar{v}^i_t},_j]_t)_{L^2(\Omega_0^f;{\mathbb R})}\\
&+( \bar{q}-q , (a_i^j)_t {v^i},_j -(\bar{a}_i^j)_t {\bar{v}^i},_j )_{L^2(\Omega_0^f;{\mathbb R})}(t) + (  \bar{q}-q , (a_i^j  -\bar{a}_i^j) {\bar{v}^i_t},_j)_{L^2(\Omega_0^f;{\mathbb R})}(t) .
\end{align*}
With the $L^2(0,T;H^3(\of;\R^3))$ control of $v_t$ we have
\begin{align*}
&\bigl|\int_0^t ( {q}-\bar{q} , [(a_i^j)_t {v^i},_j -(\bar{a}_i^j)_t {\bar{v}^i},_j]_t )_{L^2(\Omega_0^f;{\mathbb R})}\bigr|+\bigl|\int_0^t ( {q}-\bar{q} , [(a_i^j)_t -(\bar{a}_i^j)_t] {\bar{v}^i_t},_j )_{L^2(\Omega_0^f;{\mathbb R})}\bigr|\\
&\qquad\qquad\qquad\qquad\qquad\qquad\le C \sqrt{t}\ \|q-\bar{q}\|_{L^{\infty}(0,t;L^2(\of;\R))} \|v_t-\bar{v_t}\|_{L^2(0,t;H^1(\of;\R^3))},\\
&\bigl| ( {q}-\bar{q} , (a_i^j)_t {v^i},_j -(\bar{a}_i^j)_t {\bar{v}^i},_j )_{L^2(\Omega_0^f;{\mathbb R})}(t)\bigr|\\
&\qquad\qquad\qquad\qquad\qquad\qquad \le C \sqrt{t}\ \|q(t)-\bar{q}(t)\|_{L^2(\of;\R)} \|v_t-\bar{v_t}\|_{L^2(0,t;H^1(\of;\R^3))}.
\end{align*} 
The remaining terms are more delicate. We first have
\begin{align}
&\bigl|\int_0^t ( {q}-\bar{q} , (a_i^j -\bar{a}_i^j)\ {\bar{v}^i_{tt}},_j )_{L^2(\Omega_0^f;{\mathbb R})}\bigr|+\bigl|( {q}-\bar{q} , (a_i^j -\bar{a}_i^j)\ {\bar{v}^i_{t}},_j )_{L^2(\Omega_0^f;{\mathbb R})}(t) \bigr|\n\\
&\qquad\qquad\qquad\le C \int_0^t \|q-\bar{q}\|_{L^2(\of;\R)} \|a-\bar{a}\|_{L^{4}(\of;\R^9)} \|\nabla\bar{v}_{tt}\|_{L^4(\of;\R^9)}\n\\
& \qquad\qquad\qquad\qquad+ \|q(t)-\bar{q}(t)\|_{L^2(\of;\R)} \|a(t)-\bar{a}(t)\|_{L^4(\of;\R^9)} \|\nabla\bar{v}_{t}(t)\|_{L^4(\of;\R^9)}.
\label{u9}
\end{align}
The apparent problem here is that $a-\bar{a}$ is estimated
in $L^2(\of;\R^9)$ in terms of $v-\bar{v}$ in $H^1(\of;\R^3)$. Now, a bound of this quantity in $L^4(\of;\R^9)$ will require a bound of $v-\bar{v}$ in $H^2(\of;\R^3)$. In order to get such an estimate, we will bound $v-\bar{v}$ in $H^2(\of;\R^3)$ by lower order terms in $v-\bar{v}$. To do so, let us first estimate the trace of $v-\bar{v}$ on $\Gamma_0$ by using the test function $-[\z^2\ (v-\bar{v})\circ\P],_{\alpha\alpha}\circ\P^{-1}$ in the difference between the variational problems satisfied by
$v$ and $\bar{v}$. By proceeding as in Section \ref{10}, we would then get an estimate of the type, where $\delta>0$ is given:
\begin{align*}
&\int_0^t \|\z \nabla[(v-\bar{v})\circ\P],_{\alpha}\|^2_{L^2(\R^3_+;\R^9)}
+ \|\z \nabla[(\eta-\bar{\eta})\circ\P],_{\alpha}(t) \|^2_{L^2(\R^3_-;\R^9)}\\
&\le
C[\sqrt{t}+\delta]\ \int_0^t [ \|v-\bar{v}\|^2_{H^2(\of;\R^3)}+\|q-\bar{q}\|^2_{H^1(\of;\R)}] + C_{\delta} \int_0^t \|v_t-\bar{v}_t\|^2_{L^2(\Omega;\R^3)}\\
&\ \ + C \int_0^t \|\eta-\bar{\eta}\|^2_{H^2(\os;\R^3)} + C \int_0^t \|v-\bar{v}\|^2_{H^1(\Omega;\R^3)}
\end{align*}
which by patching all the charts defining $\Gamma_0$ leads to an estimate
of $v-\bar{v}$ in $L^2(0,t;H^{1.5}(\Gamma_0;\R^3))$ yielding by elliptic
regularity:
\begin{align*}
&\int_0^t [ \|v-\bar{v}\|^2_{H^2(\of;\R^3)}+\|q-\bar{q}\|^2_{H^1(\of;\R)}]+ \|\eta(t)-\bar{\eta}(t)\|^2_{H^2(\os;\R^3)}\n\\
&\le  \ 
C[\sqrt{t}+\delta] \int_0^t [ \|v-\bar{v}\|^2_{H^2(\of;\R^3)}+\|q-\bar{q}\|^2_{H^1(\of;\R)}] + C_{\delta} \int_0^t  \|v_t-\bar{v}_t\|^2_{L^2(\Omega;\R^3)} \\
&\ \ + C \int_0^t \|\eta-\bar{\eta}\|^2_{H^2(\os;\R^3)}+ C \int_0^t \|v-\bar{v}\|^2_{H^1(\Omega;\R^3)}
.
\end{align*}
Thus, with a choice of $\delta>0$ small enough, we have for $t$ small enough and the use of Gronwall's inequality,
\begin{align*}
\|\eta(t)-\bar{\eta}(t)\|^2_{H^2(\os;\R^3)}&+
\int_0^t [ \|v-\bar{v}\|^2_{H^2(\of;\R^3)}+\|q-\bar{q}\|^2_{H^1(\of;\R)}]
\\
& \le  C \int_0^t  \|v_t-\bar{v}_t\|^2_{L^2(\Omega;\R^3)}+ C \int_0^t \|v-\bar{v}\|^2_{H^1(\Omega;\R^3)} .
\end{align*}
By using this estimate in (\ref{u9}), we then get for a time small enough
\begin{align*}
&\bigl|\int_0^t ( {q}-\bar{q} , (a_i^j -\bar{a}_i^j)\ {\bar{v}^i_{tt}},_j )_{L^2(\Omega_0^f;{\mathbb R})}\bigr|+\bigl|( {q}-\bar{q} , (a_i^j -\bar{a}_i^j)\ {\bar{v}^i_{t}},_j )_{L^2(\Omega_0^f;{\mathbb R})}(t) \bigr|\n\\
&\le C \sqrt{t}\ [\int_0^t \|v_t-\bar{v}_t\|^2_{L^{2}(\Omega;\R^3)} + \int_0^t \|v-\bar{v}\|^2_{H^1(\Omega;\R^3)}+ \|q-\bar{q}\|^2_{L^{\infty}(0,t;L^2(\of;\R))}] .
\end{align*}
By putting together the estimates on $I_4$, $I_5$ and $I_6$, we have
\begin{align}
&\bigl|\int_0^t ( [a_i^j  q - {\bar a}_i^j \bar{q}]_t, [v^i,_j-\bar{v}^i,_j]_t)_{L^2(\Omega_0^f;{\mathbb R})}\bigr|\n\\
&\le C \sqrt{t}\ [\int_0^t \|v_t-\bar{v}_t\|^2_{L^{2}(\Omega;\R^3)} + \int_0^t \|v-\bar{v}\|^2_{H^1(\Omega;\R^3)}+ \|q-\bar{q}\|^2_{L^{\infty}(0,t;L^2(\of;\R))}] .
\label{u10}
\end{align}
Now, by considering the difference between the two variational forms satisfied
repectively by $(v,q)$ and $(\bar{v},\bar{q})$, and writing the difference between the pressure
terms as
\begin{align*}
\int_{\of} (a_i^j q-\bar{a}_i^j \bar{q}) \phi^i,_j=\int_{\of} a_i^j (q-\bar{q})  \phi^i,_j+\int_{\of} (a_i^j -\bar{a}_i^j) \bar{q} \phi^i,_j,
\end{align*}
the Lagrange multiplier Lemma 13 of \cite{CoSh2004} provides for all $t\in [0,T]$:
\begin{align}
\|q(t)-\bar{q}(t)\|_{L^2(\of;\R)}\le C [\ & \|(v_t-\bar{v}_t)(t)\|_{L^2(\Omega;\R^3)} + \|(v-\bar{v})(t)\|_{H^1(\Omega_0^f;\R^3)}\n\\
& + \|(\eta-\bar{\eta})(t)\|_{H^1(\os;\R^3)}+ \sqrt{t} \|v-\bar{v}\|_{L^2(0,t;H^1(\of;\R^3))}].
\label{u11} 
\end{align}
By putting together the estimates (\ref{u0})-(\ref{u11}), we then obtain for
$t_u>0$ small enough an inequality of the type:
\begin{align*}
\|v_t-\bar{v}_t\|&^2_{L^{\infty}(0,t_u;L^2(\Omega;\R^3))} + \int_0^{t_u} \|v_t-\bar{v}_t\|^2_{H^1(\of;\R^3)}+ \|v-\bar{v}\|^2_{L^{\infty}(0,t_u;H^1(\Omega_0^s;\R^3))}\le 0, 
\end{align*}
which shows that $(v,q)=(\bar{v},\bar{q})$ on $[0,t_u]$. Let 
$$T_u=\sup\{t\in [0,T]|\ (v,q)=(\bar{v},\bar{q})\ \text{on}\ [0,t]\}.$$
If $T_u<T$, we can repeat the same procedure with $T_u$ replacing $0$, which would lead to uniqueness for $[T_u,T_u+\delta t)$ as well. Thus, we have $T_u=T$, which concludes the proof of the theorem.
\end{proof}      

\section {Optimal regularity on the initial data}
\label{12}
We  first remind some extensions and regularization results on domains:
%\begin{lemma}
%\label{lift}
%Let $\Omega'$ be a domain of class of class $H^k$ ($k\ge 4$).\hfill\break
%Then, 
%there exists a linear and continuous operator $L(\Omega')$ from $\{(d,r)\in H^m(\Omega';{\mathbb R})\times H^{m+0.5}(\partial\Omega';{\mathbb R}^3)|\ \displaystyle \int_{\Omega} d=\int_{\partial\Omega} r\cdot N \}$ into $H^{m+1}(\Omega';{\mathbb R}^3)$ (for $0\le m\le k-2$), which to $(d,r)\in H^m(\Omega';{\mathbb R})\times H^{m+0.5}(\partial\Omega';{\mathbb R}^3)$ satisfying $\displaystyle\int_{\Omega'} d=\int_{\partial\Omega'} r\cdot N$ associates $u=L(\Omega')(d,r) \in H^{m+1}(\Omega';{\mathbb R}^3)$ satisfying 
%\begin{align*}
%\operatorname{div}u &=d\ \text{in}\ \Omega',\\
%u&=r\ \text{on}\ \partial\Omega'\ .
%\end{align*}
%Furthermore, the norms of this operator between those spaces stay bounded as
%the norm of the charts defining $\Omega'$ stay in a bounded set of $H^k$. 
%\end{lemma}

%We will need the following extension Lemma:
\begin{lemma}
\label{extension}
Let $\Omega'$ be a domain of class $H^4$. Then,  
there exists a linear and continuous operator $E(\Omega')$ from $H^m(\Omega';{\mathbb R}^3)$ into $H^m(\R^3;{\mathbb R}^3)$ (for each $0\le m\le 4$) such that $E(\Omega')(u)=u$ in $\Omega'$. Also, if the $H^4$ norms of a family of domains
stay bounded, the norms of the corresponding linear operators also stay bounded. 
\end{lemma}

\begin{lemma}
\label{domain}
Since $\Omega_0^s$ is of class $H^4$, let $\psi^m\in H^4(B_-(0,1);\R^3)$ ($m=1,..., N$) be a collection of charts defining a neighborhood of its boundary. We note 
$$\|\of\|_{H^4}=\sum_{m=1}^N \|\psi^m\|_{H^4(B_-(0,1);\R^3)}.$$
Then, there exists a sequence of domains ($\Omega_0^{s,n}$) of class $C^{\infty}$, so that $\os\subset\Omega_0^{s,n}$, and which are defined with a collection of charts $\psi^{m,n}\in H^4(B_-(0,1);\R^3)$ ($m=1,..., N$) so that $\sum_{m=1}^N \|\psi^m-\psi^{m,n}\|_{H^4(B_-(0,1);\R^3)}\rightarrow 0$ as $n\rightarrow\infty$. We then denote the complementary of $\overline\Omega_0^{s,n}$ in $\Omega$ by $\Omega_0^{f,n}$ and $\Gamma_0^n= \partial\Omega_0^{s,n}$. We also assume $n$ large enough so that the different connected components of $\Omega_0^{s,n}$ (if there is more than one solid) do not intersect each other or the boundary of $\Omega$. We denote $\alpha_n=\|\Omega^{s,n}\|_{H^6}$.
\end{lemma}

We now state the optimal regularity assumptions needed in our analysis, and explain the adjustements required to the previous proofs.

\begin{theorem}\label{optimal}
With the same assumptions as in Theorem \ref{main}, except for the following
concerning the regularity of the initial data:
\begin{subequations}
\label{initial}
\begin{align} 
&u_0 \in H^6(\Omega_0^f;{\mathbb R}^3)\cap H^3(\Omega_0^s;{\mathbb R}^3) 
\cap H^1_0(\Omega;{\mathbb R}^3)\cap L^2_{{div},f},\\
& f_s(0)\in H^2(\Omega_0^s;{\mathbb R}^3)\cap H^{3.5}(\Gamma_0;\R^3)\ ,\ (f_s)_t(0)\in H^1(\Omega_0^s;{\mathbb R}^3),\  (f_s)_{tt}(0)\in L^2(\Omega;{\mathbb R}^3), 
\end{align}
\end{subequations}
 the conclusion of Theorem \ref{main} holds.
\end{theorem}

\begin{remark}
We have chosen here to take different forcings for the fluid, that we still note
$f$ with the same assumptions as in Theorem \ref{main}, and the solid, in order to stress out that the higher order regularity required indeed comes from the hyperbolic scaling of the Navier-Stokes equations. The somewhat not so natural condition $f_s(0)\in H^{3.5}(\Gamma_0;\R^3)$ is made in order to get
$w_1\in H^4(\of;\R^3)$ associated to the condition $w_3\in L^2(\of;\R^3)$.
\end{remark} 
\begin{proof} The idea is to first regularize the domains, initial data and modify the forcings in an appropriate way, and then pass to the limit. 

Given $0\le \rho\in \mathcal D(B(0,1))$ with $\int_{B(0,1)} \rho=1$, we define as usual $\rho_n(x)=n^3\rho(nx)$.
 
We first notice that $u_0$, $w_1$, $q_0$ and $q_1$ still have the same regularity in $\of$ as in Theorem \ref{main}. 
We first define 
%the regularized domains $\d \Omega_{0,\delta}^{f,n}=\{x|\ d(x,(\Omega_0^{f,n})^c)>\delta\}$
%and $\d \Omega_{0,\delta}^{s,n}=\{x|\ d(x,(\Omega_0^{s,n})^c)>\delta\}$. Next,
 in $\Omega_0^{f,n}$, $u_0^n=u_0$ and $w_1^n=w_1$, $q_0^n=q_0$, $q_1^n=q_1$ (which is permitted since
$\Omega_0^{f,n}\subset \Omega_0^f$. We next define  $w_2^n$ in $\Omega_0^{f,n}$,
\begin{subequations}
\label{el1}
\begin{align}
-\nu \triangle w_2^n +\nabla q_2^n &=\rho_n\star E(\Omega_0^f) (-\nu\triangle w_2\ +\nabla q_2)\ \text{in}\ \Omega_0^{f,n},\\
\operatorname{div} w_2^n&=-[(a_i^j)_{tt}(0) u_0^i,_j+2 (a_i^j)_{t}(0) w_1^i,_j]\ \text{in}\ \Omega_0^{f,n},\\
w_2^n&=0\ \text{on}\ \partial\Omega,\\
\nu\frac{\partial w_2^n}{\partial N^n} -q_2^n N^n&=\nu\frac{\partial}{\partial N^n} \rho_n\star E(\Omega_{0}^{f}) (w_2) -\rho_n\star E(\Omega_{0}^{f})(q_2) N^n  \ \text{on}\ \Gamma_0^n,
\end{align}
\end{subequations}
where $N^n$ denotes the unit normal exterior to $\Omega_0^{f,n}$, and finally $w_3^n\in L^2(\Omega_0^{f,n};\R^3)$ by
\begin{equation*}
w_3^n=[ \nu (a_l^j a_l^k u,_k),_j - (a_i^j q,_j)_{i=1}^3 +F]_{tt}(0)\ \text{in}\ \Omega_0^{f,n},
\end{equation*}
where the time derivatives on the right-hand side are computed with the usual rules from $ u(0)=u_0^n$, $\partial_t^p u(0)=w_{p}^n$ ($p=1,2)$, $\partial_t^p q(0)=q_{p}^n$
$(p=0,1,2)$.

We next define in the solid $u_0^n$ by
\begin{subequations}
\label{el2}
\begin{align}
L^2 u_0^n  &=L^2 [\rho_n\star E(\os)(u_0)]\ \text{in}\ \Omega_0^{s,n},\\
u_0^n&=(u_0^n)^f\ \text{on}\ \Gamma_0^n,\\ 
L(u_0^n)+ \rho_n\star E(\os) ((f_s)_t(0))&=(w_2^n)^f\ \text{on}\  \Gamma_0^n,
\end{align}
\end{subequations}
where the right-hand sides of the previous boundary conditions come from the fluid regularization previously carried. Note also that
\begin{equation}
\label{el2bis}
Lu_0^n\in H^4(\Omega_0^{s,n};\R^3),
\end{equation}
(with an estimate that may blow up as $n\rightarrow\infty$) since
\begin{align*}
L (Lu_0^n)  &=L^2 [\rho_n\star E(\os)(u_0)]\ \text{in}\ \Omega_0^{s,n},\\
L(u_0^n)&=- \rho_n\star E(\os) ((f_s)_t(0))+(w_2^n)^f\ \text{on}\  \Gamma_0^n.
\end{align*}

%, and where the time derivatives on the right-hand sides at $0$ are evaluated this time with $q(0)=q_0^n$, $q_t (0)=q_1^n$ and
%$v(0)=u_0^n$, $v_t(0)=w_1^n$.  

We can then define $f_0^n$ in $\Omega_0^{s,n}$ by
\begin{align*}
L^2 f_0^n  &=L^2 [\rho_n\star E(\os) (f_s(0))]\ \text{in}\ \Omega_0^{s,n},\\
f_0^n&=(w_1^n)^f\ \text{on}\ \Gamma_0^n,\\ 
c^{mjki}\ (f_0^n,_m\cdot \text{Id},_j+f_0^n,_j\cdot \text{Id},_m) N_k^n
&=  - 2 c^{mjki}\ (u_0^n,_m\cdot\text{Id},_j+ u_0^n,_j\cdot\text{Id},_m) {u_0^n},_k^i\n\\
&\ \ +\nu\ [v^i,_k a^k_l a^j_l]_{tt}(0) N_j^n - [q a^j_i]_{tt}(0) N_j^n\ \text{on}\ \Gamma_0^n \,,
\end{align*}
with the same conventions as for the previous system for the time derivatives evaluated from $\Omega_0^{f,n}$ and $c^{mjkl}\ (u_0^n,_m\cdot\text{Id},_j+ u_0^n,_j\cdot\text{Id},_m) {u_0^n},_k^i\ N_l^n$ evaluated from $\Omega_0^{s,n}$.  
We then define in $\Omega_0^{s,n}$, 
\begin{align*}
w_1^n&=f_0^n,\\
w_2^n&=\bigl[[c^{mjkl} (\eta,_m\cdot\eta,_j-\delta_{mj}) \eta,_k],_l\bigr]_{t}(0)+\rho_n\star E(\os) ((f_s)_t(0)),\\
&=L(u_0^n)+\rho_n\star E(\os) ((f_s)_t(0)),\\
w_3^n&=\bigl[[c^{mjkl} (\eta,_m\cdot\eta,_j-\delta_{mj}) \eta,_k],_l\bigr]_{tt}(0)+ \rho_n\star E(\os) ((f_s)_{tt}(0)), 
\end{align*}
where the time derivatives on the right-hand side are evaluated with $v(0)=u_0^n$, $v_{t}(0)=w_1^n$.
We also define the regularized forcing in the solid
\begin{align*}
f^n(t)&=\rho_n\star E(\os) (f_s(t)-f_s(0))\ +f_0^n \text{in}\ \Omega_0^{s,n}.
\end{align*}
We then have $u_0^n$, $w_1^n$, $w_2^n$ in $H^1_0(\Omega;\R^3)\cap H^4(\Omega_0^{f,n};\R^3)\cap H^4(\Omega_0^{s,n};\R^3)$ and $\operatorname{div} u_0^n=0$ in $ \Omega_0^{f,n}$, $w_3^n\in L^2(\Omega;\R^3)$, with 
\begin{subequations}
\label{14.2}
\begin{align}
&\|E(\Omega_0^{f,n})(u_0^n)-u_0\|_{H^4(\of;\R^3)}+\|E(\Omega_0^{s,n})(u_0^n)-u_0\|_{H^3(\os;\R^3)}\rightarrow 0\ \text{as}\ n\rightarrow\infty, \label{14.2a}\\
&\|E(\Omega_0^{f,n})(w_1^n)-w_1\|_{H^4(\of;\R^3)}+\|E(\Omega_0^{s,n})(w_1^n)-w_1\|_{H^2(\os;\R^3)}\n\\
&+\|w_2^n-w_2\|_{H^1(\Omega;\R^3)}+\|E(\Omega_0^{f,n})(w_2^n)-w_2\|_{H^2(\Omega_0^f;\R^3)}\rightarrow 0\ \text{as}\ n\rightarrow\infty, \label{14.2c}\\
&\|u_0^n\|_{H^6(\Omega_0^{s,n};\R^3)}\le  \beta_n,\ \ \|w_1^n\|_{H^4(\Omega_0^{s,n};\R^3)}\le \beta_n,\ \ \|w_2^n\|_{H^4(\Omega_0^{s,n};\R^3)}\le \beta_n ,\label{14.2d}\\
&\|w_3^n-w_3\|_{L^2(\Omega;\R^3)}\rightarrow\infty\ \text{as}\ n\rightarrow\infty,\label{14.2e}
\end{align}
\end{subequations}
where $\beta_n$ is a given polynomial expression of $\alpha_n$ and $n$. We briefly explain how those constants appear. For instance, for the first estimate of (\ref{14.2d}), we have by elliptic regularity on (\ref{el2}) that $\|u_0^n\|_{H^6(\Omega_0^{s,n};\R^3)}$ is bounded by a sum of terms, one of which being $P(\|\Omega_0^{s,n}\|_{H^6})\ \|(w_2^n)^f\|_{H^4(\Omega_0^{f,n};\R^3)}$, $P$ being a polynomial which does not depend on $n$. Next, still by elliptic regularity on (\ref{el1}), we have that $\|(w_2^n)^f\|_{H^4(\Omega_0^{f,n};\R^3)}$ is bounded by  a sum of terms such as $\|\rho_n\star E(\of)(\triangle w_2)\|_{H^2(\Omega_0^{f,n};\R^3)}$. This particular term, by the properties of the convolution, is in turn bounded by
${n^3} \|E(\Omega_0^{f,n})(w_2)\|_{H^1(\R^3;\R^3)}$. This shows that a term of the type $P(\alpha_n){n^3} \|w_2\|_{H^1(\of;\R^3)}$ appears in the
sum of all terms bounding  $\|u_0^n\|_{H^6(\Omega_0^{s,n};\R^3)}$. Since the other terms in the sum can be dealt with similarly, this explains our estimate (\ref{14.2d}).

For the pressures, we have
\begin{align}
\label{14.3}
\|E(\Omega_0^{f,n})(q_0^n)-q_0\|_{H^3(\of;\R)}+&\|E(\Omega_0^{f,n})(q_1^n)-q_1\|_{H^3(\of;\R)}\n\\
& +\|E(\Omega_0^{f,n})(q_2^n)-q_2\|_{H^1(\of;\R)}\rightarrow 0\ \text{as}\ n\rightarrow\infty.
\end{align}
Since the initial data $u_0^n$ and forcings $f^n (0)$, $f^n_t (0)$, $f^n_{tt} (0)$ are smooth enough to ensure the regularity properties (\ref{14.2}), we then deduce that we have similarly as for theorem (\ref{main}) the existence of a solution $w_n$ of a system similar to (\ref{weakW}) with $f$, $u_0$, $\of$, $\os$ being replaced by their counterparts with an exponent $n$, and $b_{\kappa}$, $c_{\kappa}$, $d_{\kappa}$ being replaced by $b_n$, $c_n$, $d_n$ (with the choice $\d\kappa=\frac{1}{n ({\beta_n}+1)}$) given by
\begin{align*}
b_{n}(\phi)=& \frac{1}{n ({\beta_n}+1)}(c^{ijkl} {{w_{ 2}^n}},_l^k, \phi^i,_j)_{L^2(\Omega_0^{s,n};{\mathbb R})},\\
c_{n}(\phi)=&\ \frac{1}{n (\beta_n+1)} (c^{ijkl}   {w^n_1},_l^k, \phi^i,_j)_{L^2(\Omega_0^{s,n};{\mathbb R})}+ (  {w^n_2}-w_2, \phi)_{L^2(\Omega_0^{f,n};{\mathbb R})}\\
&\  + (-[(a_i^j q)_t (0)N^n_j]_{i=1}^3+ [(a_l^j a_l^k) u,_k]_t(0) N^n_j ,\ \phi)_{L^2(\Gamma_0^n;{\mathbb R}^3)}\\
&\  -(  c^{ijkl} [(\eta,_i\cdot\eta,_j-\delta_{ij}) \eta_k]_t(0) N^n_l,\ \phi)_{L^2(\Gamma_0^n;{\mathbb R}^3)},
\end{align*}
\begin{align*}
d_{n}(\phi)=&\ \frac{1}{n (\beta_n+1)} (c^{ijkl}   {u^n_0},_l^k, \phi^i,_j)_{L^2(\Omega_0^{s,n};{\mathbb R})}\\
&\  + (-[(a_i^j q) (0)N^n_j]_{i=1}^3+ [(a_l^j a_l^k) u,_k](0) N^n_j ,\ \phi)_{L^2(\Gamma_0^n;{\mathbb R}^3)}, 
\end{align*}
where the time derivatives are computed with
a velocity satisfying $u(0)=u^n_0$, $u_t(0)=w^n_1$ and a pressure such that $q(0)=q_0^n$, $q_t(0)=q_1^n$. Note that by construction, the solutions $w_n$ to these problems in $\Omega$ satisfy $w_n(0)=u^n_0$, ${w_n}_t(0)=w^n_1$, ${w_n}_{tt}(0)=w^n_2$, ${w_n}_{ttt}(0)=w^n_3$.
Next, we proceed to energy estimates similar to Section \ref{8}. The bounds obtained
are similar, except that this time the terms associated to $b_n$, $c_n$ and $d_n$ tend to zero as $n\rightarrow\infty$. This is clear from the convergence results (\ref{14.2}), (\ref{14.3}) for the integral terms associated to the fluid. The terms associated to the solid asymptotically tend to zero by properties of the convolution. For instance, with the notations of Section \ref{8}, for
$\phi_p=-[\rho_p\star(\z^2\ w_n\circ\P)],_{\ao\ao\at\at\att\att}\circ\P^{-1}$, we get after change of variables, an integration by parts in time, three integrations by parts in space:
\begin{align*}
\bigl|\int_0^t \frac{t'^2}{2 n (\beta_n+1)}&(c^{ijkl} {{w_{ 2}^n}},_l^k, \phi_p^i,_j)_{L^2(\Omega_0^{s,n};{\mathbb R})}\ dt'\bigr|\\
&\le 
\frac{C}{n (\beta_n+1)} \|w_2^n\|_{H^4(\Omega_0^{s,n};{\mathbb R}^3)}\|\eta_n\|_{L^{\infty}(0,t;H^4(\Omega_0^{s,n};{\mathbb R}^3))},
\end{align*}
and thus with our estimate (\ref{14.2d}), we have
\begin{align*}
\bigl|\int_0^t \frac{t'^2}{2 n (\beta_n+1)}(c^{ijkl} {{w_{ 2}^n}},_l^k, \phi_p^i,_j)_{L^2(\Omega_0^{s,n};{\mathbb R})}\ dt'\bigr|\le 
\frac{C}{n} \|(w_n,q_n)\|_{Z^n_t},
\end{align*}
where $Z^n_t$ denotes the same type of space as $Z_t$ with $\Omega_0^s$ and $\Omega_0^f$ being replaced by their counterparts with an exponent $n$. 
This type of estimates thus shows that this term does not change the energy inequalities
in Section \ref{8}. We can thus reproduce the arguments of Section \ref{9}, establishing that $(w_n,q_n)$ can be defined over a time $T$ independently of $n$, and that its norm in $Z^n_T$ depends solely on $N(u_0^n,(w_i^n)_{i=1}^3)+N((q_i)_{i=0}^2)+M(f^n,0,0)$ and thus, thanks to the estimates (\ref{14.2a}), (\ref{14.2c}), (\ref{14.2e}), solely
on $N(u_0,(w_i)_{i=1}^3)+N((q_i)_{i=0}^2)+M(f,0,0)$. We can then consider the sequence $(E(\Omega)(w_n),E(\Omega_0^{f,n})(q_n))$ which is bounded in a space similar as $Z_T$, but defined on $\R^3$ and extract (with respect to $n$) a weakly
convergent sequence in a space modified from $Y_T$ by replacing the condition $u\in H^1_0(\Omega;\R^3)$ by $u\in H^1(\Omega;\R^3)$ . By the classical compactness results, we next see that the weak limit $(v,q)\in Z_T$ and is a solution of (\ref{nsl}) with $f$ as forcing and $v(0)=u_0$. This solution is also unique in $Z_T$.
\end{proof}

\section{The case of incompressible elasticity}
\label{13}

In this Section, we explain how to treat the supplementary difficulties appearing when the incompressibility constraint is added in the solid. This leads to the
same system as (\ref{nsl}), with the addition of the condition $\operatorname{det}\nabla\eta=1$ {\it a.e.} in $\os$ and the addition of $[(a_i^k q),_k]_{i=1}^3$ on the left-hand side of (\ref{nsl.d}) and the addition of $-q a_i^j N_j$ (the trace of $q$ being from the solid phase in this new term) on the left-hand side of
(\ref{nsl.e}). We now state our result and explain how to overcome the additional difficulties related to this constraint.

We first update our functional frameworks. While $X_T$ and $W_T$ do not change,
$Y_T$ and $Z_T$ become respectively
\begin{align*}
Y_T= \{ (v,q)\in X_T\times &L^2(0,T;L^2(\Omega;{\mathbb R}))|\  \partial_t^n q\in L^2(0,T;H^{3-n}(\Omega_0^f;{\mathbb R})),\\
& \ \partial_t^n q\in L^2(0,T;H^{3-n}(\Omega_0^s;{\mathbb R})) (n=0,1,2)\} , 
\end{align*}
\begin{align*}
Z_T= \{ &(v,q)\in W_T\times  L^2(0,T;L^2(\Omega;{\mathbb R}))|\ \partial_t^n q\in L^2(0,T;H^{3-n} (\Omega_0^f;{\mathbb R})), \\
&\ \partial_t^n q\in L^2(0,T;H^{3-n} (\Omega_0^s;{\mathbb R}))
(n=0,1,2)|\ q_{tt}\in L^{\infty}(0,T;L^2(\Omega;\R))\}\ . 
\end{align*}
\begin{remark}
Whereas the pressure in the solid satisfies $\partial_t^n q\in L^\infty(0,T;H^{3-n} (\Omega_0^s;{\mathbb R}))$ $(n=0,1,2)$, it appears that the limit pressures $q_\kappa$ are controlled uniformly in the norm of $Z_T$ and seemingly not in these norms. Note also that whereas the velocity field is smoother in the fluid phase for the solution of our next theorem, the pressure field is actually smoother in the solid phase. Whereas our artificial viscosity smoothes the velocity field in the solid, it also interestingly makes the pressure in the solid for the regularized system less
smooth than the one associated to the solution of the constrained problem, which is source of difficulties that we shall describe later.
\end{remark} 

We now state our result:

\begin{theorem}
With the same regularity assumptions as in Theorem \ref{optimal} and assuming
that the compatibility conditions associated to our new system at $t=0$ hold
(for the sake of conciseness we do not state them here), the conclusion of
Theorem \ref{main} holds for the case where the incompressibility constraint
is added to the solid part. Furthermore, $\partial_t^n q\in L^\infty(0,T;H^{3-n} (\Omega_0^s;{\mathbb R}))$ $(n=0,1,2)$.
\end{theorem}
\begin{proof}
The extra regularity (with respect to the norm of $Z_T$) on the pressure in the solid simply comes from the equation
$$v_t - c^{mjkl} [(\eta,_m\cdot\eta,_j-\delta_{ij}) \eta,_k],_l +a_i^j q,_j= f\ \ \text{in} \ \ (0,T)\times \Omega_0^s, $$
which once the regularity for the solution $w\in W_T$ is known provides immediately the result. We now explain how to obtain a solution in $Z_T$.

The beginning of the proof follows the same lines as for the compressible elasticity case.
We first assume that the initial data satisfies the regularity assumptions of Theorem \ref{main}, and define the same smoothed problem as (\ref{smoothie}) with the corresponding updates for the incompressibility constraint. We then define the same fixed point linear problem as (\ref{linear}) where the condition $a_i^k w^i,_k=0$ in $\os$ is added (the $a_i^k$ being computed from the given $v$) and add $a_i^k q,_k$ on the left-hand side of (\ref{linear.d}) and $-q a_i^j N_j$ (the traces being taken from $\os$) on the left-hand side of (\ref{linear.e}).
 
We then proceed as in \cite{CoSh2004} to construct a solution
to this system by a penalty method (the penalty term being this time defined over $\Omega$) and get the same type of regularity result. This provides us with a solution $(w_{\kappa},q_{\kappa})$, that we also denote by $(\w,\q)$, for
the incompressible version of (\ref{smoothie}) on a time $T_{\kappa}$ shrinking to zero. As for the compressible case, $(w_\kappa,q_\kappa)$ is in $Z_{T_\kappa}$, and since our smoothed problem has a parabolic artificial viscosity, we also have for the velocity in the solid the regularity $\partial_t^n w\in L^2(0,T_{\kappa};H^{4-n}(\os;\R^3))$ ($n=0,1,2,3$) (with estimates that blow up as $\kappa\rightarrow 0$). Thus, $(w_\kappa,q_\kappa)\in\tilde Z_{T_\kappa}$ with
$$\tilde Z_t=\{(w,q)\in Z_t|\ \partial_t^n w\in L^2(0,t;H^{4-n}(\os;\R^3)) (n=0,1,2)\},$$ endowed with the norm
$$\|(w,q)\|^2_{\tilde Z_t}=\|(w,q)\|^2_{Z_t}+\kappa^2 \sum_{n=0}^2
\|\partial_t^n w\|^2_{L^2(0,t;H^{4-n}(\os;\R^3))}.$$
We next proceed as in Section \ref{8} to get energy estimates, that will be carried this time for the $\kappa$ dependent norm of $\tilde Z_t$, independently of $\kappa$ on $[0,T_{\kappa}]$, and for such a purpose it is important to keep the $\kappa^2$ factor in the definition of the norm. We could extend the sum to $n=3$, whereas it is not necessary. 

As before, the first set of estimates has to be carried on the highest order time derivative. Our energy inequality (\ref{energywttt}) has the same form, except that the integrals over $\Omega_0^f$ where $\q$ appears has to be taken this time on $\Omega$. The part over $\Omega_0^f$ is estimated as before. We now explain how to deal with the integrals set on $\os$ for the pressure, which indeed needs some justifications given that the velocity in the solid is not controlled
uniformly in $\kappa$ in a space as smooth as the velocity in the fluid, while the pressure is controlled in the same type of spaces in both phases.

\subsection{Estimates on $\w_{ttt}$}
Here $t$ denotes any time in $(0,T_\kappa)$. 
The most difficult integrals  set in $[0,t]\times\os$ and associated to the incompressibility constraint in
the solid are $\d K_1=\int_0^t\int_{\os} \q_{tt} (\a_i^j)_t \w_{ttt}^i,_j$ and $\d K_2=\int_0^t\int_{\os} \q (\a_i^j)_{ttt} \w_{ttt}^i,_j$, the others being either less difficult or similar to estimate. 
  
\noindent{Step 1.} For $K_1$, if we denote $N^s=-N$, we have
\begin{align}
|K_1|&=\bigl|-\int_0^t\int_{\os} (\q_{tt}),_j (\a_i^j)_t \w_{ttt}^i +\int_0^t\int_{\Gamma_0} \q_{tt} (\a_i^j)_t \w_{ttt}^i N_j^s\bigr|\n\\
&\le C [\int_0^t \|\q_{tt}\|_{H^1(\os;\R)}\|\e\|_{H^3(\os;\R^3)}\|\w\|_{H^3(\os;\R^3)} \|\w_{ttt}\|_{L^2(\os;\R^3)}\n\\
&\ \ \ + \int_0^t \|\q_{tt}\|_{H^{\frac{1}{2}}(\os;\R)}\|\e\|_{H^3(\os;\R^3)}\|\w\|_{H^3(\os;\R^3)} \|\w_{ttt}\|_{H^{\frac{1}{2}}(\of;\R^3)}]\n\\
&\le C \sqrt{t} \|(\w,\q)\|^4_{Z_t} + C t^{\frac{1}{2}} \sup_{[0,t]}[\ \|\q_{tt}\|^{\frac{1}{2}}_{L^2(\os;\R)}\|\w_{ttt}\|^{\frac{1}{2}}_{L^2(\of;\R^3)}] \|(\w,\q)\|^{{3}}_{Z_t}\\
&\le C t^{\frac{1}{4}} \|(\w,\q)\|^4_{Z_t} ,
\label{in1} 
\end{align}
where we have used the continuity of $\w_{ttt}$ in the sense of traces along
$\Gamma_0$ to bound the $L^2(\Gamma_0;\R^3)$ norm of $\w_{ttt}$ by means of the $H^{\frac{1}{2}}(\of;\R^3)$ norm. Note that we have also used the fact that the $L^\infty(L^2)$ norm of $\q_{tt}$ is in the definition of the norm of
$Z_t$. In order to get an estimate on this norm, we would proceed in a way similar as
to get (\ref{I8}) in Section \ref{9}.

\noindent{Step 2}. Concerning $K_2$, we have by integrating by parts in space 
\begin{align*}
K_2&=-\int_0^t\int_{\os} \q,_j (\a_i^j)_{ttt} \w_{ttt}^i + \int_0^t\int_{\Gamma_0} \q (\a_i^j)_{ttt} \w_{ttt}^i N_j^s,
\end{align*}
since our artificial viscosity provides the regularity $\w_{tt}\in L^2(0,T_{\kappa};H^2(\os;\R^3))$ and $\w_{ttt}\in L^2(0,T_{\kappa};H^1(\os;\R^3))$
(with estimates that may blow up as $\kappa\rightarrow 0$). The difficulty here comes from the second integral. Whereas as for $K_1$ we can estimate the trace of $\w_{ttt}$ on $\Gamma_0$ from the fluid, we have to take the norm of $\nabla \w_{tt}$ in $H^{-0.5}(\Gamma_0;\R^9)$, which is problematic given that the norm $Z_t$ contains
only its $L^2(\os;\R^9)$ norm. In order to circumvent this difficulty, we notice that the same formula holds if we replace $\w_{ttt}$ by $E(\of)(\w_{ttt}^f)$ (the extension to $\R^3$ of the velocity
in the fluid). Since $\w_{ttt}=\w_{ttt}^f$ on $\Gamma_0$, we have $\w_{ttt}=E(\of)(\w_{ttt}^f)$ on $\Gamma_0$, which implies:
\begin{align*}
K_2=-\int_0^t\int_{\os} \q,_j (\a_i^j)_{ttt} \w_{ttt}^i &+ \int_0^t\int_{\os} \q (\a_i^j)_{ttt} E(\of)(\w_{ttt}^f),_j^i\\
& + \int_0^t\int_{\os} \q,_j (\a_i^j)_{ttt} E(\of)(\w_{ttt}^f)^i,
\end{align*}
and thus,
\begin{align*}
K_2\le &\ C \int_0^t \|\q\|_{H^3(\os;\R)} \|\w_{tt}\|_{H^1(\os;\R^3)}
\|\e\|_{H^3(\os;\R^3)} \|\w_{ttt}\|_{L^2(\Omega;\R^3)}\\
&+ C \int_0^t \|\q\|_{H^3(\os;\R)} \|\w_{t}\|_{H^1(\os;\R^3)}
\|\w\|_{H^3(\os;\R^3)} \|\w_{ttt}\|_{L^2(\Omega;\R^3)}\\
&+ C\int_0^t \|q_0+\int_0^\cdot \q_t\|_{H^2(\os;\R)} \|\w_{tt}\|_{H^1(\os;\R^3)}
\|\e\|_{H^3(\os;\R^3)} \|\w_{ttt}\|_{H^1(\of;\R^3)}\\
&+ C\int_0^t \|q_0+\int_0^\cdot \q_t\|_{H^2(\os;\R)} \|\w_{t}\|_{H^1(\os;\R^3)}
\|\w\|_{H^3(\os;\R^3)} \|\w_{ttt}\|_{H^1(\of;\R^3)}\\
\le &\ C\ \|(\w,\q)\|^3_{Z_t} \int_0^t \|\q\|_{H^3(\os;\R)}
+ C\ \|(\w,\q)\|^2_{Z_t}\|q_0\|_{H^2(\os;\R)} \int_0^t \|\w_{ttt}\|_{H^1(\of;\R^3)}\\
& + C\sqrt{t}\ \|(\w,\q)\|^2_{Z_t}\|\q_t\|_{L^2(0,t;H^2(\os;\R))} \int_0^t \|\w_{ttt}\|_{H^1(\of;\R^3)}\\
\le &\ C\sqrt{t}\ [\ \|(\w,\q)\|^4_{Z_t}+N((q_i)_{i=0}^2)\ ] .
\end{align*}

The most difficult integral set at time $t$ on $\os$ and containing $\q$ is
$$\d K_3= \int_{\os} \q_{tt} (\a_i^j)_t \w_{tt}^i,_j,$$ for which we apparently just have an estimate of the type $|I_3|\le C \|(\w,\q)\|^2_{Z_t}$ (without any small parameter in front). We now explain how to treat this difficulty.

\noindent{Step 3.} We first notice that 
\begin {align*}
K_3=-\int_{\os} \q_{tt},_j (\a_i^j)_t \w_{tt}^i+\int_{\Gamma_0} \q_{tt} (\a_i^j)_t \w_{tt}^i N_j^s.
\end{align*}
If we could say that $\q_{tt}$ is $L^\infty(H^1)$ controlled, the $L^\infty(L^2)$ control of $\w_{tt}$ would give us a suitable bound for $K_3$. Whereas we have seen in the statement of our theorem that $q_{tt}$ for the limit solution is indeed in $L^\infty(H^1)$, we cannot seemingly get such a bound on the approximate pressures $\q_{tt}$. In order to get around this, we introduce similarly as in the previous step the extension to the solid domain of the velocity in the fluid. Since a similar
integration by parts formula holds when we replace $\w_{tt}$ by $E(\of)(\w^f_{tt})$, we deduce
\begin {align}
K_3=-\int_{\os} \q_{tt},_j (\a_i^j)_t \w_{tt}^i+\int_{\os} \q_{tt} (\a_i^j)_t E(\of)(\w^f_{tt})^i,_j + \int_{\os} \q_{tt},_j (\a_i^j)_t E(\of)(\w^f_{tt})^i.
\label{inc4}
\end{align}
The easier term to estimate is $\d K_3^2= \int_{\os} \q_{tt} (\a_i^j)_t E(\of)(\w^f_{tt})^i,_j$, for which we have for an arbitrary $\delta>0$: 
\begin{align}
|K_3^2| &\le 
C \|\q_{tt}\|_{L^2(\os;\R)} \|\text{Id}+\int_0^t \w\|_{H^3(\os;\R^3)} \|u_0+\int_0^t \w_t\|^{\frac{1}{4}}_{H^2(\os;\R^3)} \| \w\|^{\frac{3}{4}}_{H^3(\os;\R^3)} \n\\ 
&\qquad\times [ \|w_2\|_{H^1(\of;\R^3)}+\sqrt{t} \|\w_{ttt}\|_{L^2(0,t;H^1(\of;\R^3))}]\n\\
&\le C \|\q_{tt}\|_{L^2(\os;\R)}  [1+{t} \|(\w,\q)\|_{Z_t}][N(u_0,(w_i)_{i=1}^3)+t^{\frac{1}{4}}\|(\w,\q)\|^{\frac{1}{4}}_{Z_t}]\ \|(\w,\q)\|^{\frac{3}{4}}_{Z_t}\n\\
&\qquad\qquad \times [N(u_0,(w_i)_{i=1}^3)+t^{\frac{1}{2}}\|(\w,\q)\|_{Z_t}]\n\\
&\le C \|(\w,\q)\|_{Z_t} [N(u_0,(w_i)_{i=1}^3)+t^{\frac{1}{4}}\|(\w,\q)\|^{\frac{1}{4}}_{Z_t}] \|(\w,\q)\|^{\frac{3}{4}}_{Z_t}\n\\
&\qquad\times [N(u_0,(w_i)_{i=1}^3)+t^{\frac{1}{2}}\|(\w,\q)\|_{Z_t}]^2\n\\
&\le \delta \|(\w,\q)\|^2_{Z_t}+C_\delta N(u_0,(w_i)_{i=1}^3) + C_{\delta} t^{\frac{1}{4}} \|(\w,\q)\|^4_{Z_t}.
\label{inc4bis} 
\end{align}
For the first integral, the nonlinear elastodynamics equation in
$\os$ provides 
\begin{align*}
\nabla \q_{tt}= \a^{-1} [-\w_{ttt}+\kappa L\w_{tt} +c^{ijkl}[(\e,_i\cdot\e,_j-\delta_{ij})\e,_k]_{tt},_l -2\a_t\nabla \q_t-\a_{tt}\nabla \q+f_{tt}+\kappa h],
\end{align*}
leading us for $\d K_3^1=\int_{\os} \q_{tt},_j (\a_i^j)_t \w_{tt}^i$ to (since $\a^{-1}=\nabla\e$ in virtue of $\operatorname{det}\nabla\e=1$),
\begin{align}
K_3^1&=\int_{\os} [\nabla\e\ [-\w_{ttt} +c^{ijkl}[(\e,_i\cdot\e,_j-\delta_{ij})\e,_k]_{tt},_l -2\a_t\nabla \q_t-\a_{tt}\nabla \q+f_{tt}+\kappa h]]^j (\a_i^j)_t \w_{tt}^i\n\\
&\ \ + \kappa\int_{\os} [\nabla\e\ [ L\w_{tt}]]^j (\a_i^j)_t \w_{tt}^i.
\label{inc5} 
\end{align}
The integrals on the first line of this equality do not give any trouble and can be estimated in the same fashion. For instance, we have for
$$\d K_3^3=\int_{\os} [\nabla\e [c^{ijkl}[(\e,_i\cdot\e,_j-\delta_{ij})\w_t,_{lk}]]^j (\a_i^j)_t \w_{tt}^i,$$
\begin{align}
|K_3^3|&\le C \|\text{Id}+\int_0^t \w\|^4_{H^3(\os;\R^3)}\|\w_t\|_{H^2(\os;\R^3)} \|w_2+\int_0^t \w_{ttt}\|_{L^2(\os;\R^3)}\n\\
&\qquad\qquad\qquad\qquad \times \|u_0+\int_0^t \w_{t}\|^{\frac{1}{4}}_{H^2(\os;\R^3)}\|\w\|^{\frac{3}{4}}_{H^3(\os;\R^3)}\n\\
&\le \delta \|(\w,\q)\|^2_{Z_t} + C_{\delta} t^{\frac{1}{4}} \|(\w,\q)\|^{7}_{Z_t} +C_\delta  N(u_0,(w_i)_{i=1}^3)).
\label{inc6}
\end{align} 
Now, the difficult term to handle is $\d K_3^4=\kappa\int_{\os} [\nabla\e[ L\w_{tt}]]^j (\a_i^j)_t \w_{tt}^i$. We first write the divergence form $L\w_{tt}^p=\sigma,_m^{mp}(\w_{tt})$, and integrate by parts:
\begin{align*}
K_3^4&=-\kappa\int_{\os} [\nabla \e,_m[\sigma^{mp}(\w_{tt})]_{p=1}^3]^j (\a_i^j)_t \w_{tt}^i\n\\
&\ \ \ -\kappa \int_{\os} [\nabla \e[\sigma^{mp}(\w_{tt})]_{p=1}^3]^j [(\a_i^j)_t \w_{tt}^i],_m\n\\
&\ \ \ +\kappa \int_{\Gamma_0} [\nabla \e[\sigma^{mp}(\w_{tt})]_{p=1}^3]^j (\a_i^j)_t \w_{tt}^i N_m^s,
\end{align*}
leading us to
\begin{equation}
\label{inc7}
|K_3^4-\kappa \int_{\Gamma_0} [\nabla \e[\sigma^{mp}(\w_{tt})]_{p=1}^3]^j (\a_i^j)_t \w_{tt}^i N_m^s|\le C \kappa \|(\w,\q)\|^5_{Z_t},
\end{equation}
and thus by putting together (\ref{inc5}), (\ref{inc6}) and (\ref{inc7}),
\begin{align}
&\bigl|\int_{\os} \q_{tt},_j (\a_i^j)_t \w_{tt}^i-\kappa \int_{\Gamma_0} [\nabla \e[\sigma^{mp}(\w_{tt})]_{p=1}^3]^j (\a_i^j)_t \w_{tt}^i N_m^s\bigr|\n\\
&\le C \kappa \|(\w,\q)\|^5_{Z_t} +C_\delta t^{\frac{1}{4}} \|(\w,\q)\|^7_{Z_t}+\delta \|(\w,\q)\|^2_{Z_t}\\
&\qquad\quad  +C_\delta  [N(u_0,(w_i)_{i=1}^3)+N((q_i)_{i=0}^2)+M(f,\kappa g,\kappa h)].
\label{inc8}
\end{align}

Now, the apparent problem comes from the term $\sigma^{mp}( \w_{tt})$ on $\Gamma_0$ that should be taken in $H^{-0.5}(\Gamma_0;\R)$, which is troublesome since
the norm in $Z_t$ appropriate for our limit process only contains
its $L^\infty(0,t;L^2(\os;\R))$ norm. In order to circumvent this, we notice that we also have, since $E(\of)(\w_{tt})$ is at least as smooth as $\w_{tt}$ in $\os$,
\begin{align}
&\bigl|\int_{\os} \q_{tt},_j (\a_i^j)_t E(\of)(\w_{tt})^i-\kappa \int_{\Gamma_0} [\nabla \e[\sigma^{mp}(\w_{tt})]_{p=1}^3]^j (\a_i^j)_t E(\of)(\w_{tt}^f)^i N_m\bigr|\n\\
&\le C \kappa \|(\w,\q)\|^5_{Z_t} + C_{\delta} t^{\frac{1}{4}} \|(\w,\q)\|^7_{Z_t}+\delta \|(\w,\q)\|^2_{Z_t}\n\\
&\qquad\qquad  +C_\delta  [N(u_0,(w_i)_{i=1}^3)+N((q_i)_{i=0}^2)+M(f,\kappa g,\kappa h)],
\label{inc9}
\end{align}
 leading us, since $\w=E(\of)(\w^f)$ on $\Gamma_0$, to
\begin{align}
&\bigl|\int_{\os} \q_{tt},_j (\a_i^j)_t \w_{tt}^i-\int_{\os} \q_{tt},_j (\a_i^j)_t E(\of)(\w_{tt})^i\bigr|\n\\
&\le C \kappa \|(\w,\q)\|^5_{Z_t}+ C_{\delta} t^{\frac{1}{4}} \|(\w,\q)\|^7_{Z_t}+\delta \|(\w,\q)\|^2_{Z_t}\n\\
&\qquad\qquad  +C_\delta  [N(u_0,(w_i)_{i=1}^3)+N((q_i)_{i=0}^2)+M(f,\kappa g,\kappa h)].
\label{inc10}
\end{align}
Thus, by using (\ref{inc4}), (\ref{inc4bis}) and (\ref{inc10}), we have
\begin{align}
|I_3|
\le & (C \kappa + C_{\delta} {t}^{\frac{1}{4}}) \|(\w,\q)\|^7_{Z_t}+\delta \|(\w,\q)\|^2_{Z_t} \n\\
& +C_\delta  [N(u_0,(w_i)_{i=1}^3)+N((q_i)_{i=0}^2)+M(f,\kappa g,\kappa h)].
\label{inc11}
\end{align}
Thus, we finally arrive to estimates analogous to (\ref{I7}) and (\ref{I8}), with the right-hand side being of the same type as in (\ref{inc11}). 
%\begin{align}
%  \sup_{[0,t]} \|\tilde w_{ttt} \|^2_{L^2(\Omega;{\mathbb R}^3)}& 
%+\int_0^{t} \|{\tilde w}_{ttt}\|^2_{H^1(\Omega_0^f;{\mathbb R}^3)} 
%+ \sup_{[0,t]} \|{{\tilde w}_{tt}}\|^2_{H^1(\Omega_0^s;{\mathbb R}^3)}\nonumber\\
%& \le  C_{\delta} [N(u_0,(w_i)_{i=1}^3))+M(f,g,h)+N((q_i)_{i=0}^2))]\n\\
%&\ \ \ +C_{\delta} t [1+\|(\tilde w,\tilde q)\|^6_{Z_T}] + C {\delta} \|(\tilde w,\tilde q)\|^2_{Z_T} + C\kappa \|(\w,\q)\|^5_{Z_t}.
%\label{energywtttlimitinc}
%\end{align}
\subsection{Estimate on $\w_{tt}$ and $\w_t$.}

With the same arguments as in the next subsection, we have for $n=2,1$:
\begin{align}
&\|\partial_t^n \w\|^2_{L^2(0,t;H^{4-n}(\of;\R^3))}+ \|\partial_t^n \q\|^2_{L^2(0,t;H^{3-n}(\of;\R^3))}+\|\partial_t^n \e\|^2_{L^{\infty}(0,T;H^{4-n}(\os;\R^3))}\n\\
& +\kappa^2 \|\partial_t^n \w\|^2_{L^2(0,t;H^{4-n}(\os;\R^3))}+ \|\partial_t^n \q\|^2_{L^2(0,t;H^{3-n}(\os;\R^3))}\n\\
&\qquad\qquad\qquad\qquad\le C_{\delta} [N(u_0,(w_i)_{i=1}^3)+M(f,\kappa g,\kappa h)+N((q_i)_{i=0}^2)]\n\\
&\qquad\qquad\qquad\qquad\qquad\qquad+(C\kappa+C_{\delta} t^{\frac{1}{4}}) \|(\tilde w,\tilde q)\|^7_{Z_t} + C {\delta} \|(\tilde w,\tilde q)\|^2_{Z_t}. 
\label{wh3fbis}
\end{align}

 We now explain on the case of the highest space derivative how to obtain elliptic estimates independent of $\kappa$, since the addition of the pressure term does not allow to use Lemma \ref{key} directly in the present case.

\subsection{Estimate on $\e$ in $\os$.}

\subsubsection{\bf Regularity of the trace of $\e$}.
First, by proceeding as in Section \ref{8}, and as for the case of the highest order time derivative, we get an estimate for the trace similar to (\ref{tracew}), with a majorant of the same type as in (\ref{inc11}). We explain hereafter how to handle the
estimates related to the pressure in the solid in order to get this trace
estimate since difficulties different than in the higher order time derivative
problem appear in the higher order space derivative problem.

\noindent Step 1. Let $\d Q_1=\int_0^{t}\int_{\R^3_-} [{Q}\ \tilde b_i^j],_{\ao\at\att}  {[\z^2 W^i],_{\ao\at\att j}}$. \hfill\break
Then, \begin{align*}
Q_1=Q_{2}+Q_{3}+Q_4,
\end{align*}
with
\begin{align*}
Q_{2}&=\int_0^{t}\int_{\R^3_-} {Q}\ \tilde b_i^j,_{\ao\at\att}  [\z^2 W^i],_{\ao\at\att j},\\
Q_{3}&=\int_0^{t}\int_{\R^3_-} [[{Q}\ \tilde b_i^j],_{\ao\at\att}-{Q},_{\ao\at\att}\ \tilde b_i^j-{Q}\ \tilde b_i^j,_{\ao\at\att}] [\z^2 W^i],_{\ao\at\att j},\\
Q_{4}&=\int_0^{t}\int_{\R^3_-} {Q},_{\ao\at\att}\ \tilde b_i^j  [\z^2 W^i],_{\ao\at\att j}.
\end{align*}
For $Q_2$, we first notice that for $\theta=\e\circ\P$, if $\epsilon^{ijk}$ is the sign of the permutation between $\{i,j,k\}$ and $\{1,2,3\}$ if $i,j,k$ are distinct, and is set to zero otherwise,
\begin{align*}
\tilde b_i^j,_{\ao\at\att}   W^i,_{\ao\at\att j}&=\frac{1}{2}\epsilon^{mni}\epsilon^{pqj} [\theta,_p^m\theta,_q^n],_{\ao\at\att} W,_{j\ao\at\att}^i\\
&= \epsilon^{mni}\epsilon^{pqj} \theta,_{p\ao\at\att}^m\theta,_q^n  W,_{j\ao\at\att}^i\\
&\ \ +\frac{1}{2}\sum_{\sigma\in\Sigma_3} \epsilon^{mni}\epsilon^{pqj} \theta,_{p\alpha_{\sigma(1)}}^m\theta,_{q\alpha_{\sigma(2)}\alpha_{\sigma(3)}}^n  W,_{j\ao\at\att}^i\\
&\ \ +\frac{1}{2}\sum_{\sigma\in\Sigma_3} \epsilon^{mni}\epsilon^{pqj} \theta,_{p\alpha_{\sigma(1)}\alpha_{\sigma(2)}}^m\theta,_{q\alpha_{\sigma(3)}}^n  W,_{j\ao\at\att}^i\\
&= \frac{1}{2}\epsilon^{mni}\epsilon^{pqj} [\theta,_{p\ao\at\att}^m \theta,_{j\ao\at\att}^i]_t \theta,_q^n \\
&\ \ +\frac{1}{2}\sum_{\sigma\in\Sigma_3} \epsilon^{mni}\epsilon^{pqj} \theta,_{p\alpha_{\sigma(1)}}^m\theta,_{q\alpha_{\sigma(2)}\alpha_{\sigma(3)}}^n  W,_{j\ao\at\att}^i\\
&\ \ +\frac{1}{2}\sum_{\sigma\in\Sigma_3} \epsilon^{mni}\epsilon^{pqj} \theta,_{p\alpha_{\sigma(1)}\alpha_{\sigma(2)}}^m\theta,_{q\alpha_{\sigma(3)}}^n  W,_{j\ao\at\att}^i,
\end{align*}
where we have use $\epsilon^{mni}\epsilon^{pqj}=\epsilon^{nmi}\epsilon^{qpj}$ on the second equality and $\epsilon^{mni}\epsilon^{pqj}=\epsilon^{inm}\epsilon^{jqp}$ on the third one. Thus,
\begin{align*}
Q_2&=\int_0^t\int_{\R^3_-} Q \tilde b_i^j,_{\ao\at\att}[[\z^2 W]^i,_{\ao\at\att j}-\z^2 W,_{j\ao\at\att}^i]\\
&\ \ -\frac{1}{2}\int_0^t\int_{\R^3_-} \sum_{\sigma\in\Sigma_3} \epsilon^{mni}\epsilon^{pqj} [ Q \theta,_{p\alpha_{\sigma(1)}}^m\theta,_{q\alpha_{\sigma(2)}\alpha_{\sigma(3)}}^n \z^2],_{\ao} W,_{j\at\att}^i\\
&\ \ -\frac{1}{2}\int_0^t\int_{\R^3_-} \sum_{\sigma\in\Sigma_3} \epsilon^{mni}\epsilon^{pqj} [Q\theta,_{p\alpha_{\sigma(1)}\alpha_{\sigma(2)}}^m\theta,_{q\alpha_{\sigma(3)}}^n\z^2],_{\ao}  W,_{j\at\att}^i\\
&\ \ -\frac{1}{2} \int_0^t\int_{\R^3_-}\epsilon^{mni}\epsilon^{pqj}\z^2 \theta,_{p\ao\at\att}^m \theta,_{j\ao\at\att}^i [Q \theta,_q^n]_t\\
&\ \ +\frac{1}{2}\bigl[\int_{\R^3_-}\epsilon^{mni}\epsilon^{pqj}\z^2 \theta,_{p\ao\at\att}^m \theta,_{j\ao\at\att}^i Q \theta,_q^n\bigr]_0^t,
\end{align*}
showing 
\begin{align*}
|Q_2|&\le C\sqrt{t} \|(\w,\q)\|_{Z_t}^4 +\bigl|\int_{\R^3_-}\epsilon^{mni}\epsilon^{pqj}\z^2 \theta,_{p\ao\at\att}^m \theta,_{q\ao\at\att}^n Q \theta,_j^i\bigr|(t)
+N((q_i)_{i=0}^2).
\end{align*}
In order to estimate the remaining term, we notice by integrating by parts twice for $\d Q_5=\int_{\R^3_-}\epsilon^{mni}\epsilon^{pqj}\z^2 \theta,_{p\ao\at\att}^m \theta,_{q\ao\at\att}^n Q \theta,_j^i$ that
\begin{align*}
Q_5&=\epsilon^{mni}\epsilon^{pqj}\int_{\R^3_-} \z^2 \theta,_{q\ao\at\att}^m \theta,_{p\ao\at\att}^n Q\ \theta,_j^i\\ 
&-\epsilon^{mni}\epsilon^{pqj}\int_{\R^3_-} [(\z^2 Q \theta,_j^i),_p \theta,_{\ao\at\att}^m \theta,_{q\ao\at\att}^n-(\z^2 Q \theta,_j^i),_q \theta,_{\ao\at\att}^m \theta,_{p\ao\at\att}^n] \\ 
&\ \ + \epsilon^{mni}\epsilon^{pqj}\int_{x_3=0}\z^2 Q\  \theta,_{\ao\at\att}^m [\theta,_{q\ao\at\att}^n \theta,_j^i (e_3)_p -  \theta,_{p\ao\at\att}^n  \theta,_j^i (e_3)_q].
\end{align*}
Since $\epsilon^{pqj}=-\epsilon^{qpj}$, we then infer
\begin{align*}
2 Q_5=&-\epsilon^{mni}\epsilon^{pqj}\int_{\R^3_-} \theta,_{\ao\at\att}^m [(\z^2 Q \theta,_j^i),_p  \theta,_{q\ao\at\att}^n-(\z^2 Q \theta,_j^i),_q \theta,_{p\ao\at\att}^n] \\  
&+\epsilon^{mni}\epsilon^{pqj}\int_{x_3=0}\z^2 Q\  \theta,_{\ao\at\att}^m [\theta,_{q\ao\at\att}^n \theta,_j^i (e_3)_p -  \theta,_{p\ao\at\att}^n \theta,_j^i (e_3)_q].
\end{align*}
Now, if we note $\theta^f=E(\of)(\e^f)\circ\P$, we also have for 
$$\d \int_{\R^3_-}\epsilon^{mni}\epsilon^{pqj}\z^2 \theta^f,_{p\ao\at\att}^m \theta,_{q\ao\at\att}^n Q \theta,_j^i$$ a similar formula. Since $\theta,_{\ao\at\att}^m=\theta^f,_{\ao\at\att}^m$ on $\{x_3=0\}$, we then have
\begin{align*}
2 Q_5=&-\epsilon^{mni}\epsilon^{pqj}\int_{\R^3_-} (\theta-\theta^f),_{\ao\at\att}^m [(\z^2 Q \theta,_j^i),_p  \theta,_{q\ao\at\att}^n-(\z^2 Q \theta,_j^i),_q \theta,_{p\ao\at\att}^n] \\
&+ 2\int_{\R^3_-}\epsilon^{mni}\epsilon^{pqj}\z^2 \theta^f,_{p\ao\at\att}^m \theta,_{q\ao\at\att}^n Q \theta,_j^i,
\end{align*}
leading us to 
\begin{align*}
|Q_5|&\le C \|\text{Id}+\int_0^t \w\|^{\frac{1}{4}}_{H^3(\os;\R^3)} \|\e\|^{\frac{7}{4}}_{H^4(\os;\R^3)}
\|q_0+\int_0^t \q_t\|_{H^2(\os;\R)} \|\text{Id}+\int_0^t \w\|_{H^3(\os;\R^3)}\\
&\  +C \|\text{Id}+\int_0^t \w\|^{\frac{1}{4}}_{H^3(\of;\R^3)} \|\e\|^{\frac{7}{4}}_{H^4(\os;\R^3)}
\|q_0+\int_0^t \q_t\|_{H^2(\os;\R)} \|\text{Id}+\int_0^t \w\|_{H^3(\os;\R^3)}\\
&\  +C \|\text{Id}+\int_0^t \w^f\|_{H^4(\of;\R^3)} \|\e\|_{H^4(\os;\R^3)}
\|q_0+\int_0^t \q_t\|_{H^2(\os;\R)} \|\text{Id}+\int_0^t \w\|_{H^3(\os;\R^3)}\\
&\le \delta \|(\w,\q)\|^2_{Z_t} + C\sqrt{t} \|(\w,\q)\|^4_{Z_t} + C_{\delta} [N(u_0,(w_i)_{i=1}^3)+ N((q_i)_{i=0}^2)]
\end{align*}
\noindent {Step 2.} 
We see by integrating by parts with respect to the direction $\alpha_1$ that we have 
\begin{align*}
|Q_3|\le C\sqrt{t} \|(\w,\q)\|^4_{Z_t}.
\end{align*}

\noindent Step 3. Next, $Q_{4}=Q_{6}+Q_{7}$, where
\begin{align*}
Q_6&=\int_0^{t}\int_{\R^3_-} {Q},_{\ao\at\att}\ \tilde b_i^j \z^2 W^i,_{\ao\at\att j}\\
Q_7&=\int_0^{t}\int_{\R^3_-} {Q},_{\ao\at\att}\ \tilde b_i^j [ [\z^2 W]^i,_{\ao\at\att j}-\z^2 W^i,_{\ao\at\att j}].
\end{align*}
We first have
\begin{align*}
|Q_7|\le C \int_0^t [\|\q\|_{H^3(\os;\R)} \|\w\|_{H^3(\os;\R^3)}\|\e\|^2_{H^3(\os;\R^3)}]\le  C \sqrt{t} \|(\w,\q)\|^4_{Z_t}\ .
\end{align*} 
For $Q_6$ the divergence condition $\tilde b_i^j W^i,_j=0$ on $Supp\ \z$ implies 
\begin{align*}
Q_6&=\int_0^{t}\int_{\R^3_-} {Q},_{\ao\at\att} \z^2\ [\tilde b_i^jW^i,_{\ao\at\att j}-(\tilde b_i^j W^i,_j),_{\ao\at\att}],
\end{align*}
which in turn provides,
\begin{align*}
|Q_6|&\le  C \int_0^t [\|\q\|_{H^3(\os;\R)}\|\w\|_{H^3(\os;\R^3)}\|\e\|^2_{H^4(\os;\R^3)}]\le C \sqrt{t} \|(\w,\q)\|^4_{Z_t}\ ,
\end{align*}
which concludes the estimates on the pressure terms in the solid, justifying why we obtain a trace estimate similar as (\ref{tracew}) with a majorant of the type of the right-hand side of (\ref{inc11}). Now, we turn our attention to the recovery of the regularity
in the solid, which will need some justifications since we cannot directly apply Lemma \ref{key}.
\subsubsection{\bf Regularity in the incompressible solid}
First, with the introduction of
\begin{equation*}
\tilde F=-\w_{t}+c^{ijkl}[(\e,_i\cdot\e,_j-\delta_{ij})\e,_k],_l-L\e +f+\kappa h\ \text{in}\ \os.
\end{equation*}
and of $\tilde r$, the solution in $\os$ of
\begin{align*}
\frac{\kappa}{2} \tilde r_t +\tilde r&=\tilde q,\\
\tilde r(0)&=q_0,
\end{align*}
we have for the nonlinear elastodynamics
\begin{equation*}
-\frac{\kappa}{2} \nabla\e\ L\w  -\nabla\e\ L\e + \nabla[\frac{\kappa}{2}\tilde r_t+ \tilde r] =\nabla\e\ \tilde F\ \text{in}\ \os,
\end{equation*}
{\it i.e.},
\begin{equation}
\label{inc12}
\frac{\kappa}{2} [-\nabla\e\ L\e + \nabla\tilde r]_t -\nabla\e\ L\e + \nabla\tilde r =\nabla\e\ \tilde F -\frac{\kappa}{2}\nabla\w\ L\e\ \text{in}\ \os.
\end{equation}
We now apply Lemma \ref{key} to this equation, leading us to
\begin{align*}
\sup_{[0,t]}\|-\nabla\e\ L(\e)+\nabla\tilde r\|_{H^2(\os;\R^3)}\le&\  \sup_{[0,t]}\|\nabla\e\ \tilde F -\frac{\kappa}{2}\nabla\w\ L\e\|_{H^2(\os;\R^3)}\\
& +\|-L(\text{Id})+\nabla\tilde q_0\|_{H^2(\os;\R^3)},
\end{align*}
and, with $\tilde H= - L(\e)+\nabla\tilde r$, to
\begin{align}
&\sup_{[0,t]}\|\tilde H\|_{H^2(\os;\R^3)}\n\\
&\le  \sup_{[0,t]} [\|\nabla\e\ \tilde F -\frac{\kappa}{2}\nabla\w\ L\e\|_{H^2(\os;\R^3)}+
\|(\nabla\e-\text{Id})\ L(\e) \|_{H^2(\os;\R^3)}] +N((q_i)_{i=0}^2).
\label{inc13}
\end{align}
We then want to use elliptic regularity on the system:
\begin{subequations}
\label{inc14}
\begin{align}
 -L\e + \nabla \tilde r &=\tilde H,\ \text{in}\ \os,\\
\operatorname{div}\e&=(-\a_i^j+\delta_{ij})\e^i,_j +3,\ \text{in}\ \os,\\
\e&=\e_{|\Gamma_0}\ \text {on}\ \Gamma_0,
\end{align}
\end{subequations}
where the trace on $\Gamma_0$ is estimated as we explained in the previous subsection. Now, for the divergence condition in $\os$, we notice that:
\begin{align*}
[(\a_i^j-\delta_{ij})\e^i,_j],_{i_1i_2i_3}=&\ \a_i^j,_{i_1i_2i_3}\e^i,_j +(\a_i^j-\delta_{ij})\e^i,_{ji_1i_2i_3}\\
&+ \sum_{\sigma\in\Sigma_3} [\a_i^j,_{i_{\sigma(1)}}\e^i,_{ji_{\sigma(2)}i_{\sigma(3)}}+ \a_i^j,_{i_{\sigma(1)}i_{\sigma(2)}}\e^i,_{ji_{\sigma(3)}}].
\end{align*}
For the apparently problematic first term on the right-hand side, we first
notice that
\begin{align*}
\a_i^j,_{i_1}\e^i,_j&=\frac{1}{2}\epsilon^{mni}\epsilon^{pqj} (\e,_p^m\e,_q^n),_{i_1} \e,_j^i=\epsilon^{mni}\epsilon^{pqj} \e,_{pi_1}^m\e,_q^n \e,_j^i\\
&=\epsilon^{inm}\epsilon^{jqp} \e,_{ji_1}^i\e,_q^n \e,_p^m=2 \a_i^j\e^i,_{ji_1},
\end{align*}
which with the condition $\a_i^j\e^i,_j=3$, provides
$ 
0= \a_i^j,_{i_1}\e^i,_{j},
$ and thus
\begin{align*}
\a_i^j,_{i_1i_2i_3}\e^i,_j=-\a_i^j,_{i_1i_2}\e^i,_{ji_3}- \a_i^j,_{i_1i_3}\e^i,_{ji_2}-\a_i^j,_{i_1}\e^i,_{ji_2i_3}.
\end{align*}
We then deduce that
\begin{align}
\|(\a_i^j-\delta_{ij})\e^i,_j\|_{H^3(\os;\R)}(t)\le \delta \|(\w,\q)\|^2_{Z_t} + C_{\delta} N(u_0,(w_i)_{i=1}^3) + C t\ \|(\w,\q)\|^3_{Z_t}.
\label{inc15}
\end{align}
Now, with (\ref{inc13}) and (\ref{inc15}), elliptic regularity on (\ref{inc14})
provides for $$\|\e\|^2_{L^\infty(0,t;H^4(\os;\R^3))}+\|\tilde r-\frac{1}{|\os|}\int_{\os} \tilde r\|^2_{L^\infty(0,t;H^3(\os;\R))}$$ a bound of the same type as the right-hand side of (\ref{inc11}), with however the norms in $Z_t$ being replaced by the norms in $\tilde Z_t$, due to the term $\kappa \|\nabla\w\ L\e\|_{H^2(\os;\R^3)}$ appearing on the right-hand side of (\ref{inc13}), that we bound by 
\begin{align*}
C\kappa \|\w(t)\|_{H^3(\os;\R^3)} \|L(\text{Id})+\int_0^tL\w\|_{H^2(\os;\R^3)}&\le C \kappa \|\w(t)\|_{H^3(\os;\R^3)} \sqrt{t} \|\w\|_{L^2(0,t;H^4(\os;\R^3))}
\\
&\le C \sqrt{t}\|(\w,\q)\|^2_{\tilde Z_t}.
\end{align*}

We now turn our attention to the pressure, that we just need to control in $L^2(0,t;H^3(\os;\R))$. In order to do so, we notice from (\ref{inc12}) that we have for $\tilde K= \frac{\kappa}{2}[-L\w+\nabla\tilde r_t]$:
\begin{align*}
\tilde K=
\frac{\kappa}{2} [\nabla\e-\text{Id}]\ L\w +\nabla\e\ L\e-\nabla\tilde r+\nabla\e\ \tilde F-\frac{\kappa}{2}\nabla\w\ L\e,
\end{align*}
which with the previous estimate on $\e$ and $\tilde r$ shows that we have a bound on
$\|\tilde K\|^2_{L^2(0,t;H^2(\os;\R^3))}$ of the same type as the right-hand side of (\ref{inc11}), but where the norms in $Z_t$ are replaced by norms in $\tilde Z_t$ due to the estimate in $L^2(H^2)$ of $\kappa [\nabla\e-\text{Id}]\ L\w $.
Now, elliptic regularity on the system:
\begin{align*}
 -L\kappa\w + \nabla \kappa r_t &=2\tilde K,\ \text{in}\ \os,\\
\operatorname{div}\kappa\w&=\kappa (-\a_i^j+\delta_{ij})\w^i,_j,\ \text{in}\ \os,\\
\kappa\w&=\kappa\w^f_{|\Gamma_0}\ \text {on}\ \Gamma_0,
\end{align*}
provides after integrating in time an estimate for $$\d \kappa^2 [\|\w\|^2_{L^2(0,t;H^4(\os;\R^3))}+\|\tilde r_t-\frac{1}{|\os|}\int_{\os}\tilde r_t\|^2_{L^2(0,t;H^3(\os;\R))}]$$ with a bound similar as in (\ref{inc11}), still with the norms in $Z_t$ being replaced by norms in $\tilde Z_t$.

Thus, we obtain for $\|\tilde q-\frac{1}{|\os|}\int_{\os}\q\|^2_{L^2(0,t;H^3(\os;\R))}$ the same type of estimate as well. Given our estimate on $\q_{tt}$, this also implies
the same type of majoration for $\|\tilde q\|^2_{L^2(0,t;H^3(\os;\R))}$. 

Thus, we are lead to
\begin{align*}
\|(\w,\q)\|^2_{\tilde Z_t}
\le & (C \kappa + C_{\delta} {t}^{\frac{1}{4}}) \|(\w,\q)\|^8_{\tilde Z_t}+\delta \|(\w,\q)\|^2_{\tilde Z_t} \n\\
& +C_\delta  [N(u_0,(w_i)_{i=1}^3)+N((q_i)_{i=0}^2)+M(f,\kappa g,\kappa h)],
\end{align*}
which leads as in Section \ref{9} to the introduction of a polynomial, this time of degree 4, which does not bring any substantial change with respect to Section \ref{9}. Note that the addition of $C\kappa \|(\w,\q)\|^8_{Z_t}$ does not create any
difficulty since
a small  $\kappa_1$ is chosen at the same stage as $t_1$, and the conclusion is 
similar as in Section \ref{9} from the continuity of $\|(\w,\q)\|_{\tilde Z_t}$ on $[0,T_{\kappa}]$ which is established in the same way as the continuity of
$\|(\w,\q)\|_{Z_t}$. 
We then infer that there is a time of existence of $\kappa$ for our smoothed problems, with a bound on $\|(\w,\q)\|_{\tilde Z_T}$ and thus on $\|(\w,\q)\|_{Z_T}$ independent of $\kappa$. Existence follows then by weak convergence in $Y_T$ and uniqueness can be established similarly as for the compressible case in Section \ref{11}.

\end{proof}

\section*{Acknowledgments}
DC was partially supported by the National Science Foundation under grant 
NSF ITR-0313370.  SS was partially supported by the National Science Foundation 
under grants 
DMS-0105004 and NSF ITR-0313370.

\end{document}